


\input amstex


\def \BDV{BDV}

\def \BSOne{BS1}   
\def \BSTwo{BS2}   
\def \BSThree{BS3}      
\def \Bou{Bou}   
\def \BrS{BrS}
\def \CEG{CEG}
\def \CG{CGv}    
\def \CGv{CGv}  
\def \CGl{CGl}  
\def \CHK{CHK}  
\def \Cul{Cul}  
\def \CS{CS}    
\def \DaM{DaM}   
\def \DD{DD}    
\def \Dun{Dun}  
\def \DuM{DuM}    
\def \Ebe{Ebe}   
\def \Fed{Fed}
\def \GM{GM}   

\def \Gro{Gro}  
\def \GLP{GLP}
\def \GT{GT}   
\def \HamOne{Ha1}
\def \HamTwo{Ha2}
\def \Hod{Ho1}  
\def \HodTwo{Ho2} 
\def \HK{HK}
\def \HT{HT}    
\def \Iva{Iva} 
\def \JR{JR}    
\def \JS{JS}      
\def \Joh{Joh}   
\def \Jon{Jon}   
\def \Kap{Kap}    
\def \Kir{Kir}    
\def \McMOne{McM1}
\def \McMTwo{McM2}
\def \KoOne{Ko1}  
\def \KoTwo{Ko2}  
\def \McM{McM}
\def \MYOne{MY1}    
\def \MYTwo{MY2}    
\def \MS{MS}
\def \MB{MB}
\def \Mos{Mos}
\def \Mum{Mum}
\def \OtaOne{Ot1}
\def \OtaTwo{Ot2}
\def \Pe{Pe}   
\def \PoOne{Po1}
\def \PoTwo{Po2}
\def \Rag{Rag}
\def \Sak{Sak}
\def \Scott{Sc} 
\def \SOK{SOK} 
\def \Som{Som} 
\def \Sua{Sua}
\def \ToOne{To1}
\def \ToTwo{To2}
\def \Thu{Th}   
\def \ThuNotes{Th1} 
\def \ThuBull{Th2}   
\def \Wal{Wa1}  
\def \WalOne{Wa1}  
\def \WalTwo{Wa2}  
\def \Wei{Wei}   
\def \Yam{Yam}  
\def \Zhou{Zh1}
\def \ZhouTwo{Zh2}

\def\ab#1{\operatorname{ab}({#1})}
\def\norm#1{\Vert{#1}\Vert}
\def\Star#1{\operatorname{star}({#1})}
\def\Stard#1{\operatorname{star}^*({#1})}
\def\Mod#1{\vert{#1}\vert}

\define\LL{{\Cal L}}
\define\N{{\Bbb N}}
\define\NN{{\Cal N}}
\define\CC{{\Cal C}}
\define\OO{{\Cal O}}

\define\Diam{\operatorname{diam}}
\define\Hom{\operatorname{Hom}}
\define\Inj{\operatorname{inj}}
\define\Int{\operatorname{int}}
\define\Trace{\operatorname{trace}}
\define\V{\operatorname{v}}
\define\Vol{\operatorname{vol}}
\define\Length{\operatorname{Length}}

\define\subp#1{\ \par\noindent{\bf{#1.}}}
\define\sub#1{\noindent{\bf{#1.}}}

\define\iM#1{\operatorname{inj}_M(#1)}

\define\ttag#1{$(\overset{\tsize{#1}}\to{\phantom{.}})$}

\documentstyle{amsppt}

\input BoxedEPS.tex
\SetRokickiEPSFSpecial
\HideDisplacementBoxes

\exhyphenpenalty=10000
\nologo

\font\chapt=cmbx10 scaled \magstep0
\font\tit=cmbx10 scaled \magstep2

\document

\rightheadtext{ Introduction}
\leftheadtext{ Introduction}

\topmatter
\title \tit Geometrization of 3-orbifolds of cyclic type.
\endtitle 
\date May 15, 1998 \enddate
\author M. Boileau and J. Porti
\endauthor
\endtopmatter

\centerline{\chapt INTRODUCTION }

\

 A 3-dimensional orbifold is a metrizable space with coherent local
models given by quotients of $\Bbb R^3$ by finite subgroups of
$O(3)$. For example, the quotient of a 3-manifold by a properly
discontinuous group action naturally inherits a structure of a
3-orbifold. Such an orbifold is said to be \it very good \rm when the group
action is finite. For a general background about orbifolds see
\cite{\BSOne,2}, \cite{\DaM}, \cite{\Scott} and \cite{\ThuNotes, 
Ch. 13}.

The purpose of this article is to give a complete proof of Thurston's
Orbifold Theorem in the case {\it where all local isotropy groups are cyclic
subgroups of $SO(3)$}. Following \cite{\DaM}, we say that such an
orbifold is of {\it cyclic type} when in addition the ramification locus
is non-empty. Hence a 3-orbifold $\OO$ is of cyclic type if and only if its
ramification locus $\Sigma_{\OO}$ is a non-empty 1-dimensional
submanifold of the underlying manifold $\Mod\OO$, which is transverse
to the boundary $\partial\Mod\OO=\Mod{\partial\OO}$. The version of
Thurston's Orbifold Theorem proved here is the following:

\proclaim{Theorem 1 \rm (Thurston's Orbifold Theorem)} Let $\OO$ be a 
compact  connected  orientable 
irreducible 
$\partial$-incom\-pressible 3-orbifold of cyclic type. If $\OO$ is
very good, topologically atoroidal and
 acylindrical, then $\OO$ is geometric (i.e. $\OO$ admits 
 either a hyperbolic, a Euclidean, or a Seifert fibered structure).
\endproclaim

\remark{Remark}\roster\runinitem When $\partial\OO$ is a union of
toric 2-suborbifolds, the hypothesis that $\OO$ is acylindrical is
not needed.
\item If $\partial\OO\neq\emptyset$ and $\OO$ is not
$I$-fibered, then $\OO$ admits a  hyperbolic structure with finite
volume. \endroster
\endremark

We only consider {\it smooth} orbifolds, so that the local isotropy
groups are always orthogonal. We recall that an orbifold is said to be
{\it good} if it has a covering which is a manifold. 
Moreover if this covering is finite then the orbifold is said to be 
{\it very good}.

According to \cite{\BSOne,2} and \cite{\ThuNotes, Ch.13}, we
use the following terminology. 

\definition{Definitions}
We say that a compact 2-orbifold $F^2$ is respectively
{\it spherical}, {\it discal}, {\it toric} or {\it annular}
if it is the quotient by a finite smooth group action of respectivly
the 2-sphere $S^2$, the 2-disk $D^2$, the 2-torus $T^2$ or the annulus
$S^1\times [0,1]$.

A compact 2-orbifold is {\it bad} if it is not good. Such a 
2-orbifold is the union of two non-isomorphic discal 2-orbifolds
along their boundaries.

A compact 3-orbifold $\OO$ is {\it irreducible} if it does not
contain any bad 2-suborbifold and every spherical 2-suborbifold bounds
in $\OO$ a discal 3-suborbifold, where a {\it discal} 3-orbifold is a
finite quotient of the 3-ball by an orthogonal action.

A connected 2-suborbifold $F^2$ in an orientable 3-orbifold $\OO$ is
{\it compressible} if either $F^2$ bounds a discal 3-suborbifold in
$\OO$ or there is a discal 2-suborbifold $\Delta^2$ which intersects
transversally $F$ in $\partial\Delta^2=\Delta^2\cap F^2$ and such
that $\partial\Delta^2$ does not bound a discal 2-suborbifold in
$F^2$.

A 2-suborbifold $F^2$ in an orientable 3-orbifold $\OO$ is
{\it incompressible} if no connected component of $F^2$ is
compressible in $\OO$.
 The compact 3-orbifold $\OO$ is {\it $\partial$- incompressible} if
$\partial\OO$ is empty or incompressible in $\OO$.

A compact 3-orbifold is {\it topologically atoroidal} if every
incompressible toric 2- suborbifold is boundary parallel (i.e. is the
frontier of a collar neighborhood $F^2\times [0,1]\subset\OO$ of a
boundary component $F^2\subset\partial\OO$).
A compact 3-orbifold is {\it topologically acylindrical} if every
properly embedded annular 2-suborbifold is boundary parallel.

A {\it Seifert fibration} on a 3-orbifold $\OO$ is a partition of 
$\OO$ into closed 1-suborbifolds (circles or intervals with silvered
boundary) called fibers, such that each fiber has a saturated
neighborhood orbifold-diffeomorphic to $S^1\times D^2/G$, where $G$
is a finite group which acts smoothly, preserves both factors, and
acts orthogonally on each factor and effectively on $D^2$; moreover
the fibers of the saturated neighborhood correspond to the quotients
of the circles $S^1\times\{*\}$. On the boundary $\partial\OO$, the
local model of the Seifert fibration is $S^1\times D^2_+/G$, where
$D^2_+$ is a half disk.

A 3-orbifold that admits a Seifert fibration is called Seifert
fibered. Every good Seifert fibered 3-orbifold is geometric. 
Seifert fibered 3-orbifolds have been classified in \cite{\BSTwo}.

A compact orientable 3-orbifold $\OO$ is {\it hyperbolic} if its
interior is orbifold-diffeo\-morphic to the quotient of the hyperbolic
space $\Bbb H^3$ by a non-elementary discrete group of isometries. In
particular $I$-bundles over hyperbolic 2-orbifolds are hyperbolic,
since their interiors are quotients of $\Bbb H^3$ by non-elementary
Fuchsian groups. In Theorem 1, except for $I$-bundles, we prove that
when $\OO$ is hyperbolic, if we remove the toric components of the
boundary $\partial_T\OO\subset\partial\OO$, then
$\OO-\partial_T\OO$ 
has a hyperbolic structure with finite volume and
geodesic boundary. This implies the existence of a complete hyperbolic structure
on the interior of $\OO$.

We say that a compact orientable 3-orbifold is {\it Euclidean} if its
interior has a complete Euclidean structure. Thus, if a compact
orientable and $\partial$-incompressible 3-orbifold $\OO$ is Euclidean,
then  either $\OO$ is a  $I$-bundle over a 2-dimensional Euclidean
closed orbifold or $\OO$ is closed.

We say that a compact orientable 3-orbifold is {\it spherical} when
it is the quotient of $\Bbb S^3$ by the orthogonal action of 
a finite subgroup of $SO(4)$. A
spherical orbifold of cyclic type is always Seifert fibered
(\cite{Dun}, \cite{\DaM} and \cite{\Zhou}).
\enddefinition

 Thurston's conjecture asserts that the interior of a compact
irreducible orientable 3-orbifold can be decomposed along a canonical
family of incompressible toric 2-suborbifolds into geometric
3-suborbifolds.

The existence of the canonical family of incompressible toric
2-suborbifolds has been established by Jaco-Shalen \cite{\JS} and 
Johannson \cite{\Joh} for 3-manifolds and by Bonahon-Siebenmann
\cite{\BSTwo} in the case of 3-orbifolds.
 
Recall that the eight 3-dimensional geometries involved in Thurston's
conjecture are $\Bbb H^3$, $\Bbb E^3$, $\Bbb S^3$, $\Bbb
H^2\times\Bbb R$, $\Bbb S^2\times \Bbb R$, $\widetilde{SL_2(\Bbb R)}$,
$Nil$ and $Sol$. The non Seifert fibered orbifolds require a constant
curvature geometry ($\Bbb H^3$, $\Bbb E^3$ and $\Bbb S^3$) or $Sol$.
Compact orbifolds with $Sol$ geometry are fibered over a closed 
1-dimensional orbifold with toric fiber and thus are not atoroidal.

Thurston \cite{\ThuNotes,2,3,4,5} has proved his conjecture for Haken
3-manifolds (cf. \cite{\McMOne,2, \Kap, \OtaOne,2}). In 1981, Thurston
\cite{Thu2,6} announced the Geometrization Theorem for 3-orbifolds
with non-empty ramification set (without the assumptions very good and
of cyclic type), and lectured about it. Since 1986, several
useful notes about Thurston's proof (by Soma, Oshika and
Kojima \cite{\SOK} and by Hodgson
\cite{Hod}) have been circulating. In addition, in 1989 more details appeared in Zhou's thesis
\cite{\Zhou,2} in the cyclic case. However no complete proof has
appeared yet (cf.
\cite{\Kir, Prob. 3.46}).

The following corollary is a straightforward application of Meeks and
Scott's work \cite{MS} and of Theorem 1 to the geometrization of
3-manifolds with non-trivial symmetries. We state it only in the
case of 3-manifolds with empty or toric boundary for simplicity.

\proclaim{Corollary 1} Let $M$ be a compact orientable irreducible
3-manifold with zero Euler characteristic. Let $G$ be a finite group
of orientation preserving diffeomorphisms acting on $M$ with
non-trivial and cyclic stabilizers. Then there exists a (possibly empty)
$G$-invariant family  of disjoint incompressible tori which splits
$M$ into $G$-invariant geometric pieces.
\endproclaim

As a particular case, for finite group actions on geometric
3-manifolds, together with \cite{\MS} we obtain:

\proclaim{Corollary 2} Let $M$ be a closed orientable irreducible
3-manifold with a geometric structure. Every finite orientation
preserving smooth action with non-trivial cyclic stabilizers is
conjugate to a geometric action.
\endproclaim

With respect to the study of conjugacy classes of finite subgroups
of $\operatorname{Diff}^+(M)$ (cf. \cite{\Kir, Prob.  3.39}), the
following corollary is a consequence of Kojima's work \cite{\KoTwo}
and Theorem 1.

\proclaim{Corollary 3} Let $M$ be an orientable compact irreducible
3-manifold. Let $G\subset\operatorname{Diff}^+(M)$ be a finite 
subgroup with non-trivial cyclic stabilizers. Then the number of
conjugacy classes of $G$ in $\operatorname{Diff}^+(M)$ is finite.
\endproclaim

In the case of cyclic branched coverings of links in $S^3$, the
following results answers a question of Montesinos \cite{\Kir, Prob.
3.41}.

\proclaim{Corollary 4} For a given integer $p\geq 2$, there are only
finitely many links in $S^3$ with the same $p$-fold cyclic branched
covering.
\endproclaim

Using Thurston's Hyperbolization Theorem
\cite{\ThuNotes,2,3,4,5, \McMOne,2, \Kap, \OtaOne,2} and a standard
argument of doubling $\OO$ along  boundary components, the proof
of Theorem 1 reduces to the following theorem, which is the main
result of this paper.

\pagebreak

\proclaim{Theorem 2} Let $\OO$ be a closed orientable 
connected irreducible very good  3-orbifold of
cyclic type. Assume that the complement $\OO-\Sigma$ 
of the branching locus admits a complete
hyperbolic structure. Then there exists a non-empty compact essential 
$3$-suborbifold $\OO'\subseteq\OO$ which is not a
product and which is either Euclidean, Seifert fibered, $Sol$ or complete hyperbolic
with finite volume. In particular $\partial\OO'$ is either empty or a
union of toric 2-orbifolds.
\endproclaim

A compact $3$-suborbifold $\OO'\subseteq\OO$  is {\it essential} in  $\OO$ if the 2-suborbifold
$\partial\OO'$ is either empty or incompressible in $\OO$.

\remark{Remarks}
\roster
\item It follows from the proof of this theorem that the hypothesis
that $\OO$ is very good is not needed when all local isotropy groups
have order at least 4 (or 3 and $\OO$ does not contain the
toric orbifold $S^2(3,3,3)$). The orbifold
$\OO$ is seen to be very good a posteriori as a corollary of Theorem 2,
because it is geometric.
\item If the orbifold $\OO$ is topologically atoroidal, then
$\OO=\OO'$ is geometric. This is the case when all branching
indices are either  at least 4, or at least 3 and in addition $\OO$ does not contain the
toric orbifold $S^2(3,3,3)$. If the underlying manifold $\Mod\OO$ is
irreducible then $\OO$ does not contain  $S^2(3,3,3)$.
\endroster
\endremark

A straightforward corollary of the proof of Theorem 2 and of the
classification of the orientable closed Euclidean 
and spherical
3-orbifolds (cf. 
\cite{\BSOne},  \cite{\Dun}) is the following strong version of Smith
conjecture for a link in an irreducible 3-manifold.

\proclaim{Corollary 5}
Let $M$ be a compact orientable irreducible 3-manifold and $L\subset
M$ be a hyperbolic link. Then, for $p\geq 3$, any $p$-fold cyclic
covering of $M$ branched along $L$ admits a hyperbolic structure,
except when
$p=3$, $M=S^3$ and $L$ is the figure- eight knot, because in this case
the 3-fold cyclic branched covering is the 3-torus. Moreover, in all cases the
covering transformation group acts by isometries.
\endproclaim

The proof of Theorem 2 follows Thurston's original approach. His idea
was to deform the complete hyperbolic structure as far as possible
on $\OO-\Sigma$ into structures whose completion is topologically the
underlying manifold $\Mod\OO$ and has cone singularities along
$\Sigma$. These completions are called hyperbolic cone structures on 
the pair $(\Mod\OO,\Sigma)$ and their singularities are (locally)
described by cone angles. Such structures with small cone angles
are provided by Thurston's  hyperbolic Dehn filling
theorem \cite{\ThuNotes, Ch. 5} and the goal is  to study the limit
of hyperbolicity when these cone angles increase. Note that a
hyperbolic structure on $\OO$ induces a hyperbolic cone structure on
the pair
$(\Mod\OO ,\Sigma)$ with cone angles determined by the ramification
indices. Hence, if these cone angles can be reached in the space of hyperbolic
cone manifold structures, then $\OO$ is hyperbolic. Otherwise, the study
of the possible ``accidents" occuring at the limit of hyperbolicity shows
the existence of a non-empty compact essential geometric 3-suborbifold
$\OO'\subset\OO$ which is different from $\OO$ when $\OO'$ is
hyperbolic.

Our main contribution takes place in the analysis of the so called 
``collapsing cases". There we use the notion of simplicial volume due
to Gromov and a cone manifold version of his Isolation Theorem
\cite{\Gro, Sec. 3}. This gives a simpler combinatorial approach to
collapses than Thurston's original one. In particular, it spares us
the difficult  task of establishing a suitable Cheeger-Gromov theory
for collapses of cone manifolds.

However there is a price to be paid: we must work, at least when some of the
ramification indices are 2, with very good 3-orbifolds, since we
use in a crucial way the results of Meeks and Scott  \cite{\MS}.

 Another approach, in the case of closed, irreducible, geometrically
atoroidal, orientable 3-orbifolds of cyclic type, more in the spirit of
Thurston's original approach, has been announced by D.
Cooper, C. Hodgson and S. Kerckhoff \cite{\CHK}. This approach does
not require the 3-orbifold to be very good.

\

Here  is a plan of the paper. 
In Chapter I we introduce the notion of cone manifold and state the theorems
that are the main ingredient in the proof of of Thurston's Orbifold Theorem.
In Chapter II we show how to deduce Theorems 1 and 2 from 
the results quoted in Chapter I.
The remaining chapters are devoted to the proof of the theorems stated
in Chapter I. More precisely, in Chapter III we prove the Compactness
Theorem, which is a cone 3-manifold version of Gromov's Compactness 
Theorem for Riemannian manifolds of pinched sectional curvature.
In Chapter IV we prove the Local Soul Theorem, which gives a bilipschitz
approximation of the metric structure of the neighborhood of points with 
small cone-injectivity radius. By using the Compactness and the Local Soul
Theorem, in Chapters V and VI we study sequences of cone 3-manifolds with
fixed topologycal type. In Chapter V we prove Theorem A, which deals 
with the case where the cone angles are bounded above, uniformly away from $\pi$.
In Chapter VI, we prove Theorem B, which deals with the case where the
cone angles converge to $2\pi/n_i$, for some $n_i\geq 2$.

\

We wish to thank Berhnard Leeb for his interest in our work and many useful discussions.
During the writing, we benefitted also from fruitful comments by Michael Heusener,
Greg McShane and Carlo Petronio.

\leftheadtext{ }
\newpage
\rightheadtext{I \qquad   cone manifolds}
\leftheadtext{I  \qquad  cone manifolds}

\

\centerline{\smc chapter \ i}

\

\centerline{\chapt CONE \,  MANIFOLDS}

\

\

Cone 3-manifolds play a central role in the proof of
Thurston's Orbifold Theorem. Thurston has shown that
they appear naturally as a generalization of
Hyperbolic Dehn Filling on cusped hyperbolic manifolds.
Thurston's  Hyperbolic Dehn Filling Theorem
provides a family of cone 3-manifolds with small cone
angles and the proof of Thurston's Orbifold Theorem
analyzes the accidents that can occur when we increase
the cone angles in order to reach the hyperbolic metric
on the orbifold.

This chapter has two sections. In the first  we
give the basic definitions for cone 3-manifolds (of 
non-positive curvature). In the second  we
state some theorems about sequences of hyperbolic cone 3-manifolds,
which are the key steps in the proof of Thurston's
orbifold Theorem.

\head 1. Basic Definitions
 \endhead

In this paper we only consider cone 3-manifolds of
non-positive constant curvature. Moreover, we also
restrict our attention to cone 3-manifolds whose singular
set is a link and whose cone angles are less than
$2\pi$.

To fix the notations, let $\Bbb H^3_K$ be the simply
connected three-dimensional space of constant sectional
curvature $K\leq 0$. Thus $\Bbb H^3_{-1}\cong\Bbb H^3$
is the usual hyperbolic space and $\Bbb H^3_0\cong\Bbb
E^3$ is the Euclidean space.

For $\alpha\in(0,2\pi)$, let $\Bbb H^3_K(\alpha)$ be the
cone manifold of constant curvature $K\leq 0$ with a
singular line of cone angle $\alpha$, constructed as
follows.  Consider in $\Bbb H^3_K$ a solid angular
sector $S_{\alpha}$ obtained by taking the intersection of two
half spaces, such that the dihedral angle at the
(infinite) edge $\Delta$ is $\alpha$. The cone manifold
$\Bbb H^3_K(\alpha)$ is the length space
obtained when we identify the faces  of $S_{\alpha}$ by
a rotation around $\Delta$. The image of $\Delta$ in
the quotient gives the singular line $\Sigma\subset
\Bbb H^3_K(\alpha)$. The induced metric on
$\Bbb H^3_K(\alpha)-\Sigma$ is a non-singular, 
incomplete Riemannian metric of constant curvature,
whose completion is precisely $\Bbb H^3_K(\alpha)$.

In cylindrical or Fermi coordinates, the metric on 
$\Bbb H^3_K(\alpha)-\Sigma$ is:
$$
d s^2_K=
\left\{
\matrix\format \l\,&\ \l\\
d r^2+
\left(\frac{\alpha}{2\pi}
  \frac{\sinh(\sqrt{-K} r)}{\sqrt{-K}}\right)^2
d\theta^2
+ \cosh^2(\sqrt{-K} r) dh^2 &\text{ for } K<0\\
dr^2+\left(\frac{\alpha}{2\pi} r\right)^2 d\theta^2+dh^2
&\text{ for } K=0
 \endmatrix
\right.
$$
where $r\in (0,+\infty)$ is the distance from $\Sigma$,
$\theta\in [0,2\pi)$ is the rescaled angle parameter
around $\Sigma$ and $h\in\Bbb R$ is the distance along $\Sigma$. 

Having described the local models, we can now define a
cone 3-manifold.

\definition{Definition}
A {\it cone manifold} of dimension three and of
constant curvature $K\leq 0$ is a smooth $3$-manifold
$C$ equipped with a distance so that it is a complete
length space locally isometric to
$\Bbb H^3_K$
or $\Bbb H^3_K(\alpha)$ for some $\alpha\in(0,2\pi)$.
\enddefinition

The {\it singular locus} $\Sigma\subset C$ is the set
of points modeled on the singular line of some model
$\Bbb H^3_K(\alpha)$, and $\alpha$ is called the {\it
cone angle} at a singular point modeled on this
singular line. According to our definition, $\Sigma$ is
a submanifold of codimension two and the cone angle is constant along each connected component.

The topological pair $(C,\Sigma)$ is called the {\it
topological type} of the cone 3-manifold.

The induced metric on $C-\Sigma$ is a Riemannian metric
of constant curvature, which is incomplete (unless
$\Sigma=\emptyset$), and whose completion is precisely
the cone 3-manifold.

By the {\it developping map of a cone 3-manifold $C$}
with topological type $(C,\Sigma)$, we mean the
developping map of the induced metric on $C-\Sigma$:
$$
  D:\widetilde{C-\Sigma}\to \Bbb H^3_K,
$$
where $\widetilde{C-\Sigma}$ is the universal covering
of $C-\Sigma$. The associated holonomy representation
$$
\rho:\pi_1(C-\Sigma)\to\operatorname{Isom}(\Bbb H^3_K)
$$
is called the {\it holonomy representation of} $C$. 
If $\mu\in\pi_1(C-\Sigma)$ is represented by a meridian
loop around a component $\Sigma_0$ of $\Sigma$, then
$\rho(\mu)$ is a rotation of angle equal to the cone
angle of this component.

Thurston's Hyperbolic Dehn Filling theorem provides
many structures on a hyperbolic cusped 3-manifold whose
completions are precisely  cone 3-manifolds. The cone
angles of these cone 3-manifolds are not necessarily less
than $2\pi$. The complete cusped structure on $M$ is the
limit of these hyperbolic cone structures when the cone
angles approach zero. We adopt therefore the standard convention
that the cone angle at a component $\Sigma_0$ of
$\Sigma$ is zero when the end of $C-\Sigma$
corresponding to $\Sigma_0$ is a cusp, of rank 2 or 1
according to whether $\Sigma_0$ is compact or not.

We still need two more definitions.

A {\it standard ball} in a cone 3-manifold $C$ is a 
ball isometric to either a metric non-singular 
ball in $\Bbb H_K^3$ or a metric singular ball
in $\Bbb H_K^3(\alpha)$ whose center belongs to
the singular axis.

We define the
{\it cone-injectivity
radius} at a point $x\in C$ as
$$
\Inj(x)=\sup\{\delta>0\text{ such that } B(x,\delta)\text{ is
contained in a standard ball in }C\}.
 $$
We remark that the standard ball does not need to be centered
at the point $x$. With this definition, regular points close
to the singular locus do not have arbitrarily small 
cone-injectivity radius.

\head 2. Sequences of hyperbolic cone 3-manifolds
\endhead

Let $\OO$ be an orbifold as in Theorem 2: a closed  orientable 
irreducible  very good 3-orbifold of cyclic type, such that the
complement
$\OO-\Sigma$ of the ramification locus admits a complete hyperbolic
structure of finite volume.

Thurston's  Hyperbolic Dehn Filling Theorem provides a one-parame\-ter
family of hyperbolic cone 3-manifolds with topological type
$(\Mod\OO,\Sigma)$. Consider the exterior $X=\OO-\Int(\NN(\Sigma))$
of $\Sigma$ and for each component $\Sigma_i\subset\Sigma$ choose a
meridian curve $\mu_i\subset\partial\NN(\Sigma_i)$ and another
simple closed curve $\lambda_i\subset\partial\NN(\Sigma_i)$
intersecting $\mu_i$ in one point; hence $\mu_i,\lambda_i$ generate
$\pi_1(\partial\NN(\Sigma_i))$.

According to Thurston's  Hyperbolic Dehn Filling Theorem, there
exists a space of deformations of hyperbolic structures on
$\Int(X)$ parametrized by generalized Dehn coefficients
$(p_i,q_i)$, $1\leq i\leq k=\operatorname{card} (\pi_0(\Sigma))$
in an open neighborhood $U\subset
(\Bbb R^2\cup\{\infty\})^k\cong(S^2)^k$ of $(\infty,\ldots,\infty)$,
such that the structure at the $i$-th component of $\partial X$ is
described by the Dehn parameters as follows:
\roster
\item"--" When $(p_i,q_i)=\infty$, the structure at the
corresponding cusp remains complete.
\item"--" When $p_i,q_i\in\Bbb Z$ are coprime, the
completion $X((p_1,q_1),\ldots,(p_k,q_k))$ is a hyperbolic 3-manifold,
obtained by genuine Dehn filling with meridian curves $p_i\mu_i +
q_i\lambda_i$, $i=1,\ldots,k$.
\item"--"  When $p_i/q_i\in\Bbb Q\cup\{\infty\}$, let
$r_i,s_i\in\Bbb Z$ be coprime integers so that $p_i/q_i=r_i/s_i$.
Then the completion $X((p_1,q_1),\ldots,(p_k,q_k))$ is a hyperbolic
cone 3-manifold obtained by gluing solid tori with possibly
singular cores. The underlying space is the 3-manifold
$X((r_1,s_1),\ldots,(r_k,s_k))$ and the cone angle of the $i$-th
singular core is $2\pi r_i/p_i$, $i=1,\ldots,k$.
\endroster

Here we are interested in the coefficients of the form
$(p_i,q_i)=(n_i/t,0)$, where $t\in[0,1]$ and $n_i$ is the branching
index of the orbifold $\OO$ along the $i$-th component $\Sigma_i\subset\Sigma$, for
$i=1,\ldots,k$. Thurston's  Hyperbolic Dehn Filling Theorem implies
the existence of a real number $\varepsilon_0>0$ such that, for any
$t\in[0,\varepsilon_0]$, there is a deformation of the complete
hyperbolic structure on $\Int(X)$  whose completion 
$X\left((\frac{n_1}t,0),\ldots,(\frac{n_k}t,0)\right)$ is a
hyperbolic cone 3-manifold with topological type $(\Mod\OO,\Sigma)$
and cone angles $\frac{2\pi}{n_1}t,\ldots,\frac{2\pi}{n_k}t$.

The proof of Theorem 2 consists in studying the behaviour of the
hyperbolic cone 3-manifold 
$X\left((\frac{n_1}t,0),\ldots,(\frac{n_k}t,0)\right)$ while
increasing  the parameter $t\in[0,1]$. If the parameter $t=1$ can be reached so
that the cone 3-manifold 
$X\left((n_1,0),\ldots,(n_k,0)\right)$ remains hyperbolic, then the
orbifold $\OO$ itself is hyperbolic. Otherwise, there is a limit
of hyperbolicity $t_{\infty}\in[0,1]$.
Then, since the space of hyperbolic cone structures with
topological type $(\Mod\OO,\Sigma)$ and cone angles $\leq 2\pi$ is
open, one has to analyze sequences of hyperbolic cone 3-manifolds 
$X\left((\frac{n_1}{t_n},0),\ldots,(\frac{n_k}{t_n},0)\right)$ where
$(t_n)_{n\in\N}$ is an increasing sequence in $[0,t_{\infty})$
approaching $t_{\infty}$. This analysis will be carried out in detail
in Chapter 2, by using  Theorems A and B below, which are
central to the proof of Theorem 2, and should be of independent interest.
Their proofs are given respectively in Chapter V and VI.

Theorem A is used when the limit of hyperbolicity $t_{\infty}<1$,
while Theorem B is used when $t_{\infty}=1$.

\proclaim{Theorem A} 
Let $(C_n)_{n\in\N}$ be a sequence of closed hyperbolic  cone 3-manifolds
with fixed topological type $(C,\Sigma)$
such that the cone angles increase and are contained in
$[\omega_0,\omega_1]$, with
$0<\omega_0<\omega_1<\pi$.
Then there exists a subsequence $(C_{n_k})_{k\in\N}$ such that one of 
the following occurs: 
\roster
  \item"1)" The  sequence $(C_{n_k})_{k\in\N}$ converges 
geometrically to a hyperbolic cone 3-man\-ifold with topological type
$(C,\Sigma)$ whose cone angles are the limit of the cone angles of
$C_{n_k}$.
   \item"2)" For every $k$, $C_{n_k}$ contains an  embedded $2$-sphere $S_{n_k}^2\subset
C_{n_k}$
that intersects $\Sigma$ in three points, and the sum of the three cone
angles at $S_{n_k}^2\cap\Sigma$ converges to $2\pi$.
   \item"3)" There is a sequence of positive  reals $\lambda_k$ approaching
$0$ such that the subsequence of rescaled cone 3-manifolds
$(\lambda_k^{-1}C_{n_k})_{k\in\N}$ converges geometrically to a 
Euclidean cone 3-manifold of topological type $(C,\Sigma)$ and whose
cone angles are the limit of the cone angles of $C_{n_k}$.
\endroster
\endproclaim

\proclaim{Theorem B} Let $\OO$ be a closed orientable connected
irreducible very good
$3$-orbifold with topological type $(\Mod\OO,\Sigma)$ and
ramification indices $n_1,\ldots,n_k$. Assume that there exists a
sequence of hyperbolic cone 3-manifolds 
$(C_n)_{n\in\N}$  with the same topological type 
$(\Mod\OO,\Sigma)$ and
such that, for each component of $\Sigma$,  the cone angles form an
increasing sequence that converges   to $2\pi/n_i$ when $n$ approaches
$\infty$.

Then $\OO$ contains a non-empty compact essential 3-suborbifold
$\OO'\subseteq\OO$, which is not a product and which is either complete 
hyperbolic of finite volume, Euclidean,  Seifert fibered or $Sol$.
\endproclaim

As stated, these two theorems deal with geometric convergence of cone
3-man\-ifolds. Up to minor modifications, the term geometric
convergence stands for the pointed bilipschitz convergence
introduced by Gromov \cite{\GLP}. The following Compactness Theorem
plays a central role  in the proofs of Theorems A and B. It is  a
cone manifold version of Gromov's Compactness Theorem for
Riemannian manifolds with pinched sectional curvature (cf.
\cite{\GLP} and \cite{\Pe}). The proof of this Theorem is the main
content of Chapter III.

\proclaim{Theorem \rm (Compactness Theorem)} Given $a>0$ and $\omega\in
(0,\pi]$, if   $(C_n,x_n)_{n\in\N}$
 is a sequence of pointed cone
3-manifolds 
 with constant curvature in $[-1,0]$, cone angles in
$[\omega,\pi]$, and such that $\Inj (x_n)\geq a$,
then $(C_n,x_n)_{n\in\N}$
has a subsequence  that converges geometrically to a pointed cone
3-manifold
$(C_{\infty},x_{\infty})$.
\endproclaim

The Compactness Theorem is used to analyze sequences of hyperbolic
cone 3-mani\-folds which do not collapse. In the collapsing case, we
need to rescale the metric in order to apply the Compactness
Theorem. In fact we need a more precise result in the collapsing
case, analogous to the ``Local Approximation Proposition" of Cheeger 
and Gromov \cite{\CGv, Prop. 3.4},   which furnishes a description of the
(non-trivial) topology of neighborhoods of points with small
cone-injectivity radius. This is  the Local Soul Theorem,
which is the content of Chapter IV. 

\

\proclaim{Theorem \rm (Local Soul Theorem)}
Given $\omega\in (0,\pi)$, $\varepsilon>0$ and $D>1$
there exist $\delta=\delta(\omega,\varepsilon,D)>0$ and
$R=R(\omega,\varepsilon,D)>D>1$ such that, if $C$  is an
 oriented hyperbolic cone 3-manifold with cone
angles in $[\omega,\pi]$ and if $x\in C$ satisfies
$\Inj(x)<\delta$, then:

- either $C$ is $(1+\varepsilon)$-bilipchitz homeomorphic
to a compact Euclidean cone 3-man\-ifold $E$ of diameter 
$\Diam(E)\leq R \Inj(x)$;

- or $x$ has an open neighborhood $U_x\subset C$
which is $(1+\varepsilon)$-bilipschitz homeomorphic to the
normal cone fiber bundle $\NN_{\nu}(S)$, of radius $0<\nu<1$, of
the soul $S$ of a non-compact orientable Euclidean cone
3-manifold  with cone angles in $[\omega,\pi]$.
 Moreover, according to
$\dim(S)$, the Euclidean non-compact
cone 3-manifold  belongs to the following list:
\roster
\item"(I) " (when $\dim(S)=1$), $S^1\ltimes\Bbb R^2$,
$S^1\ltimes(\text{open cone disk})$ and the pillow (see
Figure IV.1), where $\ltimes$ denotes the metrically twisted product;
\item"(II) "(when $\dim(S)=2$)
\item"(i)" a product $T^2\times \Bbb R$;
$S^2(\alpha,\beta,\gamma)
\times\Bbb R$, with
$\alpha+\beta+\gamma=2\pi$;
$S^2(\pi,\pi,\pi,\pi)\times\Bbb R$;
\item"(ii)" the orientable twisted line bundle over the Klein bottle
$K^2\tilde\times \Bbb R$ or
over the projective plane with two silvered points
$\Bbb P^2(\pi,\pi)\tilde\times \Bbb R$;
\item"(iii)" a quotient by an involution
 of either  $S^2(\pi,\pi,\pi,\pi)\times\Bbb
R$, $T^2\times\Bbb R$ or $K^2\tilde\times \Bbb R$,
 that gives an orientable  bundle respectively over
either $D^2(\pi,\pi)$, an annulus, or a M\"obius strip,
with silvered boundary in the three cases (see Figure IV.2).
\endroster
Moreover the $(1+\varepsilon)$-bilipschitz
homeomorphism 
$f:U_x\to\NN_{\nu}(S)$ satisfies the inequality
$$
\max\big(\Inj(x),d(f(x),S),\Diam(S)\big)\leq \nu/D.
$$
\endproclaim

\

These two theorems, Compactness Theorem and Local Soul Theorem, are
the main ingredients in the proofs of Theorems A and B. The
assumption that cone angles are bounded above by $\pi$ is crucial
for the proof of both theorems (and cannot be removed).
Geometrically this property is related to convexity: when cone
angles are bounded above by $\pi$, then the Dirichlet polyhedron is
convex; moreover convex subsets of such cone 3-manifolds have nice
properties.

One more ingredient for the study of the collapsing case is a cone
manifold version of Gromov's Isolation Theorem \cite{\Gro, Sec. 3}
(cf. Proposition  2.3 of Chapter V). This involves the notion of simplicial
volume due to Gromov \cite{\Gro}.

An application of Theorem A together with Hamilton's Theorem 
(cf. Thm 3.2 of \cite\ZhouTwo)
is the following:

\proclaim{Proposition 1} Given $0<\omega_0<\omega_1<2\pi/3$, there
exists a positive constant $\delta_0=\delta_0(\omega_0,\omega_1)>0$
such that for any oriented closed hyperbolic cone 3-manifold with
cone angles in $[\omega_0,\omega_1]$ there is a point $x_0\in C$ with
$\Inj(x_0)\geq\delta_0>0$.
\endproclaim

\newpage
\rightheadtext{ }
\leftheadtext{II \qquad    Proof of Thurston's orbifold theorem}

\

\centerline{\smc chapter  \  ii} 

\

\centerline{\chapt PROOF \,  OF  \, THURSTON'S \, 
ORBIFOLD  \,  THEOREM}

\

\

In this chapter we prove Thurston's
Orbifold  Theorem assuming Theorems A and B. We prove first the particular case of Theorem 2
and then we deduce from this the general case (Theorem 1).

\head 1. Generalized Hyperbolic Dehn Filling
 \endhead

Let $M$ be a compact 3-manifold with non-empty boundary 
$\partial M=T^2_1\sqcup\cdots\sqcup T^2_k$ a union of tori, whose
interior is hyperbolic (complete with finite volume).
Thurston's  hyperbolic Dehn filling Theorem provides a parametrization of
a space of hyperbolic deformations of this structure on $\Int(M)$.

 To describe the deformations on the ends of $\Int(M)$, we
fix two simple closed curves 
$\mu_i$ and $\lambda_i$ on each torus
$T^2_i$ of the boundary, which generate $H_1(T^2_i,\Bbb Z)$.
The structure around the $i$-th end of $\Int(M)$ is
described by the generalized Dehn filling coefficients 
$(p_i,q_i)\in\Bbb R^2\cup\{\infty\}=S^2$, such that 
the structure at the $i$-th end is
complete if and only if $(p_i,q_i)=\infty$.
The interpretation of the coefficients $(p_i,q_i)\in\Bbb R^2$ is the following:
\roster
\item"--" If $p_i,q_i\in\Bbb Z$ are coprime, 
then the completion at the
$i$-th torus is a non- singular hyperbolic
 3-manifold, which topologically 
is the Dehn filling with surgery meridian 
$p_i\mu_i+q_i\lambda_i$.
\item"--" When $p_i/q_i\in\Bbb Q\cup\{\infty\}$, 
let $r_i,s_i\in\Bbb Z$
be coprime integers such that $p_i/q_i=r_i/s_i$. 
The completion 
is a cone 3-manifold obtained by gluing a torus with 
singular core. The surgery
meridian is $r_i\mu_i+s_i\lambda_i$ and the cone 
angle of the singular component
is $2\pi\Mod{r_i/p_i}$.
\item"--" When $p_i/q_i\in\Bbb R-\Bbb Q$,
 then the 
completion (by equivalence classes of Cauchy sequences) is not
topologycally a manifold. These singularities 
are called of Dehn type, cf. \cite{\HodTwo}.
\endroster

\proclaim{Theorem 1.1} {\rm (Thurston's  Hyperbolic Dehn Filling Theorem
\cite{\ThuNotes})} There exists a neighborhood $U\subset S^2\times\cdots\times S^2$
of $\{\infty,\cdots,\infty\}$ such that the complete hyperbolic structure on $\Int(M)$
has a space of hyperbolic deformations parametrized by $U$ via generalized
Dehn filling coefficients.
\endproclaim
 
The proof yields not only the existence of a one parameter family of cone
3-manifold structures but also gives a path of corresponding holonomies in the variety $R(M)$ of
representations of $\pi_1(M)$ into $SL_2(\Bbb C)$. The holonomy of the complete structure on
$\Int(M)$ is a representation of $\pi_1(M)$ into $PSL_2(\Bbb C)$ that can be lifted to $SL_2(\Bbb
C)$. A corollary of the proof of Thurston's  hyperbolic Dehn filling Theorem is the
following:

\proclaim{Corollary 1.2} For any real numbers $\alpha_1,\ldots,\alpha_k\geq 0$
there exist $\varepsilon>0$ and  a path
$\gamma:[0,\varepsilon)\to R(M)$, such that,
for every $t\in [0,\varepsilon)$, $\gamma(t)$ 
is a lift of the holonomy of a hyperbolic
structure on M corresponding to the generalized Dehn filling coefficients  
$$
\big( (p_1,q_1),\ldots,(p_k,q_k) \big) =
\big( (2\pi/(\alpha_1 t), 0),\ldots,(2\pi/(\alpha_k t), 0) \big).
$$
\endproclaim

When $\alpha_i t=0$, the structure at the $i$-th cusp is complete; otherwise
its completion is a cone 3-manifold obtained by adding to $T^2_i$ a  solid torus with 
meridian curve $\mu_i$ and singular core  with cone angle  $\alpha_i t$.

\rightheadtext{II \qquad    Proof of Thurston's orbifold theorem}

\head 2. The space of hyperbolic cone structures
\endhead

Let $\OO$ be an  irreducible  orientable connected  closed 
3-orbifold such that $\OO-\Sigma$  admits a complete hyperbolic
structure of finite volume. In this section we study the space of hyperbolic cone
structures with topological type $(\OO,\Sigma)$. The main result of this section,
Proposition 2.1, can be deduced from Thurston's  Hyperbolic Dehn Filling
Theorem (Theorem 1.1) and Hodgson-Kerckhoff rigidity Theorem \cite{\HK}.
Nevertheless we present here
an elementary   proof, based  only on Thurston's Dehn filling Theorem, but
independent of Hodgson-Kerckhoff rigidity Theorem.

\definition{Notation} Let 
$m_1,\ldots,m_q$ be the
ramification indices of $\OO$ along
$\Sigma$. We set
$$
\alpha=(\alpha_1,\ldots,\alpha_q)=(\frac{2\pi}{m_1},\ldots,
\frac{2\pi}{m_q}).
$$

For $t\in[0,1]$, let $C(t\alpha)$ denote the hyperbolic
cone 3-manifold having the same topological type as the
orbifold $\OO$ and cone angles 
$t\alpha=(t\alpha_1,\ldots,t\alpha_q)$ (the order of
the components of $\Sigma$ is fixed throughout this
section). With this notation, $C(0)$ is the
complete hyperbolic structure of finite volume on
$\OO-\Sigma$.
\enddefinition

Thurston's hyperbolic Dehn filling Theorem (Corollary 1.2) means that
for small values of $t>0$ the hyperbolic cone 3-manifold
$C(t\alpha)$ exists. Thurston's idea is to increase
$t$ whilst keeping $C(t\alpha)$ hyperbolic  and to study  the
limit of hyperbolicity. 

More precisely, consider the variety of representations
of $\pi_1(\OO-\Sigma)$ into $SL_2(\Bbb C)$,
$$
 R:=\Hom (\pi_1(\OO-\Sigma),SL_2(\Bbb C)).
$$
Since $\pi_1(\OO-\Sigma)$ is finitely generated, $R$ is an
affine algebraic subset of $\Bbb C^N$ (it is  not
necessarily irreducible). The holonomy representation of
the complete hyperbolic structure on $\OO-\Sigma$ lifts
to a representation $\rho_0$ into $SL_2(\Bbb C)$, that is a point of
$R$. Let $R_0$ be an irreducible  component of $R$
containing $\rho_0$.

\definition{Definition} Define the subinterval
$J\subseteq[0,1]$ to be:
$$
J:=\left\{  t\in [0,1] \,\left\vert
\matrix\format \l\\
         \text{there exists a path }\gamma: [0,t] \to R_0 \\
\text{such that, }
         \text{for every }s\in[0,t], \\
         \gamma(s)\text{ is a lift  of the holonomy}\\
          \text{of a hyperbolic cone 3-manifold }
C(s\alpha) 
\endmatrix
\right.
\right\}
$$
\enddefinition

\remark{Remark} We say \it ``a" \rm hyperbolic cone 3-manifold
$C(s\alpha)$,
since we do not use the unicity of the hyperbolic cone structure for $s>0$ 
(cf. \cite{\KoOne}).
\endremark

By hypothesis, 
$J\neq\emptyset$ because $0\in J$ (i.e. $\OO-\Sigma$ 
has a complete hyperbolic structure).

\proclaim{Proposition 2.1} The interval $J$ is open in
$[0,1]$.
\endproclaim

\demo{Proof} The fact that $J$ is open at the origin is
a consequence of Thurston's hyperbolic Dehn
filling Theorem, as seen in Corollary 1.2.

Let $\mu=(\mu_1,\ldots,\mu_q)$ be the meridians of
$\Sigma$. That is, $\mu_i\in\pi_1(\OO-\Sigma)$
represents a meridian of the $i$-th component of
$\Sigma$, for $i=1,\ldots,q$ . Note that $\mu_i$ is not unique, only
the conjugacy class of $\mu_i^{\pm 1}$ is unique. We consider
the regular map:
$$
\matrix \format\r&\,\c\,&\l\\
  Tr_{\mu}: R_0 & \rightarrow  & \Bbb C^q\\
       \rho& \mapsto& (\Trace(\rho(\mu_1)),\ldots,
               \Trace(\rho(\mu_q))). 
\endmatrix
$$

\proclaim{Claim 2.2} There exists a unique affine irreducible
curve $\CC\subset\Bbb C^q$ such that, for any
$t\in J$,
$Tr_{\mu}(\gamma([0,t]))\subset \CC$.
\endproclaim

\demo{Proof of the claim}
For $n\in\N$, consider the  Chebyshev-like polynomial
$p_n(x)=2\cos(n \arccos(x/2))$. It is related to the
classical Chebyshev polynomial by a linear change of
variable. It can also be defined inductively by the  rule
$$
\left\{
\matrix\format\l\\
p_0(x)=2, \qquad p_1(x)=x,\\
p_n(x)=x \,p_{n-1}(x)- p_{n-2}(x), \text{ for }n\in\N,
n>1.
\endmatrix
\right.
$$
We are interested in the following property:
$$
\Trace(M^n)=p_n(\Trace(M)),
\qquad\forall M\in SL_2(\Bbb C),\ \forall n\in\N.
$$

Let $\rho_0=\gamma(0)$ be the lift of the holonomy
corresponding to the complete finite volume hyperbolic structure on
$\OO-\Sigma$. Since $\rho_0$ applied to a meridian is
parabolic, $Tr_{\mu}(\rho_0)=Tr_{\mu}(\gamma(0))=(\epsilon_1 2,\ldots,
\epsilon_q 2)$, with $\epsilon_1,\ldots,\epsilon_q\in\{\pm 1\}$. 

We take $\CC$ to be the irreducible component of the
algebraic set 
$$
 \{z\in\Bbb C^q\mid 
p_{m_1}(\epsilon_1 z_1)=\cdots=p_{m_q}(\epsilon_q z_q)\}
$$
that contains $Tr_{\mu}(\rho_0)$. 

To show that the component $\CC$ is well defined and is a
curve, we use the following identity:
$$
  p'_{n}(2)=  n^2, \qquad \forall n\in\N.
$$
It follows  from this formula that
$Tr_{\mu}(\rho_0)$ is a smooth point of 
$\{
p_{m_1}(\epsilon_1 z_1)=\cdots=p_{m_q}(\epsilon_q z_q)\}$ of local dimension
$1$.  Thus $\CC$ is the only irreducible component containing
$Tr_{\mu}(\rho_0)$ and is a curve.

 Finally to prove that $Tr_{\mu}(\gamma([0,t]))\subset \CC$, we
 consider the analytic map:
$$ \matrix \format\r&\,\c\,&\l\\
    \varTheta:\Bbb C&\rightarrow& \Bbb C^q \\
       w&\rightarrow& (\epsilon_1 2\cos( w \pi/m_1),\ldots,
            \epsilon_q 2\cos( w \pi/m_q))
\endmatrix
$$
Since $\gamma(s)(\mu_i)$ is a rotation of angle $s\pi/m_i$
 and $Tr_{\mu}(\gamma(0))=(\epsilon_1 2,\ldots,
\epsilon_q 2)$,
it is clear that
$Tr_{\mu}(\gamma([0,t]))\subset\varTheta(\Bbb C)$. By
construction,
 $$
\varTheta(\Bbb C)\subset\{
p_{m_1}(\epsilon_1 z_1)=\cdots=p_{m_q}(\epsilon_q z_q)\}.
$$
Since analytic irreducibility  implies algebraic
irreducibility, $\varTheta(\Bbb
C)\subset\CC$, and the claim is proved. \qed
\enddemo

\proclaim{Claim 2.3} For every $t\in J$, there exists an
affine curve ${\Cal D}\subset R_0$ containing $\gamma(t)$
and such that the restricted map $Tr_{\mu}:{\Cal D}\to\CC$
is dominant.
\endproclaim

\demo{Proof} We distinguish two cases, according to whether
$t>0$ or $t=0$.

When $t>0$, we take an irreducible component $\Cal Z$ of 
$Tr_{\mu}^{-1}(\CC)$
that contains the path $\gamma([t-\varepsilon,t])$, for some
$\varepsilon>0$. Since 
$$
Tr_{\mu}(\gamma(s))=(\epsilon_1 2\cos(s\pi/m_1),\ldots,
\epsilon_q 2\cos(s\pi/m_q)),
$$  the rational map
$
  Tr_{\mu}:{\Cal Z} \to \CC
$
is not constant, hence dominant. By considering generic
intersection with hyperplanes we can find the curve $\Cal D$
of the claim. More precisely, we intersect $\Cal Z$ with a
generic hyperplane $H$ passing through $\gamma(t)$ and such that
it does not contain any irreducible component
$Tr_{\mu}^{-1}(Tr_{\mu}(\gamma(t)))\cap {\Cal Z}$. 
If $\dim(\Cal Z)\geq 2$, then
such a hyperplane $H$
exists because $
  Tr_{\mu}:{\Cal Z} \to \CC
$ is dominant. By construction $Tr_{\mu}:{\Cal Z}\cap H  \to
\CC$ is not constant and the dimension of ${\Cal Z}\cap H$ is
less than the dimension of $\Cal Z$. By induction we obtain a
curve $\Cal D$.

When $t=0$, we consider again an irreducible component ${\Cal Z}$ of 
$Tr_{\mu}^{-1}(\CC)$ that contains $\rho_0$. 
In this case Thurston's  hyperbolic Dehn filling Theorem
implies that 
 the restriction $
  Tr_{\mu}:X \to \CC
$ is not constant. 
By considering
intersection with generic hyperplanes as before, we obtain
the curve $\Cal D$ of the claim. \qed
\enddemo

We now conclude the proof of Proposition 2.1. Given
$t\in J$, let $\CC$ and $\Cal D$ be as in Claims 2.2 and
2.3. Since
$\Cal D$ and $\CC$ are curves, for some $\varepsilon>0$, 
the path
$$
\matrix\format\r&\,\c\,&\l\\
    g:[t,t+\varepsilon)&\to&\CC\subset\Bbb C^q\\
    s&\mapsto&(\epsilon_1 2\cos(s\pi/m_1),\ldots,\epsilon_q 2
\cos(s\pi/m_q))
\endmatrix
$$
 can be lifted through $Tr_{\mu}:{\Cal D}\to\CC$ to a map
$\tilde g:[t,t+\varepsilon)\to\Cal D$. This map $\tilde g$
is a continuation of $\gamma$.

 It remains to check that
this algebraic continuation of $\gamma$ corresponds to the
holonomy representations of hyperbolic cone 3-manifolds.
To show this, we use Lemma 1.7.2 of Canary,
Epstein and Green's Notes on notes of Thurston
\cite{\CEG}. This lemma (called ``Holonomy induces
structure")  proves  that, for $s\in
[t,t+\varepsilon)$, $\tilde g(s)$ is the holonomy of a
hyperbolic structure on the
complement $C(t\alpha)-\NN_r(\Sigma)$  of a tubular neighborhood of the singular
set, with arbitrarily small radius $r>0$. The construction of the
hyperbolic structure in a tubular neighborhood of the singular set needs
some careful but elementary analysis. In
\cite{\PoTwo} this is done in the case where the structure is deformed
 from Euclidean to hyperbolic geometry, and the
constant curvature case is somewhat simpler.
This finishes the proof of Proposition 2.1.\qed
\enddemo

\remark{Remark} It follows from the proof that, for every
$t\in J$, the path
$\gamma:[0,t]\to R_0$ of holonomies of hyperbolic cone
structures is piecewise analytic. This is useful for applying
Schl\"afli's formula  \cite{\PoTwo, Proposition 4.2}).
\endremark

The following technical lemma will be used in the next
section.

\proclaim{Lemma 2.4} The dimension of $Tr_{\mu}^{-1}(\CC)$
is 4.
\endproclaim

\demo{Proof}  The proof of 
Thurston's  hyperbolic Dehn filling Theorem
 uses the fact
that the local dimension of
 $Tr_{\mu}^{-1}(Tr_{\mu}(\rho_0))$ at $\rho_0$ is 3,
because of Weil's local rigidity Theorem.
Moreover, the dimension of the preimage of any point in $\CC$ is
at least 3, because this is the dimension of $SL_2(\Bbb C)$.
Since $Tr_{\mu}:Tr_{\mu}^{-1}(\CC)\to\CC$ is dominant,  the
dimension of
$Tr_{\mu}^{-1}(\CC)$ is 4.\qed
\enddemo

\head 3. Proof of Theorem 2  \endhead

\proclaim{Theorem 2} Let $\OO$ be a closed orientable 
connected irreducible very good  3-orbifold of
cyclic type. Assume that the complement $\OO-\Sigma$ 
of the branching locus admits a complete
hyperbolic structure. Then $\OO$ contains a non-empty compact
essential $3$-suborbifold
$\OO'\subseteq\OO$, which is not a product and which is either complete
hyperbolic of finite volume, Euclidean,  Seifert fibered or $Sol$.
\endproclaim

\demo{Proof of theorem 2} We start with the subinterval $J\subseteq [0,1]$
 as in
Section 2. Recall that  
$J$ is the set of real numbers $t\in [0,1]$ such that there is a path
$\gamma:[0,t]\to R_0$ with the property that, for every
$s\in [0,t]$, $\gamma(s)$ is the holonomy 
of a hyperbolic cone 3-manifold $C(s\alpha)$.  The
hyperbolic cone 3-manifold $C(s\alpha)$ has the same
topological type as
$\OO$ and its cone angles are $s\alpha=(s\,
2\pi/m_1,\ldots,s\,2\pi/m_q)$.

By Proposition 2.1, $J$ is open in $[0,1]$. So there are
three possibilities: either $J=[0,1]$, $J=[0,1)$, or $J=[0,t)$
with $0<t<1$.

If $J=[0,1]$  then $\OO$ is hyperbolic.
Propositions 3.1 and 3.5 deal with the cases where 
 $J=[0,t)$
with $0<t<1$ and $J=[0,1)$ respectively.

\proclaim{Proposition 3.1} If $J=[0,t)$  with $0<t<1$,
then  $\OO$ is a spherical 3-orbifold.
\endproclaim

\demo{Proof}
Fix $(t_n)_{n\in\N}$ an increasing sequence in $J=[0,t)$
converging to
$t$ and consider the corresponding sequence of cone
3-manifolds
$C_n= C(t_n\alpha)$. The cone 3-manifolds $C_n$ have the same
topological type as $\OO$, and the cone angles are
contained in some interval $[\omega_0,\omega_1]$, with
$0<\omega_0<\omega_1<\pi$, because $0<t<1$. Thus we can
apply Theorem A to the sequence $(C_n)_{n\in\N}$, and, after perhaps passing to  a
subsequence, we have three possibilities:
\roster
\item"i)" the sequence $(C_n)_{n\in\N}$ converges
geometricaly to a hyperbolic cone 3-man\-ifold  with the same
topological type;
\item"ii)" each $C_n$ contains an embedded sphere $S_n$
which intersects $\Sigma$ in 3 points and the sum of the
cone angles at these points converges to $2\pi$;
\item"iii)" there is a sequence of positive
reals $\lambda_n\to 0$ such that
$\big(\frac1{\lambda_n}C_n\big)_{n\in\N}$ converges
geometrically to a Euclidean cone 3-manifold with the same
topological type.
\endroster
We want to show that only the last possibility occurs.

If the case i) happens, we claim that 
$t\in J$; this would contradict the hypothesis
$J=[0,t)$. Let $C_{\infty}$ be the limit of the sequence $C_n$.
Since the convergence is geometric, the cone angles of 
$C_{\infty}$ are precisely $t\alpha$. Therefore
$C(t\alpha)$ is hyperbolic and it remains to show the
existence of a path of holonomy representations from $0$ to
$t$. To show that, we take a path $\gamma_n$ for each $\rho_n$
and  we prove  that the sequence of paths
$\gamma_n$ has a convergent subsequence. 

\proclaim{Lemma 3.2} The sequence of paths $\gamma_n$
has a subsequence converging to a path $\gamma_{\infty}$.
Moreover, up to conjugation,  for $n$ sufficiently large,
$\gamma_{\infty}$ is a continuation of $\gamma_n$.
\endproclaim

\demo{Proof of the lemma} Consider the algebraic affine set
$V=Tr_{\mu}^{-1}(\CC)$ and its quotient by conjugation
$X=V/PSL_2(\Bbb C)$. The space $X$ may not be
Hausdorff, but since the holonomy of a closed hyperbolic
cone 3-manifold is irreducible \cite{\PoOne, Proposition 5.4},
the points we are interested in (conjugacy classes of
holonomy representations) have neighborhoods that are
analytic (see for instance \cite{\CS} or \cite{\PoOne, Proposition
3.4}). If we remove all reducible
representations, then the quotient is analytic, even
affine \cite{\CS}, call it $X^{irr}$.
By lemma 2.4, $X^{irr}$ is a curve. Moreover, the
holonomies of hyperbolic cone structures are contained in
a real  curve of $X^{irr}$, because the traces of the
meridians are real. Hence, up to conjugation, the paths
$\gamma_n$ are contained in a  real analytic curve. Thus,
 $\rho_n$ converges to $\rho_{\infty}$, the sequence
$\gamma_n$ has a convergent subsequence, and the limit 
$\gamma_{\infty}$ is a
continuation of $\gamma_n$.
\qed
\enddemo

It follows from the lemma that the limit
$\gamma_{\infty}$ is a path of holonomy representations of
 hyperbolic cone structures. Hence $t\in J$ and
we obtain a contradiction.

Next we suppose that case ii) occurs. That is, for each
$n\in\N$,
$S_n^2\subset C_n$ is an embedded $2$-sphere which intersects
$\Sigma$ in three points and the sum of the
cone angles at these points converges to $2\pi$. Since $\Sigma$ has a finite number of
components, after passing to a subsequence, we can suppose that $S_n^2$
intersects always the same components of $\Sigma$. Let
$m_1$, $m_2$ and $m_3$ be the branching indices of the
components of $\Sigma$ which intersect $S_n^2$, counted
with multiplicity if one component of $\Sigma$ intersects
$S_n^2$ more than once. Since the sum of  the cone angles at these points
converges to
$2\pi$, we have:
$$
2\pi\leq
t\bigg(
\frac{2\pi}{m_1}+\frac{2\pi}{m_2}+\frac{2\pi}{m_3} 
\bigg)
<
2\pi\bigg(
\frac{1}{m_1}+\frac{1}{m_2}+\frac{1}{m_3} 
\bigg),
$$
because $t<1$. If we view 
 $S_n^2$ as a 2-suborbifold $F\subset \OO$, then $F$ is
spherical, because its underlying space is
$\vert F\vert\cong S^2$ and it has three singular points
with branching indices
$m_1$, $m_2$ and $m_3$, where
$\frac1{m_1}+\frac1{m_2}+\frac1{m_3}>1$. The suborbifold
$F$ cannot bound a discal
$3$-orbifold, since we assume that the ramification set of
$\OO$ is a link. This contradicts the irreducibility of the 3-orbifold
$\OO$, so the case ii) cannot happen.

Thus so far we have eliminated cases i) and ii). Now we prove from
case iii) that $\OO$ is a spherical 3-orbifold. Case iii) implies that
there is a Euclidean cone 3-manifold  $C(t\alpha)$ with the
same topological type as $\OO$ and with cone angles
$t\alpha=(t\, 2\pi/m_1,$ $\ldots,t\, 2\pi/m_q)$, where
$0<t<1$.

Let $M\to \OO$ be a very good regular covering of $\OO$
with finite deck transformation group $G$. Since $t<1$, the
Euclidean  cone 3-manifold $C(t\,\alpha)$ induces a
$G$-invariant  Euclidean cone manifold structure on $M$, with
singular angles $t\, 2\pi<2\pi$, from which we deduce (cf. (\cite{\Jon},
\cite{\GT}, \cite{\Zhou,2}):

\proclaim{Lemma 3.3} The manifold $M$ admits a non-singular
$G$-invariant Riemannian metric with non-negative
sectional curvature that is not flat.
\endproclaim

\demo{Proof}
We shall deform the singular Euclidean
metric on
$M$ (induced by $C(t\alpha)$) in a
$G$-invariant way.
Let $\Sigma_G\subset M$ be the singular set of this
Euclidean metric, which is also the set of points
where the action of $G$ is not free.
 We  deform the metric in a tubular
neighborhood of the singular set $\NN_{r_0}(\Sigma_G)$
 of radius $r_0$, for some $r_0>0$ sufficiently small. 
Around $\Sigma_G$, the local expression of the
singular Euclidean metric in Fermi (cylindrical)
coordinates is
$$
ds^2=dr^2+t^2r^2d\theta^2+dh^2,
$$
where $r\in (0,r_0)$  is the distance from $\Sigma_G$, $h$ is
the length parameter along $\Sigma_G$, and 
$\theta\in(0,2\pi)$ is the rescaled angle parameter.

The deformation that we are going to introduce depends only on
the parameter $r$, so it is $G$-invariant. This
deformation consists of replacing the above metric by a metric
of the form
$$
ds^2=dr^2+f^2(r)d\theta^2+dh^2,
$$
where $f:[0,r_0-\varepsilon)\to[0,+\infty)$ is a smooth
function which satisfies, for some $\varepsilon>0$
sufficiently small:
\roster
\item"1)" $f(r)=r$, for $r\in [0,\varepsilon)$;
\item"2)"  $f(r)=t(r+\varepsilon)$, for $r\in (r_0/2, r_0-\varepsilon)$;
\item"3)" $f$ is concave: $f''(r)\leq 0$, for all $r\in
[0,r_0-\varepsilon)$.
\endroster
Such a function $f$ exists because $0<t<1$. The first
property implies that the new metric is non-singular. 
Property 2) implies that, after reparametrization, this new
non-singular metric fits with the original singular metric
at the boundary of the tubular neighborhood
$\NN_{r_0}(\Sigma_G)$.  A classical
computation shows that the sectional curvature of the
planes orthogonal to $\Sigma_G$ is non-negative, by
property 3), and it is even positive at some point. Hence,
since the metric is locally a product, it has non-negative
sectional curvature.
\qed
\enddemo

To show that $\OO$ is spherical, we apply the
following deep theorem of Hamilton 
(\cite{\HamOne} and \cite{\HamTwo}, see also \cite{Bou}).

\proclaim{Theorem 3.4 \rm (Hamilton) \cite{\HamOne,\HamTwo}}
Let $N^3$ be a  closed 3-manifold which admits
a Riemannian metric of non-negative Ricci curvature. Then $N^3$ admits a metric
which is either spherical, flat or modelled on $\Bbb S^2\times\Bbb R$.

 Furthermore the deformation is natural, and every isometry
of the original metric is also an isometry of the new metric.
 \qed
\endproclaim

\remark{Remark} It follows from Hamilton's proof \cite{\HamTwo} that the flat
case occurs only if the initial metric on $N^3$ was already
flat, and that the case modelled on $\Bbb S^2\times \Bbb R$ occurs
only if the initial metric had reducible holonomy, contained in $SO(2)$.
\endremark

We apply Hamilton's Theorem 3.4 to the metric on $M$ given by Lemma 3.3
and we claim that  $M$ admits a G-invariant spherical metric. According to the
remark, the flat case of Hamilton's Theorem does not occur because the initial
metric was not flat. Moreover, we can also eliminate the case $S^2\times \Bbb R$,
because this case would imply that the singular Euclidean cone structure of $M$ is
of Seifert type ($M$ admits a Seifert fibration such that the singular locus is an union of
fibers). This follows easily by considering isometries of $ \Bbb S^2\times\Bbb R$
or from \cite{\PoTwo, Lemma 9.1}.

Thus
Hamilton's Theorem implies that
$M$ has a
$G$-invariant spherical metric, hence  the  3-orbifold $\OO$ is
spherical. This concludes the proof of Proposition 3.1.\qed
\enddemo
\enddemo

\remark{Remark} In this case ($0<t<1$) the hypothesis that $\OO$ is very good
is not necessary. In fact, we can apply Hamilton's Theorem to the Euclidean cone
3-manifold $C(t \alpha)$ as above to conclude that the underlying space
of $\OO$ is spherical. Thus, up to passing to a finite cover, we can assume 
that the underlying space
of $\OO$ is $S^3$. Since the ramification set is a link, $\OO$ is very good.
In particular, when all ramification indices are at least $3$
and $\OO$ does not contain the toric suborbifold $S^2(3,3,3)$
 we do not need the
hypothesis very good.
 \endremark

We complete now the proof of Theorem 2 by dealing with the case where
$J=[0,1)$:

\proclaim{Proposition 3.5}
If $J=[0,1)$ then $\OO$ 
contains a non-empty compact
 essential 3-suborbifold which is not a product and which is
geometric. 
\endproclaim

\demo{Proof}
Let $(t_n)_{n\in\N}$ be a sequence in $[0,1)$ converging to
$1$. We apply Theorem B to the corresponding sequence of
hyperbolic cone 3-manifolds $(C(t_n\alpha))_{n\in\N}$ 
whose cone angles form an increasing sequence that 
converges to $2\pi/n_i$, $i=1,\ldots,k$, when $n$ goes 
to $\infty$. By Theorem B, $\OO$ contains a non-empty compact essential
3-suborbifold $\OO'\subseteq	\OO$ which is either Euclidean,
Seifert fibered, $Sol$,  or hyperbolic (of finite volume), and which is
not a product. 

This concludes the proof of Proposition 3.5 and of Theorem 2.
\qed
\enddemo

\head 4. Proof of  Thurston's Orbifold Theorem
\endhead

\proclaim{Theorem 1 \rm (Thurston's Orbifold Theorem)} Let $\OO$ be a 
compact  connected  orientable 
irreducible 
boundary-incom\-pressible 3-orbifold of cyclic type. If $\OO$ is
very good, topologically atoroidal and
 acylindrical, then $\OO$ is geometric (i.e.  $\OO$ admits
 either a hyperbolic, a Euclidean, or a Seifert fibered structure).
\endproclaim

Throughout this section we assume that $\OO$ is a 3-orbifold which satisfies 
the hypothesis of Theorem 1. 

Let
$D\OO$  denote the double of $\OO$ along some components of
$\partial\OO$, which we call \it doubling components \rm. 
The
ramification set of $D\OO$ is denoted by $D\Sigma$. If
we double along the empty set, then we choose the convention that
$D\OO=\OO$, so that $D\OO$ is always connected.
 
First we give some results about the topology of $D\OO$ and
$D\OO-D\Sigma$ in order to reduce the general case
to the case where the hypothesis of Theorem 2 are satisfied.
Then we deduce Theorem 1 from Theorem 2.

\proclaim{Lemma 4.1} For any choice of doubling components,
\roster 
\item"(i)" $D\OO$ is irreducible, topologically acylindrical and very good;
\item"(ii)" every component of $\partial\OO$ is
incompressible in $D\OO$; 
\item"(iii)" every incompressible 
toric 2-suborbifold of $D\OO$ is parallel to
$\partial\OO\subset D\OO$.
\endroster
In particular $D\OO$ is
bound\-ary-incompressible. Furthermore,  $D\OO$ is topologically
atoroidal if and only if every doubling component is a hyperbolic
2-suborbifold. 
\endproclaim

\demo{Proof}
 Let $S\subset D\OO$ be a spherical 2-suborbifold. After
isotopy, we can suppose that $S$ is transverse to
$\partial\OO$ and that the intersection $S\cap\partial\OO$ is
minimal. We claim that $S\subset\OO$. Seeking a
contradiction, we suppose $S\not\subset\OO$.
Since
$S$ is a sphere with at  most three cone points, at least
one component of
$S\cap\OO$ is a disk  
$\Delta^2$ with at most one cone point. 
Since $\OO$ is irreducible
and the intersection $S\cap\partial\OO$ is minimal, $\partial\Delta^2$ is essential in
$\partial\OO$. Hence 
$\Delta^2$ is a compressing disk for $\partial\OO$, and we
get a  contradiction because
$\OO$ is boundary-incompressible.
Therefore $S\subset\OO$ and $S$ bounds a discal
$3$-orbifold by irreducibility of $\OO$. The same
argument applies to show that $D\OO$ does
not contain any bad $2$-suborbifold.
Hence
$D\OO$ is also irreducible.

Let $A\subset D\OO$ be a properly embedded annular 2-suborbifold.
Again we deform it so that $A\cap\partial\OO$ is
transverse and minimal. No component of $A\cap\OO$ is
a discal orbifold, because $\partial\OO$ is incompressible
and the intersection $A\cap\partial\OO$ is minimal. Hence
$A=A_1\cup\cdots\cup A_k$, where each $A_i$ is an annular
2-suborbifold properly embedded in one of the copies
of
$\OO$. If $k>1$, then, 
by minimality of the intersection, none of the annuli $A_i$
is parallel to
$\partial\OO$ nor compressible in $\OO$, contradicting the
acylindricity of $\OO$. Hence $k=1$ and $A\subset \OO$ is
not essential. This proves that $D\OO$ is topologically acylindrical.

Note that $D\OO$ is very good, because  a regular
covering $M\to\OO$ induces a regular covering $DM\to D\OO$,
where $DM$ is the double of $M$ along the components of
$\partial M$ that project to doubling components of
$\partial\OO$. Assertion (i) is proved.

To show that every component of $\partial\OO$ is
incompressible in $D\OO$, suppose that $\partial \OO$ has
a compressing disk
$\Delta^2\subset D\OO$. By making the intersection 
$\Delta^2\cap\partial\OO$ minimal, every disk component of
$\Delta^2\cap\OO$  is a compressing disk for $\partial\OO$
in
$\OO$, thus we obtain a contradiction that proves
assertion (ii).

Finally, let $F\subset D\OO$ be an incompressible toric
2-suborbifold. After an isotopy, we can again make the intersection
$F\cap\partial\OO$  transverse and minimal. If
$F\cap\partial\OO\neq\emptyset$, then the
minimality of the intersection implies that each component of $F\cap\OO$ is
an essential  annular  2-suborbifold, and we get a
contradiction. Thus $F\subset\OO$ and $F$ is parallel to
$\partial \OO$.
\qed
\enddemo

\proclaim{Lemma 4.2} 
If every doubling component is different from a
non-singular  torus,
then
the manifold $D\OO-D\Sigma$ is
irreducible and topologically atoroidal.
\endproclaim

\demo{Proof}
Let $S\subset D\OO-D\Sigma$ be an embedded 2-sphere. It
bounds a discal 3-suborbifold $\Delta^3$ in $D\OO$, because
$D\OO$ is irreducible. Since
$S\cap D\Sigma=\emptyset$, $\Delta^3$ is a 3-ball and
$\Delta^3\cap D\Sigma=\emptyset$. So $D\OO-D\Sigma$ is
irreducible.

Let $T\subset D\OO-D\Sigma$ be an embedded torus. By Lemma 4.1(iii) either 
$T$ is compressible in
$D\OO$ or  parallel to a component of  $\partial\OO\subset D\OO$. In this last case, by
hypothesis, $T$ must be boundary parallel in $D\OO$, and thus it is also boundary parallel in
$D\OO-D\Sigma$. If $T$ admits a compressing disk (a
discal 2-suborbifold), then  the irreducibility of $D\OO$
implies  that
$T$  bounds either  a solid torus or  a solid
torus with singular core $S^1\times D^2(*)$. In the former
case 
$T$ is compressible in $D\OO-D\Sigma$; in the latter case, $T$
is boundary parallel in $D\OO-D\Sigma$.
\qed
\enddemo

\proclaim{Lemma 4.3} \roster
\runinitem"(i)" For any choice of doubling
components, if
$D\OO$ is Seifert fibered then $\OO$ is Seifert fibered,
  Euclidean, or an $I$-bundle over a $2$-orbifold.
\item"(ii)" If the doubling components are non-empty, then
$D\OO$ is not $Sol$.
\endroster
\endproclaim

\demo{Proof} First suppose that $D\OO$ is Seifert fibered.
We can assume that
$\partial\OO\neq\emptyset$, otherwise the statement is trivial. We consider the
natural involution
$\tau: D\OO\to D\OO$ obtained by reflection through the
doubling components of
$\OO$. Let $DM\to D\OO$ be the covering obtained by doubling
a regular, very good covering
$M\to\OO$. The reflection $\tau: D\OO\to D\OO$ lifts to
a reflection $\tilde\tau:DM\to DM$ which commutes with
the deck transformations group of
$M\to\OO$. The fundamental group $\pi_1(DM)$ is infinite, because 
each component of $\partial\OO$  lifts to incompressible surfaces in $DM$, with infinite 
fundamental group.
Hence, we can apply
Meeks-Scott Theorem \cite{MS} and conclude that
$D\OO$ has a Seifert or Euclidean structure invariant by
the natural involution $\tau: D\OO\to D\OO$. 

If the $\tau$-invariant structure of $D\OO$ is Euclidean,
then  $\OO$ is also Euclidean. If the $\tau$-invariant
structure of
$D\OO$ is Seifert fibered, then either $\OO$ is  Seifert fibered or
$\OO$ is a
$I$-bundle over a 2-orbifold $F^2$. This proves assertion (i).

To prove assertion (ii), suppose that $D\OO$ is $Sol$. There is a
finite regular covering $M\to D\OO$ which is a manifold and fibers over
$S^1$ with fiber $T^2$ and Anosov
monodromy.
Let  $\tau:D\OO\to D\OO$ be the natural involution as above, obtained
by reflection through the doubling components of $\partial \OO$.
 Since $\tau$ is an involution whith a non-empty fixed point set,  not included in the
ramification locus, it lifts to an  involution $\tilde \tau$ of $M$ whose 
fixed point set contains a two dimensional submanifold.
By Tollefson's Theorem about finite order homeomorphisms of fiber 
bundles \cite{\ToTwo} (see also \cite{MS}), we may assume that $\tilde \tau$ preserves the 
fibration by tori. Then one can easily check that torus bundles
with Anosov monodromy cannot admit the reflection 
$\tilde \tau$. 
 \qed
\enddemo

\remark{Remarks} When $\OO$ is an $I$-bundle over a 
2-orbifold $F^2$, the following facts should
be noted:
\roster
\item"(i)" The 2-orbifold $F^2$ is either Euclidean or
hyperbolic, because $\OO$ is irreducible. 
In particular,
the interior of $\OO$ has a complete Euclidean or 
hyperbolic structure.
\item"(ii)" acylindricity of $\OO$ restricts the
possibilities for $F^2$.  
\item"(iii)" The manifold $D\OO-D\Sigma$ is Seifert fibered.
\endroster
\endremark

\proclaim{Lemma 4.4} \roster \runinitem"(i)" For any choice
of doubling components, either $D\OO-D\Sigma$ 
 has an incompressible boundary or $\OO$ is Seifert fibered.
\item"(ii)" If every doubling component is
different from a non-singular torus, then either
$D\OO-D\Sigma$  is  topologically acylindrical or
$D\OO$ is Seifert fibered.
\endroster
\endproclaim

\demo{Proof} To prove the first assertion, we assume that
$D\OO-D\Sigma$  has a compressible boundary. Then it is a
solid torus, because $D\OO-D\Sigma$ is irreducible and its
boundary is a union of tori. Hence the underlying space of $D\OO$ is a
Lens space and its singular set is the core of
one of the solid tori of a  genus one Heegard splitting.
As the 3-orbifold $D\OO$ is very good, it cannot be the
product  $S^1\times S^2(*)$, where $S^2(*)$ is a 2-sphere
with one cone point (a bad 2-orbifold). Hence $D\OO$ has
$S^3$ as universal covering. Since the orbifold group
$\pi_1^o(D\OO)$
is finite, we have that $\partial\OO=\emptyset$ and  $D\OO=\OO$, by Lemma
4.1 (ii). Thus
$\OO$ is Seifert fibered.

Suppose that  $D\OO-D\Sigma$ contains an essential annulus.
We claim that $D\OO$ is Seifert fibered.
By Lemma 4.2 and the Characteristic Submanifold
Theorem \cite{\JS}, \cite{\Joh}, $D\OO-D\Sigma$ is Seifert fibered.
We consider a component $\Sigma_i$ of $\Sigma$ and a solid
torus neighborhood $\NN(\Sigma_i)$. If the fiber of the
Seifert fibration of $D\OO-D\Sigma$ is not homotopic to the meridian of
$\Sigma_i$ in the torus $\partial \NN(\Sigma_i)$, then this
Seifert fibration  can be extended to $\NN(\Sigma_i)$ 
so that
$\Sigma_i$ is a fiber. 

We suppose now that the  fiber of the
Seifert fibration of $D\OO-D\Sigma$  is homotopic to the
meridian of $\Sigma_i$. If  the base 2-orbifold of the 
Seifert fibration in $D\OO-D\Sigma$ is different from a
disk or a disk with one cone point, then $D\OO-D\Sigma$
contains an essential annulus which is vertical and its
boundary is in
$\partial
\NN(\Sigma_i)$. In particular, the union of this annulus
with two meridian disks (with one cone point) of
$\NN(\Sigma_i)$ gives an incompressible spherical
2-suborbifold in $D\OO$,  contradicting the
irreducibility of $D\OO$. Hence $D\OO-D\Sigma$ is a solid
torus, and, as we have already shown above, this implies that
$D\OO=\OO$ is Seifert.
\qed
\enddemo

\proclaim{Proposition 4.5} If Thurston's Orbifold
Theorem holds when  no component of 
$\partial\OO$ is a non-singular torus, then it holds
in general.
\endproclaim

\demo{Proof}
We decompose  the boundary of $\OO$ in three parts:
$$\partial
\OO=
\partial_T\OO\sqcup\partial_{SE}\OO\sqcup\partial_{H}\OO,
$$
where:
\roster
\item"-"
$\partial_T\OO$ is the union of the boundary
components homeomorphic to a torus
\item"-"
$\partial_{SE}\OO$ is the union of the singular Euclidean
boundary components
\item"-"
$\partial_{H}\OO$ is the union of  the hyperbolic boundary
components
\endroster
We assume $\partial_T\OO\neq\emptyset$.
We double along the hyperbolic components 
$\partial_{H}\OO$:
$$
D\OO=\OO\underset{\partial_{H}\OO}\to\cup\OO.
$$

Since $D\OO$ is very good, we fix a  regular covering
$p:DM\to D\OO$ of finite order which is a manifold, and let
$G$ denote its group of deck transformations.
Since $D\OO$ is irreducible, topologically atoroidal and 
boundary-incompressible  (lemma 1.1), the Equivariant  Sphere
and Loop Theorems (\cite{\DD},
\cite{\MYOne,2},
\cite{\JR}) imply that 
$DM$ is also irreducible, topologically atoroidal, and 
 boundary-incompressible. Since by hypothesis
$\partial (DM)\neq\emptyset$, Thurston's Hyperbolization theorem for Haken
3-manifolds \cite{\ThuNotes,2,3,4,5, \McM1,2, \OtaOne,2} implies that $DM$ 
is either Seifert fibered or hyperbolic. If  $DM$ is Seifert fibered, Meeks-Scott
Theorem \cite{\MS} implies that
$D\OO$ is also Seifert fibered, and thus $\OO$ is geometric (Lemma 1.3).
 Therefore we can assume  that $DM$ is hyperbolic.

Let $\gamma=\{\gamma_1,\ldots,\gamma_r\}$ be a family of
simple closed curves, one on each torus component
of
$\partial_T (D\OO)$. Let $D\OO(\gamma)$ denote the 3-orbifold
obtained by generalized Dehn filling with meridian curves
$\gamma=\{\gamma_1,\ldots,\gamma_r\}$. 
Generalized Dehn filling means that the filling solid
tori may have  ramified cores. Moreover,  we choose the branching indices of these filling cores
so that the generalized Dehn filling $D\OO(\gamma)$ lifts to a genuine  Dehn filling of $DM$.

We consider a sequence 
of families of simple closed  curves 
$(\gamma^n)_{n\in\N}=
(\{\gamma_1^n,\ldots,$ $\gamma_r^n\})_{n\in\N}$ such that,
for each $n\in\Bbb N$,
$\gamma^n$ gives precisely one curve on
each component of
$\partial_T D\OO$, and for each $i=1,\ldots,r$, the
curves of the sequence $(\gamma_i^n)_{n\in\N}$ represent
different homotopy classes on the $i$-th torus boundary
component. For  $n\in\Bbb N$ sufficiently large, 
the orbifold $D\OO(\gamma^n)$ has a regular covering obtained by Dehn filling of $DM$, which
we may assume to be hyperbolic by  Thurston's
Hyperbolic Dehn Filling Theorem. 
Then, by the proof of Smith conjecture
\cite{\MB} and the Equivariant  Loop Theorem \cite{\MYOne,2},
$D\OO(\gamma^n)$ is irreducible and topologically
atoroidal, for $n\in\N$ sufficiently large. Moreover, by construction, no component of
$\partial D\OO(\gamma^n)$ is a non-singular torus. Hence, for $n\in\N$ sufficiently large,
$D\OO(\gamma^n)$ is geometric  by hypothesis, and so it is hyperbolic.

For each $n\in\N$ sufficiently large, choose a point $x_n\in D\OO(\gamma^n)$
so that $\Inj(x_n)>\varepsilon_3$, where $\varepsilon_3>0$ 
is the 3-dimensional Margulis constant. By the Compactness
Theorem (Chap. III) there is a subsequence of the sequence $( D\OO(\gamma^n), x_n)$ which
converges geometrically to a hyperbolic 3-orbifold. Moreover the limit is non-compact and gives 
a hyperbolic structure on the interior of the 3-orbifold $D\OO$, because the sequence of
coverings of $ D\OO(\gamma^n)$ converges geometrically to the interior of 
$DM$ by Thurston's hyperbolic
Dehn filling Theorem. 

Since $D\OO$ is hyperbolic, the next lemma shows that $\OO$ is
also hyperbolic.
\qed
\enddemo

\remark{Remark} The hypothesis that $\OO$ is very  good is not necessary
for the proof of Proposition 4.5: we can show directly that $D\OO$ is either Seifert fibered or
topologically acylindrical by using Bonahon-Siebenmann Characteristic Toric Splitting for
3-orbifolds (\cite{\BSOne}). Then, for $n\in\N$ sufficiently large, it follows that the 3-orbifold 
$D\OO(\gamma^n)$ is irreducible, topologically atoroidal and acylindrical; hence by hypothesis it is
geometric and thus already very good. \endremark

\proclaim{Lemma 4.6} If $D\OO$ is hyperbolic with finite
volume,  then
$\OO$ is hyperbolic and the doubling components are totally
geodesic.
\endproclaim

\demo{Proof} Note that if $D\OO$ is hyperbolic with
finite volume, then the doubling components are precisely the
hyperbolic pieces of $\partial\OO$. We  
assume that $\partial\OO$ has hyperbolic components,
otherwise $\OO=D\OO$ and there is nothing to prove. 
Consider the reflection $\tau_0:D\OO\to
D\OO$  through the doubling components.
By Mostow-Prasad Rigidity Theorem, $\tau_0$ is homotopic (in the
orbifold sense) to an isometric involution $\tau_1:D\OO\to
D\OO$. We claim that these two involutions are in fact conjugate. This
claim implies that the 3-orbifold $\OO$ is hyperbolic and that the hyperbolic
components of $\partial\OO$ are totally geodesic. 

To prove the claim, we use an unpublished argument of Bonahon and
Siebenmann \cite{\BSThree}. Since both $\tau_0$ and $\tau_1$ preserve the
ramification set, they induce respectively an involution
of $D\OO-D\Sigma$. By Lemmas 4.2 and 4.4 and Thurston's
hyperbolization theorem, the manifold $D\OO-D\Sigma$ is
hyperbolic. Hence, by Mostow-Prasad Rigidity Theorem, the
restrictions of $\tau_0$ and $\tau_1$ to $D\OO-D\Sigma$
are respectively homotopic to the isometric involutions  $g_0$ and $g_1$  on
$D\OO-D\Sigma$. By  Waldhausen's and Tollefson's Theorems
(\cite{\Wal}, \cite{\ToOne}), there exist two homeomorphisms
$h_0,h_1:D\OO-D\Sigma\to D\OO-D\Sigma$ isotopic to the
identity such that the restrictions 
$\tau_0\vert_{D\OO-D\Sigma}=h_0g_0h_0^{-1}$ and
$\tau_1\vert_{D\OO-D\Sigma}=h_1g_1h_1^{-1}$. Therefore, the
involutions  $g_0$ and $g_1$ on $D\OO-D\Sigma$ can be extended 
respectively to involutions $\bar g_0,\bar
g_1:D\OO\to D\OO$. It remains to show that $\bar g_0=\bar
g_1$ on $D\OO$.

The map $f=\bar g_0\bar g_1^{-1}$ is homotopic to the identity on $D\OO$
in the orbifold sense; moreover $f$ is of finite order,
because its restriciton to
$D\OO-D\Sigma$ is an isometry. Since $f$ is homotopic to
the identity, it lifts to a homeomorphism $\tilde f:\Bbb
H^3\to\Bbb H^3$ whose extension to the sphere at infinity
$\partial\Bbb H^3\cong S^2$ is the identity. Since $f$ is of finite order, so
is $\tilde f$, because $\tilde f^n$ is an isometry  of $\Bbb H^3$ whose extension 
to the sphere at infinity is the identity. Since the identity is the only periodic  map of the
ball which is the identity on the boundary, $\bar g_0=\bar g_1$. Hence the involutions $\tau_0$
and $\tau_1$ are conjugate on $D\OO$
\qed
\enddemo

\demo{Proof of Theorem 1} We use the  previous results of this section 
to make some reductions of the general case. First,
  by Proposition 4.5, we can
assume that no component of $\partial\OO$ is a
non-singular torus.

Let $D\OO$ be the double of $\OO$ along all boundary
components. In particular $\partial (D\OO)=\emptyset$. By
Lemma 4.1,
$D\OO$ is irreducible and very good; moreover every
incompressible Euclidean 2-suborbifold is singular and
parallel to a doubling component.
By Lemma 4.2, $D\OO-D\Sigma$ is
irreducible and atoroidal. Furthermore, by Lemmas  4.4 and 4.3, we can
assume that
$D\OO-D\Sigma$ is also topologically acylindrical and has an incompressible
boundary. Hence, by Thurston's hyperbolization theorem
\cite{\ThuNotes,2,3,4,5,6, \McM1,2, \OtaOne,2}
$D\OO-D\Sigma$ has a complete hyperbolic structure of finite
volume, and we can apply Theorem 2 to $D\OO$.

By Theorem 2, $D\OO$ contains a non-empty compact essential 3-suborbifold 
$\OO'\subset D\OO$ which is not a product and which is either 
Euclidean, Seifert fibered, $Sol$ or complete hyperbolic with finite volume. We distinguish 
two cases, according to whether $\OO$ is closed or not.

If $\partial\OO=\emptyset$ , then
$\OO=\OO'$, because 
the boundary $\partial\OO'$ is either empty or a union of
incompressible toric 2-suborbifolds, and $\OO$ is topologically atoroidal.
Thus $\OO$ is either Seifert fibered, Euclidean or hyperbolic; it cannot be
 $Sol$ by atoroidality.

Next we suppose $\partial\OO\neq\emptyset$.
Note that in this case $\OO'$ cannot be $Sol$ by Lemma 4.3 (ii).
By Lemma 4.1, every component of $\partial\OO'$ is isotopic to
a Euclidean component of $\partial\OO$.
Therefore $\OO'$ is obtained by cutting open $D\OO$ along 
some (perhaps none) components of $\partial\OO$. This implies
that $\OO$ can be isotoped into $\OO'$, because $\OO'$ is connected and not a product.
 Moreover   after isotopy
we can assume that
either
$\tau(\OO')=\OO'$ or $\tau(\OO')\cap\OO'=\emptyset$,
where  $\tau:D\OO\to D\OO$ is the
reflection through $\partial\OO$. 
There
are three possibilities:
\roster
\item"-"
If $\tau(\OO')\cap\OO'=\emptyset$ then $\OO=\OO'$ is
Euclidean, Seifert or hyperbolic, possibly with cusps.
\item"-"
If $\tau(\OO')=\OO'$ and $\OO'$ is hyperbolic, 
then, by atoroidality, $\partial\OO$ has hyperbolic
components and
$\OO'$ is the double of
$\partial\OO$ along the hyperbolic boundary components.
By Lemma 4.6,
$\OO$ is hyperbolic, with some boundary components totally
geodesic and possibly some boundary components cusped.
\item"-"
If $\tau(\OO')=\OO'$ and $\OO'$ is Euclidean or
Seifert fibered, then $\OO'$ is the double of $\OO$ along some
boundary components.
Lemma 4.3 (i) implies that $\OO$ is Seifert fibered, Euclidean or an
$I$-bundle over a 2-orbifold. The
$I$-bundle case is not possible, because it would imply
that $D\OO-D\Sigma$ is Seifert. Hence $\OO$ is Seifert
fibered or Euclidean.
\endroster

This finishes the proof of Theorem 1. \qed
\enddemo

\leftheadtext{ }
\newpage

\

\centerline{\smc chapter \  iii} 

\

\centerline{\chapt A \,  COMPACTNESS  \, THEOREM \,  FOR \,  CONE \,  3-MANIFOLDS}

\vglue.2cm

\centerline{\chapt WITH  \, CONE  \, ANGLES \,  BOUNDED \,  ABOVE \,  BY \,  $\pi$}

\

\

\leftheadtext{III \qquad    Compactness theorem}
\rightheadtext{III \qquad    Compactness theorem }

The purpose of this chapter  is to establish a version of the Gromov
compactness theorem for sequences of Riemaniann  manifolds (cf.
\cite{\GLP} and \cite{\Pe}) in the context of cone 3-manifolds.

Before stating the main theorem we need some definitions.

\definition{Definition} For $\varepsilon\geq 0$, a map $f\!:\! X\to Y$
between two metric spaces is \it $(1+\varepsilon)$-bilipschitz \rm if:
$$
\forall x_1,x_2\in X, \qquad (1+\varepsilon)^{-1} d(x_1,x_2)
\leq d(f(x_1),f(x_2)) \leq (1+\varepsilon) d(x_1,x_2).
$$
\enddefinition 

\remark{Remark} A $(1+\varepsilon)$-bilipischitz map is always an
embedding. Hence one can also define a 
$(1+\varepsilon)$-bilipischitz map as an embedding $f$ such that $f$
and $f^{-1}$ have Lipschitz constant $1+\varepsilon$.
 A map is $1$-bilipischitz if and only if it is an
isometric embedding.
\endremark

\definition{Definition} A sequence of pointed cone 3-manifolds
$\{(C_n,x_n)\}_{n\in\N}$ \it converges geometrically \rm to a pointed
cone 3-manifold $(C_{\infty},x_{\infty})$ if, for every $R>0$ and
$\varepsilon>0$, there exists an integer $n_0$ such that, for
$n>n_0$, there is a $(1+\varepsilon)$-bilipschitz map
$f_n:B(x_{\infty},R)\to C_n$  satisfying:
\roster
\item"(i)" $d(f_n(x_{\infty}),x_n)<\varepsilon$,
\item"(ii)" $B(x_n,R-\varepsilon)\subset f_n(B(x_{\infty},R))$, and
\item"(iii)" $f_n(B(x_{\infty},R)\cap\Sigma_{\infty})=\left(
f_n(B(x_{\infty},R))\right)\cap\Sigma_n$.
\endroster
\enddefinition

\remark{Remark} By definition, the following
inclusion is also satisfied:
$$
f_n(B(x_{\infty},R))\subset B(x_n,R(1+\varepsilon)+\varepsilon).
$$
\endremark

\definition{Definition} For a cone 3-manifold $C$, we define the
\it
cone-injectivity radius
\rm at
$x\in C$:
$$
\Inj(x)=\sup\{\delta>0\text{ such that } B(x,\delta)\text{ is
contained in a standard ball in }C\}.
 $$ We recall that a standard ball is isometric  to either a
non-singular  metric ball in $\Bbb H_K^3$, or to
a singular metric ball in $\Bbb H^3_K(\alpha)$.
The definition does not assume 
the ball to be centered at
$x$, in order to avoid cone-injectivity radius close to zero
for non-singular points close to the singular locus.

Given $a>0$ and $\omega\in (0,\pi]$,
$\CC_{[\omega,\pi],a}$  is the set of pointed cone  3-manifolds
$(C,x)$ with constant curvature in $[-1,0]$, cone angles in
$[\omega,\pi]$, and such that $\Inj (x)\geq a$.
\enddefinition

This chapter is devoted to the proof of the following result:

\proclaim{Theorem \rm (Compactness Theorem)} For $a>0$ and $\omega\in
(0,\pi]$, the closure of  $\CC_{[\omega,\pi],a}$ in
$\bigcup\limits_{b>0}\CC_{[\omega,\pi],b}$ is compact for  the
geometric convergence topology. \endproclaim

This theorem says that any sequence $\{(C_n,x_n)\}_{n\in\N}$ of
pointed cone 3-manifolds in $\CC_{[\omega,\pi],a}$ admits a 
subsequence that converges geometrically to a pointed cone 3-manifold 
in
$\CC_{[\omega,\pi],b}$ for some $b>0$.

The proof of the Compactness Theorem occupies Sections 1 to 4. In
Section 5 we give some properties of the geometric convergence.

The main steps in the proof of the Compactness Theorem are the
following ones.
First we show that $\CC_{[\omega,\pi],a}$ is relatively compact in
the space $\LL$ of locally compact metric length  spaces equipped 
with the Hausdorff-Gromov topology (Proposition 2.1).
 Next, in proposition
2.3, we show that,   for a sequence of pointed cone 3-manifolds of
$\CC_{[\omega,\pi],a}$ that converges in $\LL$, the limit is a
pointed cone 3-manifold in  $\CC_{[\omega,\pi],b}$ for some $b>0$.
Finally we show that Hausdorff-Gromov convergence in
$\CC_{[\omega,\pi],a}$ implies geometric convergence  (Proposition
3.1).

The second and third step of the proof rely on the following technical
result (Proposition 2.2): given a radius $R>0$,
and constants $a>0$ and $\omega\in (0,\pi]$,  for any pointed cone
3-manifold
$(C,x)\in\CC_{[\omega,\pi],a}$, the cone-injectivity radius of each
point in the ball $B(x,R)$ has a positive uniform lower bound, which
only depends on the constants $R$, $a$ and $\omega$. The proof of
this result is postponed untill section 4.

This chapter is organized as follows. Section 1 is
devoted to the Dirichlet polyhedron and the Bishop-Gromov Inequality.
In Section 2 we show that
every sequence in 
$\CC_{[\omega,\pi],a}$ has a subsequence that converges to a cone
3-manifold for the Hausdorff-Gro\-mov topology, assuming Proposition 2.2.
In Section 3 we show that the convergence is in fact geometric.  In
Section 4 we prove Proposition 2.2 using the Dirichlet polyhedron.
Finally, in Section 5 we show some basic properties of the geometric
convergence.

\head 1. The Dirichlet polyhedron   \endhead

In this section we first describe the Dirichlet polyhedron and give some
elementary facts about minimizing paths. Then we prove 
Bishop-Gromov inequality. The Dirichlet polyhedron for cone
3-manifolds is also considered in \cite{\Sua} and minimizing paths for
cone 3-manifolds are also studied in \cite{\HT}.

\definition{Definition} Let $C$ be a  cone 3-manifold
of curvature $K\leq 0$ and $x\in C-\Sigma$. We define the \it
Dirichlet polyhedron centered at $x$: \rm 
$$
D_x=\{y\in C-\Sigma\mid\text{there exists a unique minimizing path
between } y\text{ and } x\}.
 $$
\enddefinition

The open set $D_x$ is star shaped with respect to $x$, hence it can be
locally isometrically embedded in the space of constant sectional
curvature $\Bbb H_K^3$ as a star shaped domain. The following
proposition explains why it is called a polyhedron.

\proclaim{Proposition 1.1} The open domain $D_x$ is the interior of a
solid polyhedron $\overline D_x$ of $\Bbb H_K^3$. Moreover the cone
3-manifold $C$ is isometric to the quotient of $\overline D_x$ under
some
face identifications.
 \endproclaim

In order to prove this proposition we need first to understand the
minimizing paths from $x$ to points in $C-D_x$. First we recall a
well known fact about minimizing paths in cone 3-manifolds with cone
angles less than $2\pi$ (cf. \cite{\HT} for a proof).

\proclaim{Lemma 1.2} Let $C$ be a cone
3-manifold with cone angles less than $ 2\pi$ and  singular set
$\Sigma$. Let
$\sigma$ be a  minimizing path between two points in $ C$.
If $\sigma\cap\Sigma\neq\emptyset$, then either $\sigma\subset\Sigma$,
or $\sigma\cap\Sigma$ is one or both of the end-points of $\sigma$.
\qed \endproclaim

Recall that a subset
$A\subset C$ is called \it convex \rm if every  minimizing path
between two points of $A$ is itself contained in $A$. For instance,
a subset with only one point is convex.
 The following lemma will be used in 
the proof of Proposition 1.1 in the case where $A=\{x\}$, but it will
be used more generally in Section 4 of this chapter and in Chapter
IV.

\proclaim{Lemma 1.3} Let $C$ be a cone 3-manifold of non-positive
curvature, $A\subset C$  a convex subset and $y\in C$. Then the following
hold:
 \roster
\item"(i)" There exist a finite number of minimizing paths from
$y$ to $A$.
\item"(ii)" {Minimizing paths to $A$ with origin
 close to $y$ are obtained
by perturbation:} \newline 
for every $\varepsilon>0$ there exists a
neighborhood
$U\subset C$ of $y$ such that, for every $z\in U$ and every 
 minimizing  path $\sigma_z$ from $z$ to $A$, there exists 
a minimizing path $\sigma_y$ from $y$ to $A$ such that 
$\sigma_z\subset\NN_{\varepsilon}(\sigma_y)$, where
$\NN_{\varepsilon}(\sigma_y)$ is the set of points whose distance 
 to $\sigma_y$ is less
than $\varepsilon$.
 \endroster
\endproclaim

\demo{Proof of Lemma 1.3}
We prove (i) by contradiction. We assume that we have an infinite
sequence $\{\sigma_n\}_{n\in\N}$ of different minimizing paths
between $y$ and $A$. Since $\Length(\sigma_n)$ is constant, up to
taking a subsequence, $\{\sigma_n\}_{n\in\N}$ converges to a path
$\sigma_{\infty}$. Let $\varepsilon>0$ be sufficiently small so that
the developping map around  a tubular neighborhood
$\NN_{\varepsilon}(\sigma_{\infty})$ is defined. By considering
developping maps, we have an infinite sequence of different
minimizing paths between one point and a convex subset in $\Bbb H^3_K$
or $\Bbb H^3_K(\alpha)$, and this is not possible when $K\leq 0$.

We also prove (ii) by contradiction: we assume that there is
$\varepsilon>0$ and a sequence of points $z_n\in C$ such that $z_n\to y$
and every
$z_n$ has a minimizing path $\sigma_n$ to $A$ not contained in
$\NN_{\varepsilon}(\sigma_y)$, for any minimizing path $\sigma_y$
between $y$ and $A$. Since $z_n\to y$, there is a
subsequence of 
$\{\sigma_n\}_{n\in\N}$ that converges to a path
$\sigma_{\infty}$. Moreover, by taking limits in the  inequality
$$
d(y,A)\geq \Length(\sigma_n)-d(y,z_n),
$$
we get that $d(y,A)\geq\Length(\sigma_{\infty})$. Therefore
$\sigma_{\infty}$ is a minimizing path from $y$ to $A$ whose
$\varepsilon$-tubular neighborhood
$\NN_{\varepsilon}(\sigma_{\infty})$ contains infinitely many
$\sigma_n$, and we get a contradiction.
 \qed \enddemo

\demo{Proof of Proposition 1.1}  We describe locally $C-D_x$ by using
Lemma 1.3. Let $y\in C-D_x$, we consider six different cases.

\it Case 1. \rm Consider first the case where $y\not\in\Sigma$ and
there are preciselly two minimizing paths $\sigma_1$ and
$\sigma_2$ between $y$ and $x$. Take $\varepsilon>0$ so that the
developping map is defined around the $\varepsilon$-neighborhoods
$\NN_{\varepsilon}(\sigma_1)$ and $\NN_{\varepsilon}(\sigma_2)$.
By 
Lemma 1.3. there is an open neighborhood $y\in U_y\subset C$ so 
that all minimizing paths from points $z\in U_y$ to $x$ are in  one
of these tubular neighborhoods 
$\NN_{\varepsilon}(\sigma_1)$ or $\NN_{\varepsilon}(\sigma_2)$. 
Since the set of points in $\Bbb H^3_K$ equidistant from two given
different points is a plane, by using developping maps we conclude
that $(C-D_x)\cap U_y$ is the ``bisector" plane between $\sigma_1$
and $\sigma_2$. This case corresponds to the interior of two faces 
of
$\overline D_x$ identified.

\it Case 2. \rm Next consider the case where $y\not\in\Sigma$  and
there are  $n\geq 3$ minimizing paths $\sigma_1,\ldots,\sigma_n$
from $y$ to $x$ satisfying the following property
$$
\text{ there exists } v\in T_y C,\ v\neq 0\text{ such that }
\langle \sigma_1'(0),v \rangle=
\cdots =\langle \sigma_n'(0),v \rangle,
\tag*$$
where the minimizing paths are parametrized by arc length (in
particular $\norm{\sigma_i'(0)}=1$).
 Property \ttag* means that the vectors
$\sigma_1'(0),\ldots,\sigma_n'(0)$ can be ordered in such a way 
that, if $P_i$ denotes the ``bisector" plane between $\sigma'_i(0)$
and 
$\sigma'_{i+1}(0)$, then the intersection $P_1\cap\cdots\cap P_n$ is
a line  (generated by the vector $v$). See Fig. III.1. When $n=3$
property
\ttag* always holds.

\midinsert
 \centerline{\BoxedEPSF{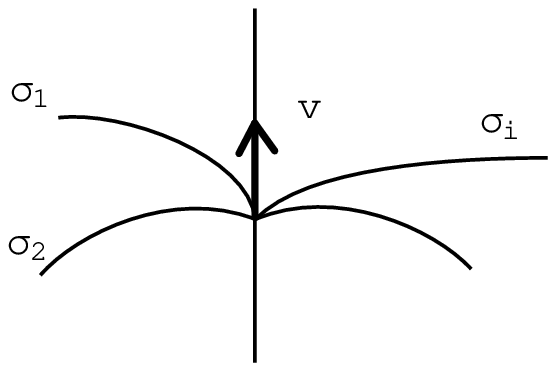 scaled 750}}
   \botcaption{Figure III.1}
    \endcaption
     \endinsert

An argument similar to case 1 shows that, for some neighborhood $U_y$
of $y$, $(C-D_x)\cap U_y$ is the union of $n$ half planes bounded by
the same line. These are precisely the ``bisector" half planes
between the $n$ pairs of paths $\sigma_i$ and $\sigma_{i+1}$. This
case corresponds to the interior of several edges of $\overline D_x$
identified. Note that the dihedral angles are less than
$\pi$ by construction. 

\it Case 3. \rm To finish with the non-singular possibilities, consider
the case where $y\not\in\Sigma$ and there are $n\geq 4$ minimizing
paths  $\sigma_1,\ldots,\sigma_n$ from $y$
to $x$ that do not satisfy  property \ttag* above. This case is
treated as the previous ones and corresponds to $n$ vertices of 
$\overline D_x$
identified.

\it Case 4. \rm When $y\in\Sigma$ and $y$ has only one minimizing
path $\sigma$ to $x$. This case corresponds to the interior of an edge
of  $\overline D_x$, whose dihedral angle equals the cone angle of
$\Sigma$ at $y$. The two adjacent faces of this edges are identified
by a rotation around this edge. 

\it Case 5. \rm When $y\in\Sigma$ and $y$ has $n\geq 2$ minimizing
paths $\sigma_1,\ldots,\sigma_n$ to $x$ that satisfy property \ttag*
of case 2. In this case, the vector $v$ of property \ttag* is
necessarily tangent to $\Sigma$ and it corresponds to $n$ edges of
$\overline D_x$ that are identified to get a piece of $\Sigma$.

\it Case 6. \rm Finaly, consider the case where $y\in\Sigma$ and there
are  $n\geq 2$ minimizing paths $\sigma_1,\ldots,\sigma_n$ between $y$
and $x$ that do not satisfy property \ttag*.  It can be shown that
this case corresponds to   $n$ vertices of 
$\overline D_x$
identified.
  \qed
\enddemo

\proclaim{Corollary 1.4} If the cone angles of $C$ are less than or
equal to $\pi$, then for every $x\in C-\Sigma$ the Dirichlet polyhedron
$D_x$ is convex.
\endproclaim

\demo{Proof} It suffices to show that the dihedral angles of
$\overline D_x$ are less than or equal to $\pi$. We have seen in the
proof of Proposition 1.1 that this is true for dihedral angles of
non-singular edges. For singular edges, this follows from the
hypothesis about cone angles, because dihedral angles are bounded
above by cone angles.\qed \enddemo

We next define the Dirichlet polyhedron centered at singular points.
Recall that  $\Bbb
H^3_K(\alpha)$ denotes the simply connected space of curvature $K$
with a singular axis of cone angle $\alpha$.

\definition{Definition} Let $C$ be a  cone 3-manifold
of curvature $K\leq 0$ and $x\in \Sigma\subset C$. We define the \it
Dirichlet polyhedron centered at $x$: \rm 
$$
D_x=\left\{y\in C\,\left\vert\matrix\format\l\\
\text{there exists a unique minimizing path }
\sigma \text{ between } y\text{ and } x, \\  
\text{and, in addition, if } y\in\Sigma \text { then }
\sigma\subset \Sigma.\endmatrix\right.\right\} 
 $$\enddefinition

As in the non-singular case, $D_x$ is   open,  
star shaped and it can be locally isometrically embedded in $\Bbb
H^3_K(\alpha)$.

\remark{Remark} It is possible to work 
 in the non-singular space $\Bbb H^3_K$ by using the following
construction. Let
$S_{\alpha}$ be an infinite sector of $\Bbb H^3_K$ of
dihedral  angle $\alpha$
and consider the quotient map $p:S_{\alpha}\to
\Bbb H^3_K(\alpha)$ that identifies the faces of $S_{\alpha}$ by a
rotation around its axis. Then we look at the inverse image
$p^{-1}(D_x)\subset S_{\alpha}\subset\Bbb H^3_K$. As in Proposition
1.2, the set $\overline{p^{-1}(D_x)}$ is a solid polyhedron and $C$
is the quotient of $\overline{p^{-1}(D_x)}$ by isometric face
identifications. The point $x$ is in the boundary of 
$\overline{p^{-1}(D_x)}$, although the polyhedron is star shaped with
respect to $x$. As in Corollary 1.4, if the cone angles of $C$ are
bounded above by $\pi$, then $\overline{p^{-1}(D_x)}$ is convex (and
so is $\overline D_x\subset \Bbb H^3_K(\alpha)$).
\endremark

The following lemma will be used in the proof of Lemma 4.4.

\proclaim{Lemma 1.5} Let $C$ be a cone 3-manifold with cone angles
less than or equal to $\pi$. If $x$ is in a compact component
$\Sigma_0$ of $\Sigma$, then $D_x$ is contained in a region of $\Bbb
H^3_K(\alpha)$ bounded by two planes orthogonal to the singular axis
of  $\Bbb H^3_K(\alpha)$. Moreover, the distance between these two
planes is bounded above by $\Length(\Sigma_0)$. \endproclaim

\demo{Proof} By convexity, it suffices to study the  points
$y\in\Sigma_0$ that have at least two minimizing paths to $x$, one of
them contained in $\Sigma_0$. 
When we embed $D_x$ into $
\Bbb H^3_K(\alpha)$, these are the points that will correspond to the
intersection of $\partial \overline D_x$ 
with the singular axis of $
\Bbb H^3_K(\alpha)$. 

 As
in the proof of Proposition 1.1 the local geometry of $\partial \overline 
D_x$ will be given by ``bisector" planes between the minimizing paths form
$y$ to $x$. We distinguish two cases.

First we consider the case where there are precisely two minimizing
paths $\sigma_1,\sigma_2$ between $y$ and $x$ and
$\sigma_1,\sigma_2\subset\Sigma_0$. In this case
$\sigma_1\cup\sigma_2=\Sigma_0$ and the point $y$ is obtained by
identifying the two points of the intersection of $\partial \overline 
D_x$  with the axis of  $
\Bbb H^3_K(\alpha)$. The bisector plane to $\sigma_1$ and $\sigma_2$
passing  through $y$ is the plane orthogonal to $\Sigma_0$, therefore
the lemma is clear in this case.

In the second case, among all the minimizing paths between $y$ and
$x$, one of them $\sigma_1\subset \Sigma_0$ but at least another
$\sigma_2\not\subset\Sigma_0$. Moreover, we may assume that one of
the faces of $\partial\overline  D_x$ is given by the bisector plane
between
$\sigma_1$ and $\sigma_2$.
 In this case, consider the projection $p:S_{\alpha}\to\Bbb H^3_K(\alpha)$
described in the remark above. The preimage $p^{-1}(\sigma_1)$ is
contained in the axis of the sector $S_{\alpha}$, and we choose
the projection $p$ so that $p^{-1}(\sigma_2)$ is contained in the
bisector plane of
$S_{\alpha}$. That is, $p^{-1}(\sigma_2)$ defines the same angle
between both faces of $S_{\alpha}$. Since $\alpha\leq\pi$, the region
of $S_{\alpha}$ bounded by the bisector plane between
$p^{-1}(\sigma_1)$ and $p^{-1}(\sigma_2)$ is contained in the region
of $S_{\alpha}$ bounded by the plane orthogonal to the axis of
$S_{\alpha}$ passing through $p^{-1}(y)\cap p^{-1}(\sigma_2)$. By
convexity, it follows that
$p^{-1}(\overline D_x)$ is contained in the region of $S_{\alpha}$ bounded
by this orthogonal plane. Therefore, $p^{-1}(\overline D_x)\subset
S_{\alpha}$ is contained in a region bounded by two planes orthogonal to
the axis of
$S_{\alpha}$ and the distance of these planes is bounded above by
$\Length{\Sigma_0}$. This finishes the proof of the lemma.
 \qed
\enddemo

Finally, we prove Bishop-Gromov inequality as an application of
the Dirichlet polyhedron.

For $r\geq 0$, let $\V_K(r)$ denote the volume
of the ball of radius $r$ in $\Bbb H^3_K$ the simply connected 
3-space of curvature $K\leq 0$. 

\proclaim{Proposition 1.6. \rm(Bishop-Gromov inequality)} Let $C$ be a
cone 3-manifold of curvature $K\leq 0$ and let $x\in C$. If $0< r\leq R$
then:
$$
\frac{\Vol \left(B(x,r)\right)}{\V_K(r)}\geq \frac{\Vol
\left(B(x,R)\right)}{\V_K(R)}.
$$
\endproclaim

\demo{Proof} The proof follows from the fact that the Dirichlet
polyhedron is star shaped. Namely, if $\Bbb B(x,r)$ is the ball
of radius $r$ in $\Bbb H^3_K$ and $D_x$ is the Dirichlet polyhedron
centered at $x$, then
$
\Vol(B(x,r))=\Vol(\Bbb B(x,r)\cap \overline D_x) \text{ and }
\V_K(r)=\Vol(\Bbb B(x,r)).
$
Since $\overline D_x$ is star shaped, the function 
$
r\mapsto \Vol(\Bbb B(x,r)\cap \overline D_x) /\Vol(\Bbb B(x,r))
$
is decreasing in $r$, and the proposition is proved.
 \qed \enddemo

\proclaim{Corollary 1.7}  Let $C$ be a cone 3-manifold of curvature
$K\in [-1,0]$. Given $\varepsilon>0$ and $R>0$, the number of
disjoint balls of radius $\varepsilon>0$ that can be contained  in a
ball of radius $R$  in
$C$ has a uniform upper-bound, independent of $C$. 
\qed
\endproclaim

\head 2. Hausdorff-Gromov convergence for cone 3-manifolds
\endhead

We first recall some well known definitions.

\definition{Definitions} For  $\varepsilon>0$, an \it
$\varepsilon$-approximation
\rm
 between two pointed compact metric spaces
$(X,x)$ and $(Y,y)$ is a distance $d$ on the disjoint union $X\sqcup Y$
whose restrictions coincide with the original distances on $X$ and
$Y$, and such that $X$ (resp. $Y$) belongs to a
$\varepsilon$-neighborhood of  $Y$ (resp. $X$) and
$d(x,y)\leq\varepsilon$.

 Let $(X,x)$ and $(Y,y)$ be two pointed
compact metric spaces. The \it Hausdorff-Gro\-mov distance \rm
$d_H((X,x),(Y,y))$   is defined as:
$$
d_H((X,x),(Y,y))=\inf\{\varepsilon>0\mid \exists\text{ a }
\varepsilon\text{-approximation between }
(X,x)\text{ and }(Y,y)  \}.
$$
\enddefinition
 
\remark{Remark} By \cite{\GLP, Prop. 3.6}, two pointed compact metric
spaces are isometric by an isometry respecting base points if and
only if their Hausdorff-Gromov distance is zero (see also 
\cite{\BrS}).

Moreover, since $d_H$ verifies the triangle inequality, it is a
distance on the set of pointed compact metric spaces.
\endremark

A cone 3-manifold is a complete metric length space  (cf. Ch. I): the
distance between two points is the infimum of the lengths of paths
joining both points.

In the sequel, $\LL$ will denote the set of complete locally compact
pointed length spaces. Thus we have the inclusion 
$\CC_{[\omega,\pi],a}\subset\LL$.
In a complete locally compact metric length space,  closed
balls are compact (see for instance \cite{\GLP, Thm.
1.10}).  Hence the following definition makes sense.

\definition{Definition}  A sequence $(X_n,x_n)$  in $\LL$ \it
converges for the Hausdorff-Gromov topology \rm to
$(X_{\infty},x_{\infty})\in\LL$ if  for every $R>0$  the Hausdorff-Gromov
distance between the closed balls of radius $R$
$d_H\big(\overline{B(x_n,R)},\overline{B(x_{\infty},R)}\big)$ tends
to zero as $n$ goes to infinity.
\enddefinition

The following proposition is the first step in the proof of the
Compactness Theorem.

\proclaim{Proposition 2.1} The space $\CC_{[\omega,\pi],a}$ is
relatively compact in $\LL$ for the Hausdorff- Gromov topology.
\endproclaim

\demo{Proof}
It is a consequence of Gromov relative compactness criterion
\cite{\GLP, Prop. 5.2} for sequences of pointed complete
locally compact metric
spaces and the fact that the space $\LL$ is closed for the
Hausdorff-Gromov topology \cite{\GLP, Prop. 3.8 and 5.2}.

By Gromov's relative compactness criterion, a sequence
$(X_n,x_n)$ in
$\LL$ has a convergent subsequence if and only if, for every $R>0$ and for every $
\varepsilon>0$, the number of disjoint balls with radius $\varepsilon$
included in the ball $B(x_n,R)$ is uniformly bounded above
independently of $n$. In our case, such a uniform bound follows from
Bishop-Gromov inequality for cone 3-manifolds with constant
curvature  $K\in [-1,0]$ (Corollary 1.7 above, see also
\cite{\HT}).
\qed
\enddemo

We are now stating a key result for the remaining of the proof of
the Compactness Theorem. This result needs the fact that
cone angles are bounded above by $\pi$, and  is not true anymore for
cone angles bigger than $\pi$.

\proclaim{Proposition 2.2 \rm (Uniform lower bound for
cone-injectivity radius)} Given $R>0$, $a>0$ and $\omega\in (0,\pi]$,
there exists a uniform constant $b=b(R,a,\omega)>0$ such that, for every
pointed cone 3-manifold $(C,x)\in\CC_{[\omega,\pi],a}$, the
cone-injectivity radius at any point of $B(x,R)\subset C$ is bigger
than $b$.  \endproclaim

The proof of this proposition is rather long, so we postpone it to
Section 4. We will use it in the proof of the following proposition
as well as in Section 3.

\proclaim{Proposition 2.3} Let $(C_n,x_n)$ be a sequence of pointed
cone 3-manifolds in $\CC_{[\omega,\pi],a}$ that converges to
$(X_{\infty},x_{\infty})$ in $\LL$ for the Hausdorff-Gromov topology.
Then the limit $(X_{\infty},x_{\infty})$ is a pointed cone 3-manifold
in $\CC_{[\omega,\pi],b}$ for some $b>0$. Moreover the curvature of
$X_{\infty}$ is the limit of the curvatures of
$C_n$.
\endproclaim

\remark{Remark} The cone-injectivity radius is lower
semi-continuous, and it could happen that $\Inj(x_{\infty})<a$.
\endremark

\demo{Proof of Proposition 2.3} Let $(C_n,x_n)$ be a sequence in 
$\CC_{[\omega,\pi],a}$ that converges to
$(X_{\infty},x_{\infty})\in \LL$ for the Hausdorff-Gromov topology.
Since $X_{\infty}$ is a complete, locally compact, metric
length space,  we have to show  that $X_{\infty}$ is locally
isometric to a cone 3-manifold of constant sectional curvature.

Let $y\in X_{\infty}$; choose $R=d(y,x_{\infty})+1$. From
Hausdorff-Gromov convergence, for $n$ large enough we have an
$\varepsilon_n$-approximation between  the closed balls
$\overline{B(x_{\infty},R)}$ and 
$\overline{B(x_{n},R)}$, with $\varepsilon_n\to 0$. We take
$y_n\in\overline{B(x_n,R)}$ such that $d_n(y,y_n)<\varepsilon_n$.

By Proposition 2.2, there is a uniform constant $b>0$ independent of
$n$ such that $\Inj(y_n)\geq b$, for every $n$ large enough. Since both
$C_n$ and $X_n$ are length spaces, $d_n$ induces a
$3 \varepsilon_n$-approximation between the compact balls
$\overline{B (y_{n},b)}$ and $\overline{B (y,b)}$. By taking
a subsequence if necessary, there are three cases to be considered:

\it Case 1. \rm For every $n\in\N$, $B(y_n,b)$ is a standard
non-singular ball
 (i.e.  $B(y_n,b)$ is isometric to a metric ball in
the space 
$\Bbb H^3_{K_n}$, where $K_n$ is the curvature of $C_n$).
Since $K_n\in [-1,0]$, up to a subsequence $K_n$ converges to
$K_{\infty}\in [-1,0]$. Moreover, since 
$$
\lim_{n\to\infty}d_H(\overline{B (y_{n},b)},\overline{B (y,b)})
=0,
$$ the unicity of the
Hausdorff-Gromov limit  for compact spaces \cite{\GLP, Prop. 3.6}
shows that the ball
$\overline{B (y_,b)}$ must be isometric to a metric ball in
the space of constant curvature $\Bbb H^3_{K_{\infty}}$.

\it Case 2. \rm For every $n\in\N$,
$B(y_n,b)$ is contained in a standard singular ball, but the distance
between
$y_n$ and $\Sigma$ is bounded below, uniformly away from zero. In
this case, since the cone angles are also bounded below by
$\omega>0$, there exists a uniform constant $b'>0$ such that
$B(y_n,b')\cap\Sigma=\emptyset$ and $B(y_n,b')$ is isometric to a
metric ball in $\Bbb H^3_{K_n}$, for every $n\in\N$. Thus we are in
the first case and we can conclude that $B(y,b')$ is isometric to a
metric ball in $\Bbb H^3_{K_{\infty}}$, where
$K_{\infty}= \lim\limits_{n\to \infty}  K_n$.

\it Case 3. \rm For every $n\in\N$, 
$B(y_n,b)\cap\Sigma\neq\emptyset$
and the distace between $y_n$ and $\Sigma$ tends to zero. In this
case we replace $y_n$ by $y_n'\in\Sigma$ so that $d(y_n,y_n')\to 0$.
The ball $B(y_n',b)$ is isometric to a singular ball
in the space $\Bbb H^3_{K_n}(\alpha_n)$ of constant curvature $K_n$
with a
singular axis, where $K_n$ is the curvature of $C_n$ 
and
$\alpha_n$ the cone angle at $y_n'$. Since $K_n\in[-1,0]$ and
$\alpha_n\in[\omega,\pi]$, up to a subsequence we may assume that
$K_n\to K_{\infty}\in[-1,0]$ and $\alpha_n\to\alpha_{\infty}\in
[\omega,\pi]$.  As in case 1, the fact that the Hausdorff-Gromov
distance $d_H(\overline{B (y_{n}',b)},\overline{B (y,b)})$
tends to zero and the unicity of the Hausdorff-Gromov limit imply
that the ball $B(y,b)$ is isometric to a singular metric ball in 
$\Bbb H^3_{K_{\infty}}(\alpha_{\infty})$ and $y\in\Sigma_{\infty}$.

This achieves the proof of Proposition 2.3 \qed
\enddemo

The following corollary is a direct consequence of the proof of
Proposition 2.3 and will be used later in the proof of Proposition
3.1.

\proclaim{Corollary 2.4}
Let $(C_n,x_n)$ be a sequence of pointed cone 3-manifolds that
converges to the pointed cone 3-manifold $(C_{\infty},x_{\infty})$
for the Hausdorff-Gromov topology. Given  $y\in C_{\infty}$, choose
$R\geq d(x_{\infty},y)+1$ and $d_n$  an $\varepsilon_n$-approximation
between $\overline{B(x_n,R)}$ and $\overline{B(x_{\infty},R)}$, with
$\varepsilon_n\to 0$. Then  there  exists a sequence
$y_n\in B(x_n,R)$ such that $d_n(y_n,y)\to 0$ as $n\to\infty$, and
$y_n\in\Sigma$ if and only if
$y\in\Sigma$. Moreover, when $y\in\Sigma$, the sequence of cone angles
at $y_n$ converges to the cone angle at $y$. \qed
\endproclaim

\head 3. Hausdorff-Gromov convergence implies geometric convergence
\endhead

The goal of this section is to prove the following:

\proclaim{Proposition 3.1} If a sequence $(C_n,x_n)$ of pointed cone
3-manifolds converges in $\CC_{[\omega,\pi],a}$ to a pointed cone
3-manifold $(C_{\infty},x_{\infty})\in\CC_{[\omega,\pi],a}$ for the
Hausdorff-Gromov topology, then it also converges geometrically.
\endproclaim


\demo{Proof of Proposition 3.1} Fix  a radius $R>0$. Let $T$ be a
compact triangulated subset of the underlying space of $C_{\infty}$
such that $B(x_{\infty},12 R)\subset T$. By subdividing the
triangulation we may assume that:
\roster
\item"(i)" all simplices are totally geodesic,
\item"(ii)" $\Sigma_{\infty}\cap T$ belongs to the $1$-skeleton
$T^{(1)}$, and
\item"(iii)" the base point $x_{\infty}$ is a vertex of $T^{(0)}$.
\endroster
Let $R'>0$ be such that $T\subset B(x_{\infty},R')$ and let $d_n$ be a
$\varepsilon_n$-approximation between the compact balls 
$\overline{B(x_{\infty},R')}$ and $\overline{B(x_{n},R')}$, with 
$\lim\limits_{n\to\infty}\varepsilon_n=0$.

Let $T^{(0)}=\{ z^0_{\infty},\ldots, z^r_{\infty}\}$ be the vertices
of $T$, with $z^0_{\infty}=x_{\infty}$. We choose some
points
$ z^0_{n},\ldots, z^r_{n}\in B(x_n,R')$ such that
$\lim\limits_{n\to\infty} d_n(z^i_n,z^i_{\infty})=0$, for
$i=0,\ldots,r$. It follows from Corollary 2.4  that one can choose
$z^i_n\in\Sigma_n$ if and only if $z^i_{\infty}\in\Sigma_{\infty}$.

For a simplex $\Delta$ of $T$, $\Star{\Delta}$ denotes the star of
$\Delta$, and $\Stard{\Delta}$, the union of simplices of $T$ that
intersect $\Delta$ but not the singular set $\Sigma$. With this
notation we have the following:

\proclaim{Lemma 3.2} It is possible to geodesically subdivide  the
triangulation
$T$ so that any simplex $\Delta$ satisfies the following properties:
\roster
\item"(i)" $\Star\Delta$ is included in a standard ball of
$C_{\infty}$.
\item"(ii)" Let $\{ z^{i_1}_{\infty},\ldots,z^{i_s}_{\infty}\}$ be
the vertices of $\Star{\Delta}$. For $n$ sufficiently large, 
$\{ z^{i_1}_{n},$ $\ldots, z^{i_s}_{n}\}$ belongs to a standar ball in
$C_n$,
\item"(iii)" If $\Delta\cap\Sigma_{\infty}=\emptyset$ then
$\Stard\Delta$ is included in a non-singular standard ball.
\item"(iv)" If $\Delta\cap\Sigma_{\infty}=\emptyset$ and 
$\{ z^{i_1}_{\infty},\ldots,z^{i_t}_{\infty}\}$ are the vertices of
$\Stard\Delta$ then, for  $n$ sufficiently large,
$\{ z^{i_1}_{n},\ldots,z^{i_t}_{n}\}$ belongs to a non-singular
standard ball in $C_n$.
 \endroster
\endproclaim

\remark{Remark} It is worthwile to recall that a standard ball in a
cone 3-manifold $C$ with constant sectional curvature $K$ is
isometric either to a non-singular metric ball in $\Bbb H_K^3$, or to
a singular metric ball in $\Bbb H^3_K(\alpha)$ whose center lies in the singular axis.
\endremark

\demo{Proof of Lemma 3.2}
From  Proposition 2.2, there are two
constants $r_1>0$ and $r_2>0$ such that for any $y\in
B(x_{\infty},R')$ or $y_n\in B(x_n,R')$:
\roster
\item"(a)" If $y\in\Sigma_{\infty}$, then $B(y,r_1)$
is a standard singular ball in $C_{\infty}$, and if  $y_n\in\Sigma_n$,
then $B(y_n,r_1)$
is a standard singular ball in $C_{n}$.
\item"(b)" If $d(y,\Sigma_{\infty})>\frac{r_1}8$, then 
$B(y,r_2)$ is a non-singular standard ball 
in $C_{\infty}$, and if $d(y_n,\Sigma_{n})>\frac{r_1}8$,
then $B(y_n,r_2)$ is a non-singular standard ball 
 $C_n$.
\endroster

Using geodesic barycentric subdivision, we first achieve that for each
simplex
$\Delta$ of $T$, the diameter
$\Diam(\Delta)\leq\frac18\inf\{r_1,r_2\}$.  Thus
$\Diam(\Star\Delta)\leq\frac38\inf\{r_1,r_2\}$, which implies
assertions (i) and (ii).

To prove (iii) and (iv) we introduce the constant $r_3(T)>0$
depending on $T$:
$$
r_3(T)=\inf\{ d(\Delta,\Sigma_{\infty})\text{, for any simplex }
\Delta\subset T \text{ such that }
\Delta\cap\Sigma_{\infty}=\emptyset\}.
 $$

Proposition 2.2 and the fact that cone angles are bounded below by
$\omega>0$ imply the existence of a constant $r_4=r_4(r_3,\omega)>0$ 
depending on $r_3$ and $\omega$ such that for any point $y\in B(x_{\infty},R')$, if
$d(y,\Sigma_{\infty})>r_3/2$ then $B(y,r_4)$ is a non-singular
standard ball in $C_{\infty}$, and for any point $y_n\in
B(x_{n},R')$, if $d(y_n,\Sigma_{n})>r_3/2$ then $B(y_n,r_4)$ is
a non-singular standard ball in $C_{n}$.

Next we will subdivide $T$ in such a way that the constants $r_3$ and
hence $r_4=r_4(r_3,\omega)$ do not change, but the diameter of any simplex
not meeting $\Sigma_{\infty}$ becomes less than $\frac{r_4}8$. This
will imply properties (iii) and (iv) of the lemma. 

The process for subdividing a simplex $\Delta$ of $T$ is the
following:
\roster
\item"(a)" When $\Delta\cap\Sigma_{\infty}=\emptyset$, we apply
geodesic baricentric subdivision to $\Delta$.
\item"(b)" When $\Delta\subset\Sigma_{\infty}$, we do not subdivide
$\Delta$.
\item"(c)" When $\emptyset\neq\Delta\cap\Sigma_{\infty}\neq\Delta$,
we express $\Delta$ as a joint $\Delta=\Delta_0*\Delta_1$, where
$\Delta_0\subset\Sigma_{\infty}$ and
$\Delta_1\cap\Sigma_{\infty}=\emptyset$. Then we apply geodesic
barycentric subdivision $\Delta_1'$ to $\Delta_1$ and consider
$\Delta'=\Delta_0*\Delta_1'$ (see Figures III.2 to III.4, which
describe this process).
\endroster
This process of subdivision makes the diameter of any simplex
disjoint from $\Sigma_{\infty}$ arbitrarily small without
decreasing
its distance to $\Sigma_{\infty}$. \qed
 \enddemo

\midinsert
 \centerline{\BoxedEPSF{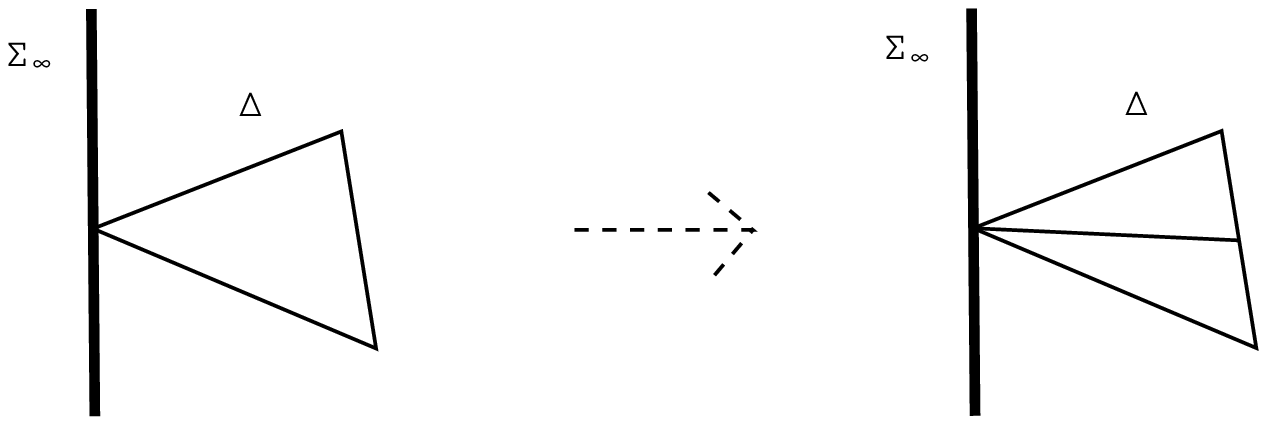 scaled 750}} 
   \botcaption{Figure III.2}
    \endcaption
     \endinsert

\midinsert
 \centerline{\BoxedEPSF{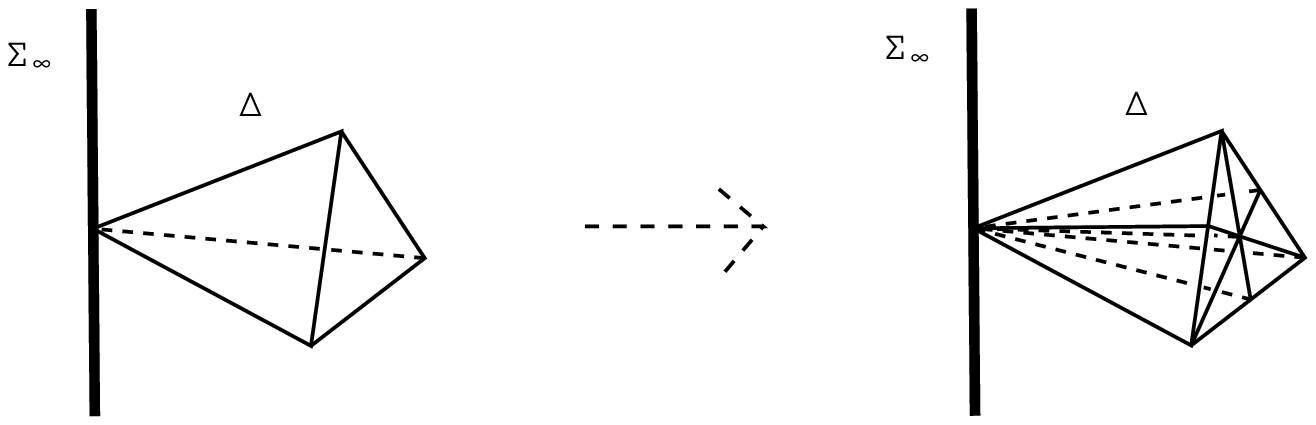 scaled 750}}
   \botcaption{Figure III.3}
    \endcaption
     \endinsert

\midinsert
 \centerline{\BoxedEPSF{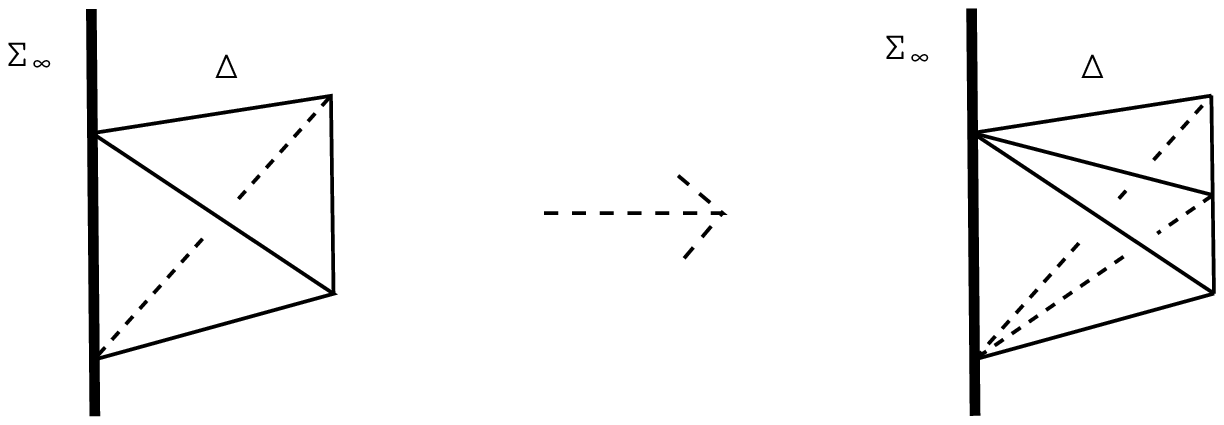 scaled 750}}
   \botcaption{Figure III.4}
    \endcaption
     \endinsert

Now we define $g_n:T\to C_n$ by mapping $z^i_{\infty}$ to $z^i_n$,
$g_n(z^i_{\infty})=z^i_n$ for $i=1,\ldots,r$, and  we extend
$g_n$ piecewise-linearly on each simplex of $T$. To show that the
restriction of $g_n$ to $\overline{B(x_{\infty},R)}$ is a
$(1+\varepsilon_n)$-bilipischitz map with $\varepsilon_n\to 0$,
we need the following lemma:

\proclaim{Lemma 3.3} For $n$ sufficiently large, $g_n:T\to C_n$ is a
well defined map having the following properties:
\roster
\item"(i)" $g_n(T)$ is a  geodesic polyhedron in $C_n$, and
$g_n(T\cap\Sigma_{\infty})=g_n(T)\cap\Sigma_n$;
\item"(ii)" $\forall x,y\in B(x_{\infty},6 R)$,
$d(g_n(x),g_n(y))\leq(1+\delta_n) d(x,y)$ with $\delta_n\to 0$;
\item"(iii)" the restriction of $g_n$ induces a homeomorphism from
$\Int(T)$ onto its image 
$g_n(\Int(T))$.
 \endroster
\endproclaim

\demo{Proof of Lemma 3.3}
Let $\Delta$ be a 3-simplex in $T$ such that
$\Delta\cap\Sigma_{\infty}=\emptyset$. Up to permutation of indices,
let $\{z^1_{\infty},z^2_{\infty},z^3_{\infty},z^4_{\infty}\}$  denote
the vertices of $\Delta$. By Lemma 3.2, there exists $n_1$ such that
for $n\geq n_1$, the points 
$\{z^1_{n},z^2_{n},z^3_{n},z^4_{n}\}$ are contained in a non-singular
standard ball in $C_n$. By construction, the sequence
$\{d(z^i_n,z^j_n)\}_{n\in\N}$ tends to $d(z^i_{\infty},z^j_{\infty})$
as $n\to\infty$, for any $i,j\in\{1,2,3,4\}$. Moreover, the sectional
curvature $K_n$ of $C_n$ tends to $K_{\infty}$. It follows that the
the bijection between 
$\{z^1_{\infty},z^2_{\infty},z^3_{\infty},z^4_{\infty}\}$ and
$\{z^1_{n},z^2_{n},z^3_{n},z^4_{n}\}$ extends linearly to a map from
the geodesic simplex $\Delta$ onto the non-degenerated geodesic
simplex in $C_n$ with vertices $\{z^1_{n},z^2_{n},z^3_{n},z^4_{n}\}$.
That is, $g_n(\Delta)$ is a well defined non-degenerated simplex such
that the restriction map $g_n\vert_{\Delta}:\Delta\to g_n(\Delta)$ is
$(1+\delta_n)$-bilipschitz, with $\delta_n\to 0$.

If $\Delta\cap\Sigma_{\infty}\neq\emptyset$ then, by Corollary 1.4, 
for $n$ sufficiently large, $z^i_n\in\Sigma_n$ if and only if
$z^i_{\infty}\in\Sigma_{\infty}$, $i=1,2,3,4$. By using the same
method as in the non-singular case and Lemma 3.2, one shows that $g_n$ is well defined
on
$\Delta$ and that $g_n(\Delta)$ is a non-degenerated totally
geodesic simplex in $C_n$ such that
$g_n(\Delta\cap\Sigma_{\infty})=g_n(\Delta)\cap\Sigma_n$.
We remark that minimizing paths between two points are not necesarily unique in singular
balls, but they are unique if at least one of the points lies in $\Sigma$.
Hence, $g_n$ is well defined on $\Delta$, because of assertions (iii)
and (iv) of Lemma 3.2.
 Moreover the
restriction map $g_n\vert_{\Delta}:\Delta\to g_n(\Delta)$ is
$(1+\delta_n)$-bilipschitz with $\delta_n\to 0$. This proves property
(i).

To prove property (ii), consider $x,y\in B(x_{\infty},6 R)$. Let 
$\sigma$ be a minimizing path between
 $x$ and $y$. Since $\sigma\subset B(x_{\infty},12
R)\subset T$, the inequality $d(g_n(x),g_n(y))\leq  (1+\delta_n)
d(x,y)$ follows from the fact that, for any 3-simplex $\Delta$ of
$T$,
$g_n:\Delta\to g_n(\Delta)$ is a $(1+\delta_n)$-bilipschitz map, 
with
$\delta_n\to 0$. It suffices to break up $\sigma$ into pieces
$\sigma\cap\Delta$ and to use the fact that $\sigma$ is minimizing.

Finally we prove property (iii). Note that the
restriction of $g_n$ from $\Star\Delta$ onto   $g_n(\Star\Delta)$ is
a homeomorphism. This follows from
the construction of $g_n$ by
piecewise-linear extension and the fact that, for $n$ sufficiently
 large and for any simplex $\Delta$ of $T$, $\Star\Delta$ and
$g_n(\Star\Delta)$ are contained in standard balls.
Thus it remains to show that the restricition of $g_n$ to
$\Int(T)$ is injective for $n$ sufficiently large. 

 Suppose that 
$x,y\in\Int(T)$
are two points
 such that $g_n(x)=g_n(y)$;
we claim that $x=y$.
 Let $\Delta_x$ and
$\Delta_y$ be the simplices of $T$ containing $x$ and $y$ 
respectively.
 We claim first that $\Delta_x\cup\Delta_y$ is
contained in a standard ball in $C_{\infty}$. 
Recall that the diameter of the simplices is choosen to be small 
with respect to the lower bound of the cone-injectivity radius on
$T$. Thus we prove the claim by showing  that the diameter of
$\Delta_x\cup\Delta_y$ is also small. 
To show this, we first remark that the diameter of
$g_n(\Delta_x\cup\Delta_y)$ is small, because 
$g_n(\Delta_x)\cap g_n(\Delta_y)\neq\emptyset$ and
$\Diam(g_n(\Delta))\leq(1+\delta_n)\Diam(\Delta)$, with $\delta_n\to 0$.
In particular, $g_n(\Delta_x\cup\Delta_y)$ is contained in a standard
ball. Moreover, the limit
$\lim\limits_{n\to\infty} d(z_n^i,z_n^j)= d(z_{\infty}^i,z_{\infty}^j)
$
means that the distance between vertices of $g_n(\Delta_x\cup\Delta_y)$
converges to the distance between vertices of $\Delta_x\cup\Delta_y$,
therefore $\Diam(\Delta_x\cup\Delta_y)$ is small.

Finally, since the restriction 
$g_n\vert_{\Delta_x\cup\Delta_y}$ is
defined by geodesic linear extension in a standard ball, 
$g_n\vert_{\Delta_x\cup\Delta_y}$ is
injective; thus $x=y$. This proves Lemma 3.3.
\qed
\enddemo

\proclaim{Lemma 3.4} For $n$ sufficienly large, two points in
$g_n(B(x_{\infty},R))$ are joined by a minimizing geodesic contained
in $g_n(B(x_{\infty},5 R))$.
\endproclaim

\demo{Proof} Since two points in $B(z_n^0,2R)\subset C_n$ are joined
by a minimizing geodesic contained in $B(z_n^0,4R)$, it suffices to
show the following inclusions for $n$ large enough:
\roster
\item"(a)" $g_n(B(x_{\infty},R))\subset B(z^0_n,2 R)$,
\item"(b)" $B(z^0_n,4R)\subset g_n(B(x_{\infty},5R))$,
\endroster
where we recall that $z^0_n=g_n(x_{\infty})$ and
$x_{\infty}=z^0_{\infty}\in T^{(0)}$.

Property (a) follows from Lemma 3.3 (ii), which states that for $x,y\in
B(x_{\infty}, 6R)$, $d(g_n(x),g_n(y))\leq (1+\delta_n) d(x,y)$, with
$\delta_n\to 0$. 

For $n$ sufficiently large, the restriction $g_n:B(x_{\infty}, 6 R)
\to g_n(B(x_{\infty}, 6 R))$ is a homeomorphism and hence 
$g_n(\partial B(x_{\infty},5 R))=\partial g_n(B(x_{\infty}, 5R))$.
Thus  the inclusion (b) will follow from the inequality
$d(z^0_n,g_n(\partial B(x_{\infty}, 5R)))>4R$, for $n$ sufficiently
large.

Let $y$ be a point in $\partial B(x_{\infty},5 R)$, that is 
$d(x_{\infty},y)=5 R$. Set: 
$$
r_0=\sup\{\Diam (\Delta)\mid\Delta\text{
is a 3-simplex of }T\}.
$$
Since $y\in\partial B(x_{\infty},5R)\subset T$, there exists a vertex
$z^{i}_{\infty}\in T^{(0)}$ such that $d(z^i_{\infty},y)\leq r_0$.
So we write:
$$
d(z^0_n,g_n(y))\geq d(z^0_n,z^i_n)-d(z^i_n,g_n(y)).
 $$
By Lemma 3.3, $d(z^i_n,g_n(y))\leq (1+\delta_n) d(z^i_{\infty},y)\leq 2
r_0$. Moreover,
$$
\lim_{n\to\infty}d(z^0_n,z^i_n)=d(z^0_{\infty},z^i_{\infty})=
d(x_{\infty},z^i_{\infty})\geq d(x_{\infty},y)- d(z^i_{\infty},y)
\geq 5 R-r_0.
$$
Summarizing these inequalities we conclude that 
$ d(z^0_n,g_n(y))\geq 5 R- 4r_0
$
and it suffices to choose $r_0<R/4$ by using the proof of 
Lemma 3.2. This
achieves the proof of inclusion b) and of Lemma 3.4. \qed 
 \enddemo

The following lemma concludes the proof of proposition 3.1.

\proclaim{Lemma 3.5} For any real $\varepsilon>0$, there is an
integer $n_0$ such that, for $n\geq n_0$:
\roster
\item"(i)" The restriction $g_n\!:\! B(x_{\infty},R)\to C_{n}$ is
$(1+\varepsilon)$-bilipschitz,
\item"(ii)" $d(g_n(x_{\infty}),x_n)<\varepsilon/2$, 
\item"(iii)" $B(x_n,R-\varepsilon)\subset g_n(B(x_{\infty},R))$.
\endroster
\endproclaim

\demo{Proof} By Lemma 3.3 (ii), there exists a sequence
$\delta_n\to 0$ such that
$$ \qquad\qquad
d(g_n(x),g_n(y))\leq (1+\delta_n) d(x,y) \qquad \forall x,y,\in
B(x_{\infty}, 6 R),\forall n\geq n_0.
$$
By choosing $n$ sufficiently large we may assume
$\delta_n<\varepsilon$, hence property (i) will follow from the
following inequality
$$
 \qquad\qquad
(1+\varepsilon)^{-1} d(x,y)\leq d(g_n(x),g_n(y))
, \qquad
\forall n\geq n_0, \forall x,y\in B(x_{\infty},R).
$$
To prove this inequality, given $x,y\in  B(x_{\infty},R)$, we choose a
minimizing path $\sigma$ between $g(x)$ and $g(y)$ that is
contained in $g_n(B(x_{\infty}, 5R))$, by Lemma 3.4. 
Since $g_n:B(x_{\infty},5 R)\to g_n(B(x_{\infty},5R))$ is  a
homeomorphism, $\tilde \sigma=g_n^{-1}(\sigma)$ is a path joining
$x$ and $y$. The map $g_n$ is constructed in the proof of Lemma 3.3 
so that its restriction to each simplex $\Delta$ of $T$
is  $(1+\delta_n)$-bilipschitz. Then, by breaking $\tilde \sigma$ into
pieces $\tilde \sigma\cap \Delta$ we prove that
$(1+\varepsilon)^{-1}\Length(\tilde\sigma)\leq\Length(\sigma)$, and
the claimed inequality follows.

Property (ii) follows from the construction, because the Hausdorff-Gromov
distance between the pointed balls
$\big(\overline{B(x_n,R')},x_n\big)$ and
$\big(\overline{B(x_{\infty},R')},x_{\infty}\big)$ goes to zero,
and the points $x_{\infty}$ and $g_n(x_{\infty})$ are arbitrarily
close in the Hausdorff-Gromov approximations.

Next we prove  property (iii). Let $y_n\in B(x_n,
R-\varepsilon)$. By property (ii), $y_n\in B(g_n(x_{\infty}),R)=
B(z_n^0,R)$. Moreover, in the proof of Lemma 3.4 (inclusion (b))
we have seen that $B(z_n^0,R)\subseteq g_n(B(x_{\infty},5R))$.
Hence we can choose a point $y\in B(x_{\infty},5 R)$  such that $g_n(y)=y_n$.
Then, for $n$ large, we have
$$
d(x_{\infty},y)\leq (1+\delta_n) \big(d(y_n,x_n)+d(x_n,g_n(x_{\infty}))
\big)
\leq (1+\delta_n) (d(x_n,y_n)+\varepsilon/2),
$$
with $\delta_n\to 0$. For $n$ sufficiently large so that
$\delta_n<\varepsilon/(2R)$, we conclude  that $y\in
B(x_{\infty},R)$. Hence $B(x_n,
R-\varepsilon)\subseteq g_n(B(x_{\infty},R))$ and the lemma is proved. \qed
\enddemo
 \enddemo

\head 4. Uniform lower bound for the cone-injectivity radius
\endhead

The goal of this section is to prove Proposition 2.2. For convenience
we recall the statement.

\proclaim{Proposition 2.2 \rm (Uniform lower bound for
cone-injectivity radius)} Given $R>0$, $a>0$ and $\omega\in (0,\pi]$,
there exists a uniform constant $b=b(R,a,\omega)>0$ such that, for any
pointed cone 3-manifold $(C,x)\in\CC_{[\omega,\pi],a}$, the
cone-injectivity radius at any point of $B(x,R)\subset C$ is bigger
than $b$.  \endproclaim

\remark{Remarks} (1) We recall the definition of cone-injectivity
radius at a point $x\in C$:
$$
\Inj(x)=\sup\{\delta>0\text{ such that } B(x,\delta)\text{ is
contained in a standard ball in }C\}.
 $$
Note that the definition does not assume the ball to be centered at
$x$, otherwise regular points near the singular locus would have
arbitrarily small cone-injectivity radius; such points are contained
in larger standard balls centered at nearby singular points.
Moreover, if $x\in\Sigma$ then the standard   ball in the
definition can be assumed to be centered at $x$.

(2) Proposition 2.2 implies that there is a uniform lower bound for
the radius of a tubular neighborhood of the singular locus $\Sigma$.
In particular the singular locus can not cross itself when the cone
angles are $\leq\pi$. The proof of Proposition 2.2 is based on volume
estimates using the convexity of the Dirichlet polyhedron (Corollary
1.4).
\endremark

The proof of Proposition 2.2 is divided in two propositions, the
first one deals with the case of singular points, the second one with
the case of regular points.

\proclaim{Proposition 4.1}
Given $R>0$, $a>0$ and $\omega\in (0,\pi]$, there exist constants
$\delta_1=\delta_1(R,a,\omega)>0$ and
$\delta_2=\delta_2(R,a,\omega)>0$ (depending only on $R$, $a$ and
$\omega$) such that any pointed cone 3-manifold
$(C,x)\in\CC_{[\omega,\pi],a}$ satisfies:
\roster
\item"(i)" any component $\Sigma_0$ of the singular locus
$\Sigma\subset C$ that intersects $B(x,R)$ has length
$\Mod{\Sigma_0}\geq\delta_1$,
\item"(ii)" $\NN_{\delta_2}(\Sigma)\cap B(x,R)= \{ y\in B(x,R)
\mid d(y,\Sigma)<\delta_2\}$ is a tubular neighborhood of
$\Sigma\cap B(x,R)$.
\endroster
\endproclaim

\proclaim{Proposition 4.2}
Given $R>0$, $a>0$ and $\omega\in (0,\pi]$, there exists a constant 
$\delta_3=\delta_3(R,a,\omega)>0$ (depending only on $R$, $a$ and
$\omega$) such that for any pointed cone 3-manifold
$(C,x)\in\CC_{[\omega,\pi],a}$, if $y\in B(x,R)\subset C$ and
$d(y,\Sigma)>\min(\delta_1,\delta_2)$ then $\Inj(y)>\delta_3$ (where
$\delta_1$ and $\delta_2$  are the constants given in Proposition
4.1).
\endproclaim

The proof of Proposition 4.2 is given in \cite{\KoOne, Prop. 5.1.1}.
It may  be proved also by perturbing the singular metric on the
tubular neighborhood of $\Sigma\cap B(x,R)$ with radius
$\min\{\delta_1,\delta_2\}$ to get a Riemannian metric with pinched
sectional curvature, with a pinching constant depending only on
$\delta_1$ and $\delta_2$; then we are in the case of non-singular
Riemannian metrics for which the result is well known (cf.
\cite{\GLP}, \cite{\Pe}). Therefore we only give the proof of
Proposition 4.1. 

\demo{Proof of Proposition 4.1 (i)} The proof follows from the
following volume estimations.

\proclaim{Lemma 4.3} Given $a>0$ and $\omega\in (0,\pi]$,  there
exists a constant $c_1=c_1(a,\omega)>0$ such that for any pointed
cone 3-manifold $(C,x)\in\CC_{[\omega,\pi],a}$, $\Vol(B(x,1))\geq c_1$.
\endproclaim

\proclaim{Lemma 4.4} Given $R>0$ there is a constant $c_2=c_2(R)>0$
such that, if $C$ is a cone 3-manifold of curvature $K\in[-1,0]$ and if
$\Sigma_0$ is a component of the singular locus of $C$, then for any
$y\in\Sigma_0$, $\Vol (B(y,R+1))\leq c_2(R) \Mod{\Sigma_0}$, where 
$\Mod{\Sigma_0}$ is the length of $\Sigma_0$.
\endproclaim

\demo{Proof of Proposition 4.1 (i) from Lemmas 4.3 and 4.4} Let
$\Sigma_0$ be a component of $\Sigma$ that intersects $B(x,R)$ and
$y\in\Sigma_0\cap B(x,R)$. By Lemmas 4.3 and 4.4 we have:
$$
c_1=
c_1(a,\omega)\leq \Vol(B(x,1))\leq \Vol(B(y,R+1)) \leq c_2(R)
\Mod{\Sigma_0}.
$$
Therefore $\Mod{\Sigma_0}\geq \delta_1=c_1/c_2$.\qed
\enddemo

We now give the proofs of Lemmas 4.3 and 4.4.

\demo{Proof of Lemma 4.3}  Let $(C,x)\in\CC_{[\omega,\pi],a}$; in
particular $\Inj(x)>a$.  Because of the definition of the
cone-injectivity radius, we distinguish two cases,
according to wehether $B(x,a)$ is contained in a
singular standard ball or in a  non-singular one.

\it Non-singular case. \rm When $B(x,a)$ is a non-singular standar
ball,  by taking $a_0=\inf\{1,a\}$ we have $\Vol(B(x,1))\geq \Vol
(B(x,a_0))\geq \frac43\pi  a_0^3$, because the curvature $K\leq 0$.

\it Singular case. \rm When $B(x,a)$ is contained in a standard
singular ball, there exists a point $z\in\Sigma$ and $a'\geq a$ such
that $B(z,a')$ is a singular standard ball that contains $B(x,a)$. We
may assume that $d(x,z)=d(x,\Sigma)$.  We distinguish again two
subcases.

If $d(x,z)\leq 1/2$, then by taking $a_0=\inf\{1/2,a\}$, we have
that $B(z,a_0)\subset B(x,1)$ and thus
$\Vol(B(x,1))\geq\Vol(B(z,a_0))\geq\frac23\omega a_0^3$, because  the
cone angles are bounded below by $\omega$ and  the curvature $K\leq 0$.
 
If $d(x,z)=d(x,\Sigma)\geq 1/2$, an elementary trigonometric argument
shows that we can find  a constant $b=b(\omega,a)>0$ 
such that $B(x,b)$ is a non-singular standar ball. This constant
depends only on $\omega$ and $a$, because the curvature $K\in [-1,0]$. As
in the non-singular case, by taking
$b_0=\inf\{b,1\}$, we have the inclusion $B(x,b_0)\subset B(x,1)$ and
the inequality $\Vol( B(x,1))\geq \frac43\pi b_0^3$.

This finishes the proof of Lemma 4.3. \qed
\enddemo

\demo{Proof of Lemma 4.4} Let $\Sigma_0$ and $y\in\Sigma_0$ be as in
the statement of Lemma 4.4. Consider $D_y$  the Dirichlet polyhedron
centered at $y$. By Lemma 1.5, $D_y$ is contained in the region  of
$\Bbb H_K^3(\alpha)$ bounded by two planes orthogonal to the
singular axis of $\Bbb H_K^3(\alpha)$ and the
distance between them is less that or equal to the length
$\Mod{\Sigma_0}$. Therefore we have:
$$
\Vol(D_y\cap B(y,R+1))\leq 2\pi\Mod{\Sigma_0}\sinh_K^2(R+1)
$$
where
$\sinh_K(r)=\sinh(\sqrt{-K}r)/\sqrt{-K}$ if $K<0$ and $\sinh_0(r)=r$.
Since $K\in[-1,0]$, $\sinh_K(r)\leq\sinh(r)$ and we conclude:
$$
\Vol(D_y\cap B(y,R+1))\leq \left(2\pi\sinh^2(R+1)\right)
\Mod{\Sigma_0}. $$
This inequality proves Lemma 4.4.\qed
\enddemo

\demo{Proof of Proposition 4.1 ii)}
Let $\sigma$ be a minimizing arc between two points of
$B(x,R)\cap\Sigma$ which is not contained in $\Sigma$, in particular 
$\sigma\cap\Sigma=\partial\sigma$. We assume that $\sigma$ has
minimal length among all such possible arcs. Proposition
4.1 ii) will follow from Lemma 4.3 and the following:

\proclaim{Lemma 4.5} There exists a constant $c_3=c_3(R)>0$ depending
only on $R$ such that $\Vol(\NN_{R+1}(\sigma))<c_3(R)\Mod{\sigma}$,
where $\NN_{R+1}(\sigma)=\{ y\in C\mid d(y,\sigma)\leq R+1\}$ and
$\Mod\sigma$ is the length of $\sigma$.
\endproclaim

The proof of Proposition 4.1 ii) from Lemmas 4.3 and 4.5 is similar to
the proof of Propostition 4.1 i). From the inclusion $B(x,1)\subset
\NN_{R+1}(\sigma)$ and the inequalities of Lemmas 4.3 and 4.5 we conclude
that $\Mod\sigma>c_1/c_3$. Thus it suffices to choose
$\delta_2=\frac12 c_1/c_3$ in Proposition 4.1 ii). 

 The remaining of this section is
devoted to the proof of Lemma 4.5.

\demo{Proof of Lemma 4.5}
Let $D_{\sigma}$ be the open subset of $C$ defined as:
$$
D_{\sigma}=\{ y\in C-\Sigma\mid\text{ there is a unique minimizing
arc between } y\text{ and }\sigma\}.
 $$
The open set $D_{\sigma}$ is perhaps not convex, but it is star
shaped with respect to $\sigma$. So $D_{\sigma}$ may be isometrically
embedded in $\Bbb H^3_K$, the space of constant sectional curvature
$K\in[-1,0]$.

\proclaim{Claim 4.6} The set $C-D_{\sigma}$ has Lebesgue
measure zero.
\endproclaim

\demo{Proof} Since $\Sigma$ is $1$-dimensional, it suffices to show
that $C-(\Sigma\cup D_{\sigma})$ has measure zero.  Given $z\in 
C-(\Sigma\cup D_{\sigma})$ there are several  but only a finite
number of minimizing paths from $z$ to $\sigma$, by Lemma 1.2.
Moreover, by the same lemma, there is a neighborhood $U_z$ of $z$ such
that for every $y\in U_z\cap C-(\Sigma\cup D_{\sigma})$ the minimizing
paths between $y$ and $\sigma$ are in tubular neighborhoods of the
minimizing paths between $z$ and $\sigma$. Therefore, by using
developping maps along these tubular neighborhoods and the fact that
the set of points in $\Bbb H_K^3$ that are equidistant from two
geodesics has measure zero, we conclude that $U_z\cap C-(\Sigma\cup
D_{\sigma})$ has measure zero. This implies in particular that
$C-(\Sigma\cup D_{\sigma})$ itself  has measure zero and the claim is
proved.\qed \enddemo

This claim implies that $\Vol(\NN_{R+1}(\sigma))=\Vol(D_{\sigma}\cap
\NN_{R+1}(\Sigma))$. We next use the fact that $D_{\sigma}\cap
\NN_{R+1}(\Sigma)$ can be isometrically embedded in $\Bbb H^3_K$ to
get an upper bound for its volume.

Let $\{p,q\}=\sigma\cap\Sigma$ be the end-points of $\sigma$. The
proof of Lemma 4.5 is divided in three cases: Lemmas 4.7, 4.8 and
4.9.

\proclaim{Lemma 4.7} If $\sigma$ is orthogonal to $\Sigma$ at $p$ and
$q$, then Lemma 4.5 holds.
\endproclaim

\demo{Proof} Since the cone angles of $C$ are $\leq \pi$, the
orthogonality hypothesis implies that $D_{\sigma}$ is contained in the
subspace of $\Bbb H^3_K$ bounded by the two planes orthogonal to
$\sigma$ at its end-points $p$ and $q$. Therefore, as in the proof of
Lemma 4.4, we obtain the following inequality: $$
\Vol(\NN_{R+1}(\Sigma)\cap D_{\sigma})\leq
2\pi\Mod\sigma\sinh_K^2(R+1)\leq 2\pi\sinh^2(R+1)\Mod\sigma
,
 $$
because $K\in [-1,0]$. \qed
\enddemo

It may happen that $\sigma$ is not orthogonal to $\Sigma$ at $p$ or
$q$, when $p$ or $q$ belong to the boundary of $B(x,R)$. Let
$\theta$ and $\phi$ in $[0,\pi/2]$ be the angles between $\sigma$ and
$\Sigma$ at $p$ and $q$ respectively.


\proclaim{Lemma 4.8} 
If
$\max\{\cos(\theta),\cos(\phi)\}
\leq2\Mod\sigma$ then Lemma 4.5 is true.
\endproclaim

\demo{Proof} By assumption, $D_{\sigma}$ is contained in the union
$S_{\sigma}\cup S_p\cup S_q\subset \Bbb H^3_K$, where $S_{\sigma}$ is
the subspace of $\Bbb H^3_K$ bounded by the two planes orthogonal to
$\sigma$ at $p$ and $q$, $S_p$ is a solid angular sector with axis
passing through $p$ and dihedral angle $\pi/2-\theta$ at the axis, and
$S_q$ is a solid angular sector with axis passing through $q$ and
dihedral angle $\pi/2-\phi$ at the axis. One of the faces of $S_p$
(and of $S_q$) is a face of $S_{\sigma}$ and the other contains a
piece of $\Sigma$ (cf. Fig. III.5).

\midinsert
 \centerline{\BoxedEPSF{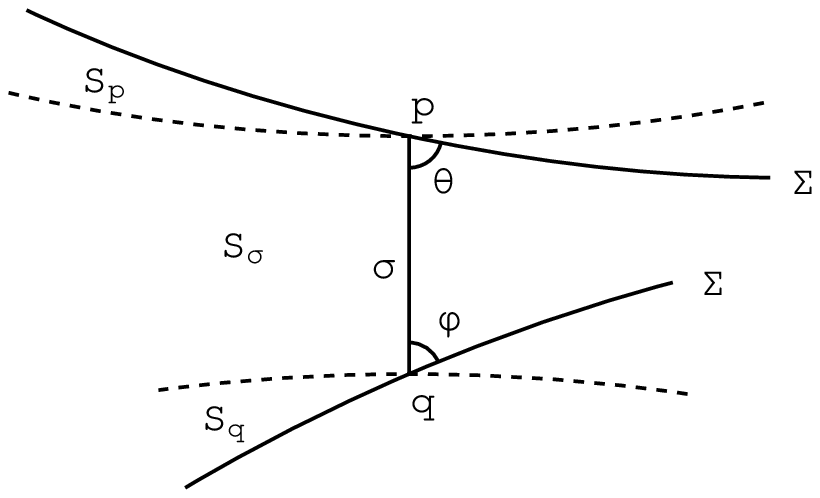 scaled 750}}
   \botcaption{Figure III.5}
    \endcaption
     \endinsert

Since $$\Vol(\NN_{R+1}(\Sigma)\cap D_{\sigma})
=  \Vol(\NN_{R+1}(\Sigma)\cap S_{\sigma})+
 \Vol(\NN_{R+1}(\Sigma)\cap S_{p})+
 \Vol(\NN_{R+1}(\Sigma)\cap S_{q}),$$ to prove Lemma 4.5 it suffices
to get a suitable upper bound for each one of these three volumes.

For $
\Vol(\NN_{R+1}(\Sigma)\cap S_{\sigma})$ the same upper bound as in
Lemma 4.7 goes through: 
$$
\Vol(\NN_{R+1}(\Sigma)\cap S_{\sigma})\leq
2\pi\sinh^2(R+1)\Mod{\sigma}.
$$

For $\Vol(\NN_{R+1}(\Sigma)\cap S_{p})$  we use the volume of the
sector of angle $\pi/2-\theta$:
$$
\Vol(\NN_{R+1}(\Sigma)\cap S_{p})\leq\frac{\pi/2-\theta}{2\pi}
\V_{K}(R+1),
$$
where $\V_{K}(R+1)$ is the volume of the ball
of radius $R+1$ in $\Bbb H^3_K$.
 Moreover, since $K\geq-1$, $\V_{K}(R+1)
\leq \V_{-1}(R+1)\leq\pi\sinh(2R+2)$. Therefore:
 $$
\Vol(\NN_{R+1}(\Sigma)\cap S_{p})\leq(\pi/2-\theta)\frac12\sinh(2R+2).
$$
Since
$\lim\limits_{\theta\to\pi/2}\cos(\theta)/(\pi/2-\theta)=1$, there is a
constant $\lambda>0$ such that  $\pi/2-\theta\leq \lambda
\cos(\theta)
$. Therefore the hypothesis $\cos(\theta)\leq 2\Mod\sigma$ gives:
$$
\Vol(\NN_{R+1}(\Sigma)\cap S_{p})\leq \lambda\sinh(2R+2)\Mod\sigma .
$$
The same bound can be applied to $\Vol(\NN_{R+1}(\Sigma)\cap
S_{q})$. This proves Lemma 4.8. \qed
\enddemo

To achieve the proof of Lemma 4.5 we need the following:

\proclaim{Lemma 4.9} There is a universal constant $\mu>0$ such that if
$\max\{\cos\theta,\cos\phi\}> 2\Mod\sigma$, then
one can find a minimizing path $\sigma'$ between two singular points,
that satisfies the following:
 \roster
\item"(i)" $\Int(\sigma')\cap\Sigma=\emptyset$,
\item"(ii)" $\Mod{\sigma'}\leq\Mod\sigma$,
\item"(iii)" $\sigma\subset\NN_{\mu}(\sigma')$, 
\item"(iv)" if $\theta'$ and $\phi'$ are the angles between $\sigma'$
and $\Sigma$ at the end-points of $\sigma'$ then
$\max\{\cos\theta',\cos\phi'\}\leq2\Mod{\sigma}
$ (note that we are using $\Mod{\sigma}$ instead of $\Mod{\sigma'}$).
\endroster \endproclaim

Assuming Lemma 4.9, if the minimizing arc $\sigma$ does not fulfils
the hypothesis of Lemmas 4.7 and 4.8, then we apply the bounds
obtained in these lemmas to the region  $
\NN_{R+1+\mu}(\sigma')$ that contains
$\NN_{R+1}(\sigma)$. Thus we
obtain the bound
$$
\Vol(\NN_{R+1}(\sigma))\leq \Vol(\NN_{R+1+\mu}(\sigma'))\leq c_3(R+\mu)
\Mod{\sigma}
$$
where $c_3(R+\mu)>0$ depends only on $R$ because $\mu$ is universal. This
inequality completes the proof of Lemma 4.5. We next prove Lemma 4.9.

\demo{Proof of Lemma 4.9}
By Lemmas 4.7 and 4.8,  we can assume  
$\max\{\cos\theta,\cos\phi\}>2\Mod\sigma$.
 For $\varepsilon>0$
sufficiently small there is a homotopy
$\{\sigma_t\}_{t\in [0,\varepsilon)}$ of $\sigma=\sigma_0$ such that, for
any
$t\in [0,\varepsilon)$, $\sigma_t$ is a geodesic arc between two points of
$\Sigma$ satisfying the following:
\roster
\item"(1)" $\Int(\sigma_t)\cap\Sigma=\emptyset$, for all
$t\in [0,\varepsilon)$;
\item"(2)" the length $\Mod{\sigma_t}$ is decreasing with $t$;
\item"(3)" the angles $\theta_t,\phi_t\in [0,\pi/2]$
between $\sigma_t$ and $\Sigma$ are increasing with $t$.
 \endroster
When we increase the parameter $t$, we end up with one of the following
possibilities.
\roster
\item"(a)" either for some parameter $t_0$ we reach a path
$\sigma_{t_0}$ that satisfies (i) and (ii), and moreover
$\max\{\cos\theta_{t_0},\cos\phi_{t_0}\}
\leq2\Mod\sigma$;
\item"(b)" or before reaching such a $t_0$ the homotopy crosses
$\Sigma$: there is $t_1>0$ such that
$\Int(\sigma_{t_1})\cap\Sigma\neq\emptyset$ and, for any
$t\in[0,t_1)$, 
$\max\{\cos\theta_t,\cos\phi_t\}>2\Mod{\sigma}$.
\endroster
Both possibilities happen at bounded distance, because of
the following claim:

\proclaim{Claim 4.10}  In both cases, $d(\sigma_t,\sigma)\leq 1$,
where $t\leq t_0$ in case (a) and $t\leq t_1$ in case (b). 
\endproclaim

\demo{Proof} By using developping maps, we embed the homotopy
$\{\sigma_t\}_{t\in[0,\varepsilon)}$ locally isometrically in $\Bbb
H^3_K$. In particular, the pieces of $\Sigma$ and the arcs $\sigma_t$
are embedded as  geodesic arcs. 

Up to permutation, we can assume that $\cos\theta\geq\cos\phi$, where
$\theta$ and
$\phi$ are the angles at $p$ and $q$ respectively.
 Let $p_t\in\Sigma$
be the end-point of $\sigma_t$ obtained by moving $p$ along $\Sigma$
and let $p_t'$ its orthogonal projection to the geodesic of $\Bbb
H^3_K$
containing $\sigma$. This projection $p_t'$ lies between $p$ and $q$
(hence it is contained in $\sigma$) by construction of the
homotopy (see Fig. III.6). In particular
$d(p,p_{t}')\leq\Mod\sigma$. 

\midinsert
 \centerline{\BoxedEPSF{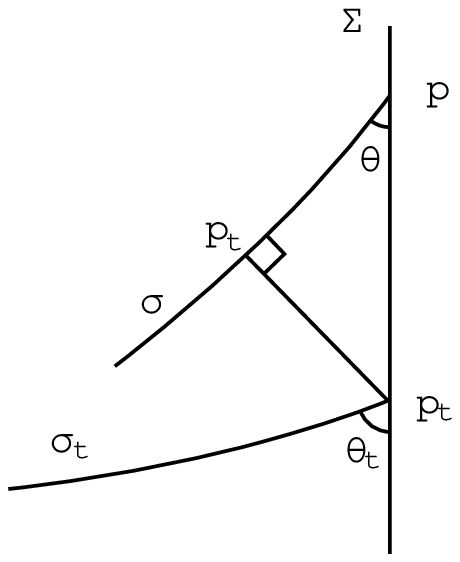 scaled 750}}
   \botcaption{Figure III.6}
    \endcaption
     \endinsert

Next we consider the right-angle
triangle with vertices $p$, $p_t$ and $p_t'$ (cf. Fig. III.6) and we
apply the trigonometric formula for $\cos\theta$:
 $$
\cos\theta=\tanh_K(d(p,p_t'))/\tanh_K(d(p,p_t))\leq
  \tanh_K\Mod\sigma/\tanh_K(d(p,p_t)),
$$
where $\tanh_K(r)=\tanh(\sqrt{-K} r)/\sqrt{-K}$ for $K>0$ and
$\tanh_0(r)=r$. We recall that $\tanh_K(\Mod x)\leq\Mod x$.

Since $\cos\theta\geq 2\Mod\sigma\geq 2\tanh_K\Mod\sigma$,
 we obtain
that $\tanh_K(d(p,p_t))\leq 1/2$. Moreover $1/2\leq \tanh_K1$ because
$K\in [-1,0]$. Thus the motonicity of $\tanh_K$ implies  that $d(p,p_t)\leq
1$, and  in particular $d(\sigma,\sigma_t)\leq 1$. This proves the claim.
\qed
\enddemo

From this claim, we deduce that $\sigma_t\subset\NN_2(\sigma)$,
because $\Mod{\sigma_t}\leq\Mod\sigma\leq 1$. In particular,
if case (a) above occurs, the path $\sigma_{t_0}$ satisfies the
conclusion of Lemma 4.9 and we are done. Hence we assume that case
(b) happens. Let $\sigma_{t_1}$ be the path coming from the homotopy
that intersects $\Sigma$ in its interior and consider $p_1$ and
$q_1$ the nearests two distint points of $\Sigma\cap\sigma_{t_1}$.
We obtain in this way two points $p_1$ and $q_1$ on $\Sigma$ joined
by a minimizing arc $\sigma_1$ such that
$\sigma_1\cap\Sigma=\{p_1,q_1\}$. Moreover, by Claim 4.10,
$d(\sigma_1,\sigma)\leq 2$ and, by the choice of
$p_1$ and
$q_1$,
$$
\Mod{\sigma_1}\leq \frac12\Mod{\sigma_{t_1}}\leq \frac12\Mod{\sigma}
\leq\frac14.
$$

\midinsert
 \centerline{\BoxedEPSF{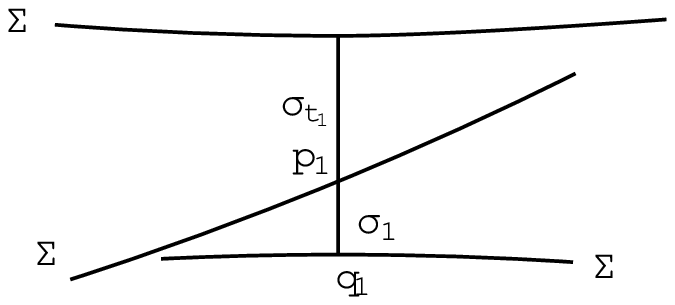 scaled 750}}
  \botcaption{Figure III.7}
   \endcaption
    \endinsert

Let $\theta_1,\phi_1\in [0,\pi/2]$ be the angles between
$\sigma_1$ and $\Sigma$ at $p_1$ and $q_1$ respectively. We assume
again that
$
\max\{\cos\theta_1,\cos\phi_1\}>2\Mod\sigma,
$
otherwise the minimizing path $\sigma_1$ would satisfy the conclusion
of Lemma 4.9 and we would be done. 

By iterating this process, we construct two sequences of points $p_n$
and $q_n$ on $\Sigma$ such that $p_0=p$, $q_0=q$, $p_n\neq q_n$ and
there is a minimizing path $\sigma_n$ between $p_n$ and $q_n$ such
that $\sigma_n\cap\Sigma=\{ p_n,q_n\}$ and
$\Mod{\sigma_n}\leq\frac12\Mod{\sigma_{n-1}}$. Moreover, if
$\theta_n,\phi_n\in[0,\pi/2]$ are the angles between $\sigma_n$ and
$\Sigma$ at $p_n$ and $q_n$, then we make the choice
$\cos\theta_n\geq \cos\phi_n$. There are two possibilities:
\roster
\item"--" either  $\cos\theta_n\leq 2\Mod\sigma$ and the
sequences stop at $n$.
\item"--" or 
$\cos\theta_n>2\Mod\sigma$ and the sequences
go on.
\endroster
The following claim shows that the sequences stop at uniformly bounded
distance.

\proclaim{Claim 4.11} There is a universal constant $\eta>0$ such that
$d(p,p_n)<\eta$ and $d(p,q_n)<\eta$, whenever $p_n$ and $q_n$ are
defined.
\endproclaim

The claim implies that the sequences stop, otherwise
$(p_n)_{n\in\N}$ would have a convergent subsequence in the compact
ball $B(p,\eta)$, contradicting the fact that the cone-injectivity
radius of the limit point is positive. Hence,  for the value of $n$
where the sequences stop, the path
$\sigma_n$ satisfies the conclusions of Lemma 4.9, because it is at
a uniformly bounded distance of $\sigma$. 

To achieve the proof of Lemma 4.9 it remains to prove Claim 4.11.

\demo{Proof} Let $p_n$ and $q_n$ be two points of the sequences on
$\Sigma$, and $\sigma_n$ the minimizing arc between them such that
$\sigma_n\cap\Sigma=\{p_n,q_n\}$ and
$\Mod{\sigma_n}\leq 2^{-n}\Mod\sigma$. We also assume
that $\cos\theta_n>2\Mod\sigma$, so that
$p_{n+1}$ and $q_{n+1}$ are defined. These points 
are constructed  by considering a homotopy of $\sigma_n$ as in the
beginning of the proof of Lemma 4.10. This homotopy gives:
\roster
\item"--" either a path
$\sigma'_n$ that crosses
$\Sigma$, and the points $p_{n+1}$ and $q_{n+1}$ are the  two
nearests different points in $\sigma'_n\cap\Sigma$;
\item"--" or a path $\sigma'_n$ such that the angles $\theta_{n+1}$
and $\phi_{n+1}$ between $\sigma'_n$ and $\Sigma$ satisfy
$
\max\{\cos\theta_n',\cos\phi_n\}\leq
2\Mod\sigma$. In this case $p_{n+1}$ and
$q_{n+1}$ are the end-points of $\sigma_{n+1}=\sigma_n'$ and the
sequences stop at $n+1$.
\endroster
In both cases we have:
$$
\max\{d(p_n,p_{n+1}),d(p_n,q_{n+1})\}\leq
\Mod{\sigma_n'}+d(p_n,\sigma_n'),
$$
and $\Mod{\sigma_n'}\leq \Mod{\sigma_n}\leq
\Mod\sigma/2^n\leq1/2^{n+1}$. The trigonometric argument of Claim 4.10
applies here to give the following inequality,
$$
\cos\theta_n\leq\tanh_K\Mod{\sigma_n}/\tanh_K(d(p_n,\sigma_n')).
$$
Combining this with the hypothesis $\cos\theta_n>
2\Mod{\sigma}$ we get:
$$
\tanh_K(d(p_n,\sigma_n'))
\leq\frac{\tanh_K\Mod{\sigma_n}}{2\Mod\sigma}
\leq\frac{\Mod{\sigma_n}}{2\Mod\sigma}
\leq\frac1{2^{n+1}}
$$
Since $\tanh_K(x)=x+O(\Mod x^3)$, it follows from 
this inequality that there is a universal constant $\eta_0>0$ such that
$
d(p_n,\sigma_n')\leq \eta_0/2^{n+1}
$. Summarizing these inequalities we obtain:
$$
\max\{d(p_n,p_{n+1}),d(q_n,q_{n+1})\}\leq (\eta_0+1)/2^{n+1}
$$
and
$$
\max\{d(p,p_{n+1}),d(q,q_{n+1})\}\leq
\sum_{i=0}^{n}(\eta_0+1)/2^{i+1}< \eta_0+1.
$$
It suffices to take $\eta=\eta_0+1$ to achieve the proof of Claim
4.11.\qed
\enddemo
\enddemo
\enddemo
\enddemo
\enddemo

\head 5. Some properties of geometric convergence
\endhead

In this section we study properties of  sequences
of pointed cone 3-manifolds  in $\CC_{[\omega,\pi],a}$ that converge
geometrically. During all the section we will assume
$\omega\in (0,\pi]$ and $a>0$.

\proclaim{Proposition 5.1} Let $(C_n,x_n)$  be a sequence in
$\CC_{[\omega,\pi],a}$ that converges geometrically to a pointed cone
3-manifold $(C_{\infty},x_{\infty})$. Then:
\roster
\item"(i)" the curvature of $C_n$ converges to the 
curvature of $C_{\infty}$;
\item"(ii)"
$\Inj(x_{\infty})\leq\lim\limits_{n\to\infty}\!\inf\, \Inj(x_n)$.
\endroster
\endproclaim

\demo{Proof} Property (i) has been proved in Proposition 2.3. It
follows also from the fact that the sectional curvature may be computed
from small geodesic triangles.

To prove Property (ii) we distinguish two cases, according to whether
the cone- injectivity radius at $x_{\infty}$ is estimated using
singular or non-singular balls. Let $r_{\infty}=\Inj(x_{\infty})$. We
first assume that for any $0<\varepsilon<r_{\infty}$, the ball
$B(x_{\infty},r_{\infty}-\varepsilon)$ is standard and non-singular.
Geometric convergence implies that for $n$ sufficienly large,
$B(x_n,r_{\infty}-2\varepsilon)$ is standard in $C_n$, hence
$\Inj(x_n)\geq\Inj(x_{\infty})-2\varepsilon$. A similar argument applies
in the case of singular standard balls. 
\qed
\enddemo

\proclaim{Proposition 5.2} Let $(C_n,x_n)$ be a sequence in
$\CC_{[\omega,\pi],a}$ that converges geometrically to a pointed cone
3-manifold $(C_{\infty},x_{\infty})$. For any compact subset $A\subset
C_{\infty}$ there exists  $n_0>0$ such that  for $n\geq n_0$ 
there is an embedding $f_n:A\to C_n$ with the following properties:
\roster
\item"(i)" $f_n(A)\cap\Sigma_n=f_n(A\cap\Sigma_{\infty})$;
\item"(ii)" the cone angles at $f_n(A)\cap\Sigma_n$ approach the cone
angles at $A\cap\Sigma_{\infty}$ as $n$ goes to infinity.
\endroster
\endproclaim

\proclaim{Corollary 5.3} If the limit $C_{\infty}$ of a geometrically
convergent sequence $(C_n,x_n)$  in
$\CC_{[\omega,\pi],a}$ is compact, then for $n$ sufficiently large
$C_n$ has the same topological type as $C_{\infty}$ (i.e. the pairs 
$(C_n,\Sigma_n)$ and $(C_{\infty},\Sigma_{\infty})$ are homeomorphic)
and the cone angles of $C_n$ converge to those of $C_{\infty}$. \qed
\endproclaim

\demo{Proof of Proposition 5.2} Since $A$ is compact, there exists $R>0$
 such that $A\subset B(x_{\infty},R)$. The definition of geometric
convergence implies
 that for any $\varepsilon>0$ and for $n\geq n_0$  ($n_0$
depending on $R$, $\varepsilon$ and the sequence) there exists a
$(1+\varepsilon)$-bilipschitz map $f_n:B(x_{\infty},R)\to C_n$ such
that 
$$
f_n(B(x_{\infty},R)\cap\Sigma_{\infty})=f_n(B(x_{\infty},R))\cap\Sigma_n.
$$
This proves Property (i) of the proposition. Moreover, by taking
$\varepsilon\to 0$ we get Property (ii). \qed
\enddemo

\proclaim{Proposition 5.4} Let $(C_n,x_n)$ be a sequence in
$\CC_{[\omega,\pi],a}$ that converges geometrically to a pointed cone
3-manifold $(C_{\infty},x_{\infty})$. Let $A\subset C_{\infty}$ be
a compact subset. Let $\rho_n$ (resp. $\rho_{\infty}$) be the holonomy
of the cone 3-manifold $C_n$ (resp. $C_{\infty}$). Then, we can choose
the holonomy representations  and the embeddings $f_n$ of Proposition
5.2 so that:
\roster
\item"(i)" If the curvature of $C_n$ does not depend on $n$, then
$\forall\gamma\in\pi_1(A-\Sigma_{\infty},x_{\infty})$,
$$
\lim_{n\to\infty}\rho_n(f_{n*}(\gamma))=\rho_{\infty}(\gamma).
$$
\item"(ii)" If the curvature $K_n\in [-1,0)$ of $C_n$ converges to $0$
as $n\to\infty$, then
$\forall\gamma\in\pi_1(A-\Sigma_{\infty},x_{\infty})$,
$$
\lim_{n\to\infty}\rho_n(f_{n*}(\gamma))=
\operatorname{ROT}(\rho_{\infty}(\gamma)),
 $$
where $\operatorname{ROT}:\operatorname{Isom}(\Bbb E^3)\to O(3)$ is the
surjective morphism whose kernel is the subgroup of translations.
\endroster
\endproclaim

\demo{Proof}
We first prove assertion (i). Given 
$\gamma\in\pi_1(A-\Sigma_{\infty},x_{\infty})$, we realize 
$\gamma$ by a piecewise geodesic loop $\sigma:[0,l]\to A-\Sigma$.
Let $0=t_0<t_1<\cdots < t_r=l$ be the associated partition. That is, each
piece $\sigma([t_i,t_{i+1}])$ is a minimizing geodesic. 
We take
each geodesic piece $\sigma([t_i,t_{i+1}])$ small enough to be contained
(as well as its image by $f_n$) in non-singular standard balls in
$C_{\infty}$ (or in $C_n$). Let $\sigma_n:[0,l]\to C_n-\Sigma_n$ be the
piecewise geodesic path such that $\sigma_n(t_i)=f_n(\sigma(t_i))$ and
each
piece $\sigma_n([t_i,t_{i+1}])$ is a minimizing geodesic. By
construction the loop
$\sigma_n$ is homotopic to $f\circ \sigma$ relatively to the  base point.
Moreover, by  using $(1+\varepsilon)$-bilipschitz maps with $\varepsilon\to
0$, the length $\Mod{\sigma_n([t_i,t_{i+1}])}$ converges to 
$\Mod{\sigma([t_i,t_{i+1}])}$ and the angle at $\sigma_n(t_i)$ 
$\,\angle (\sigma_n([t_i,t_{i+1}]),\sigma_n([t_{i+1},t_{i+2}]))$ converges
to 
$\angle (\sigma([t_i,t_{i+1}]),\sigma([t_{i+1},t_{i+2}]))$.
 It follows that $\rho_{\infty}(\gamma)$ is the limit of
$\rho_n(f_{n*}(\gamma))$.

Assertion (ii) is proved in \cite{\PoOne, Prop. 5.14 (i)}. 
\qed
\enddemo

\newpage
\rightheadtext{ }
\leftheadtext{IV \qquad    Local Soul}

\

\centerline{\smc chapter \  iv} 

\

\centerline{\chapt LOCAL  \, SOUL  \, THEOREM \,  FOR \,  CONE \, 3-MANIFOLDS}

\vglue.2cm

\centerline{\chapt WITH \,  CONE  \, ANGLES  \,  LESS \,  THAN \,  OR  \, EQUAL \,  TO
\, 
$\pi$}

\

\

The goal of this chapter is to describe the metric
structure of a neighborhood of a point with sufficiently
small injectivity radius in a hyperbolic cone 3-manifold with
 cone angles bounded above
by $\pi$. This description is crucial to study collapsing
sequences of cone 3-manifolds in the proofs of Theorems A
and B.

We need first the following definition.

\definition{Definition} 
Let $C^3$ be a  cone 3-manifold and $D$ a
cone manifold of dimension less than $3$, possibly with
silvered boundary $\partial D$. A surjective map
$p\!:\! C\to D$ is said to be a \it cone fiber bundle \rm if 
\roster
\item"-" on $D-\partial D$, the restriction of $p$ is a
locally trivial fiber bundle with fiber a cone manifold.
Moreover, if $\dim(D)=2$ then
$p(\Sigma_C)=\Sigma_D$
\item"-" on $\partial D$, the restriction of $p$ is an orbifold fibration. In
particular, the fiber over a point of $\partial D$  is an orbifold with
cone angles equal to $\pi$.
\endroster
\enddefinition

\proclaim{Theorem \rm (Local Soul Theorem)}
Given $\omega\in (0,\pi)$, $\varepsilon>0$ and $D>1$
there exist $\delta=\delta(\omega,\varepsilon,D)>0$ and
$R=R(\omega,\varepsilon,D)>D>1$ such that, if $C$  is an
 oriented hyperbolic cone 3-manifold with cone
angles in $[\omega,\pi]$ and if $x\in C$ satisfies
$\Inj(x)<\delta$, then:

- either $C$ is $(1+\varepsilon)$-bilipchitz homeomorphic
to a compact Euclidean cone 3-man\-ifold $E$ of diameter 
$\Diam(E)\leq R \Inj(x)$;

- or $x$ has an open neighborhood $U_x\subset C$
which is $(1+\varepsilon)$-bilipschitz homeomorphic to the
normal cone fiber bundle $\NN_{\nu}(S)$, of radius $0<\nu<1$, of
the soul $S$ of a non-compact orientable Euclidean cone
3-manifold  with cone angles in $[\omega,\pi]$.
 Moreover, according to
$\dim(S)$, the Euclidean non-compact
cone 3-manifold  belongs to the following list:
\roster
\item"(I) " (when $\dim(S)=1$), $S^1\ltimes\Bbb R^2$,
$S^1\ltimes(\text{open cone disk})$ and the pillow (see
Figure IV.1), where $\ltimes$ denotes the metrically twisted product;
\item"(II) "(when $\dim(S)=2$)
\item"(i)" a product $T^2\times \Bbb R$;
$S^2(\alpha,\beta,\gamma)
\times\Bbb R$, with
$\alpha+\beta+\gamma=2\pi$;
$S^2(\pi,\pi,\pi,\pi)\times\Bbb R$;
\item"(ii)" the orientable twisted line bundle over the Klein bottle
$K^2\tilde\times \Bbb R$ or
over the projective plane with two silvered points
$\Bbb P^2(\pi,\pi)\tilde\times \Bbb R$;
\item"(iii)" a quotient by an involution
 of either  $S^2(\pi,\pi,\pi,\pi)\times\Bbb
R$, $T^2\times\Bbb R$ or $K^2\tilde\times \Bbb R$,
 that gives an orientable  bundle respectively over
either $D^2(\pi,\pi)$, an annulus, or a M\"obius strip,
with silvered boundary in the three cases (see Figure IV.2).
\endroster
Moreover the $(1+\varepsilon)$-bilipschitz
homeomorphism 
$f\!:\! U_x\to\NN_{\nu}(S)$ satisfies the inequality
$$
\max\big(\Inj(x),d(f(x),S),\Diam(S)\big)\leq \nu/D.
$$
\endproclaim

\midinsert
 \centerline{\BoxedEPSF{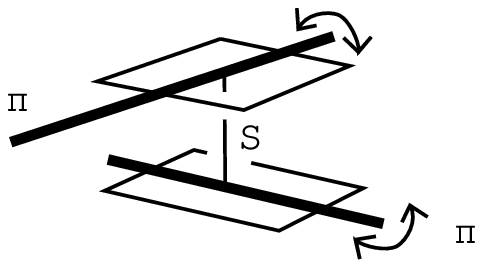 scaled 750}}
    \botcaption{Figure IV.1} The pillow. Its soul is the
interval $[0,1]$ with silvered boundary.
    \endcaption
     \endinsert

\midinsert
 \centerline{\BoxedEPSF{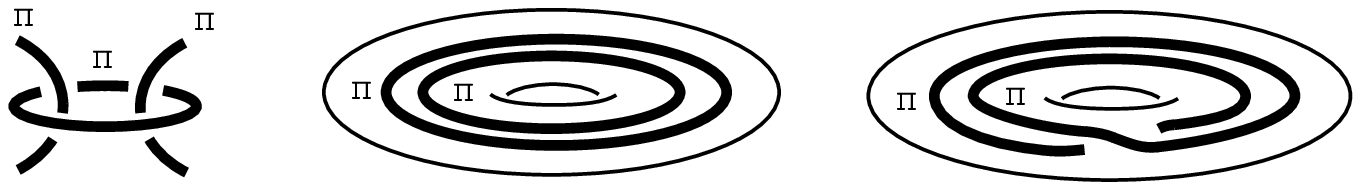 scaled 750}}
    \botcaption{Figure IV.2} From left to right, the
non-compact
Euclidean cone 3-manifolds with soul $D^2(\pi,\pi)$, an
annulus and a M\"obius strip, with silvered boundary in
every case. They are the respective quotients of 
$S^2(\pi,\pi,\pi,\pi)\times\Bbb
R$, $T^2\times\Bbb R$ and $K^2\tilde\times \Bbb R$ by an involution.
    \endcaption
     \endinsert

The Euclidean cone 3-manifolds $E$ in the Local Soul  Theorem
are called the
\it local models.
\rm
We call $S$ the \it soul, \rm because, in each case,
 $S$ is a
totally convex cone submanifold of the
local model $E$, and $E$ is isometric to the normal cone
fiber bundle of $S$.

The first step of the proof is Thurston's Classification Theorem
of non-compact orientable Euclidean cone 3-manifolds. We
need in fact a simpler classification, because Thurston's
classification includes general singular locus, and we
consider here only  the case where the singular locus 
is a 1-dimensional submanifold (cf.
\cite{\SOK},
\cite{\Hod} and \cite{\Zhou}).

\rightheadtext{IV \qquad    Local Soul}

\head 1. Thurston's classification Theorem of non-compact Euclidean
cone 3-manifolds
\endhead

In this section we prove the following result:

\proclaim{Theorem 1.1  \rm (Thurston)} Let $E$ be a
non-compact orientable Euclidean cone $3$-manifold with
cone angles less than or equal to $\pi$ and  a 
1-dimensional submanifold as singular locus. 
Then either $E=\Bbb R^3$,   $E=\Bbb R^3(\alpha)=
\Bbb R\times(\text{open cone disk})$, or $E$ is one  of
the local models given in the
Local Soul Theorem.
\endproclaim

We prove this theorem by using the Soul Theorem for
Euclidean  cone 3-man\-ifolds, which is just a version of Cheeger-Gromoll
Soul Theorem for Riemannian manifolds with non-negative
curvature  \cite{CGl} (see also \cite{\Sak}).
We need first some definitions.

\definition{Definition}
Let $C$ be a cone 3-manifold with singular locus $\Sigma\subset
C$ a 1-dimensional submanifold and curvature $K\in[-1,0]$.
\roster
\item"--" The \it silvered points \rm of $C$ are the points
of
$\Sigma$ having cone angle $\pi$.
\item"--" A path $\gamma\!:\! [0,l]\to C$ is \it geodesic \rm if
it is locally minimizing.
\item"--" A path $\gamma\!:\! [0,l]\to C$ is \it s-geodesic \rm
if it is locally minimizing except for some $t\in (0,l)$ where
$\gamma(t)$ is silvered and the following happens: there
exist $\varepsilon>0$ and a neighborhood $U$ of $\gamma(t)$
such that $\gamma\!:\! (t-\varepsilon,t+\varepsilon)\to U$ lifts
to a minimizing path in the double  cover $\tilde
U\to U$ branched along $\Sigma\cap U$. See Figure IV.3.
\item"--" A subset $S\subset C$ is \it totally s-convex \rm
if every s-geodesic path with end-points in $S$ is
contained in $S$.
\endroster
\enddefinition

\midinsert
 \centerline{\BoxedEPSF{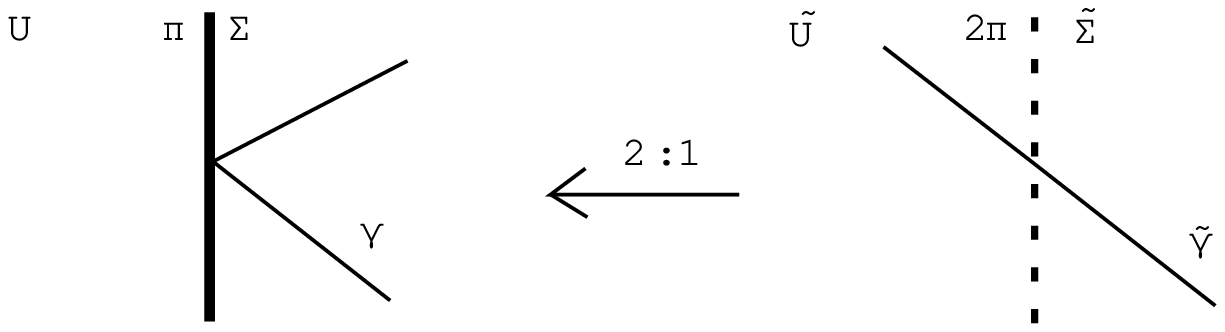 scaled 750}}
    \botcaption{Figure IV.3} Example of s-geodesic.
    \endcaption
     \endinsert

The notion of s-geodesic generalizes the notion of
geodesic, thus total s-convexity is stronger than  usual total
 convexity, as shown in the following example.

\remark{Example} Let $A$ be a totally s-convex subset of a
cone 3-manifold $C$ and $\gamma$ a geodesic path from 
$p\in A$ to $q\in\Sigma$. If $q$ is a silvered point and
$\gamma$ is orthogonal to $\Sigma$, then $\gamma$ is
contained in $A$, because the path $\gamma
* \gamma^{-1}$ is a s-geodesic with end points in $A$.  In
particular, a totally s-convex set intersects all singular
components having cone angle $\pi$.
\endremark

\proclaim{Theorem 1.2 \rm (Soul Theorem)} Let $E$ be a
non-compact Euclidean cone 3-man\-ifold with cone
angles
$\leq\pi$. Then
$E$ contains a compact totally s-convex cone submanifold
$S\subset E$ of dimension $0$, $1$ or $2$, with silvered or
empty boundary.
Moreover
$E$ is isometric to the normal cone fiber bundle of
$S$. 
\endproclaim

Given a cone submanifold $S\subset E$, there is an
$\varepsilon>0$ such that the tubular neighborhood of radius
$\varepsilon>0$,
$\NN_{\varepsilon}(S)$, is a cone bundle over $S$. When we
say that $E$ is isometric to the normal
cone fiber bundle of
$S$, we mean that we can choose  the radius
$\varepsilon=\infty$ and that the metric has the local product
structure of the bundle.

The cone submanifold $S$ is called the soul of $E$.
 As we said above, this theorem is a cone manifold version
of Cheeger-Gromoll's theorem \cite{CGl}.
We postpone the proof of the Soul Theorem to Sections 2 and 3, and
now we use it to prove Thurston's Classification Theorem
1.1.

\demo{Proof of Theorem 1.1}
Let $E$ be an orientable non-compact Euclidean cone 3-man\-ifold and
$\Sigma$ its singular locus. Note that every finite 
covering
$\tilde S\to S$ (possibly branched at silvered
points) induces a covering 
$\tilde E\to E$. Moreover, $\tilde S$ is the soul of 
$\tilde E$, because  $\tilde S$
is totally s-convex and  $\tilde E$ is the normal cone fiber bundle of 
$\tilde S$. Passing to finite  coverings will help us to simplify
the proof.

We distinguish three cases,
according to the dimension of the soul $S\subset E$.

When $\dim(S)=0$, then $S$ is a point $p$. For
$\varepsilon>0$ sufficiently small, the ball
$B(p,\varepsilon)$ of radius $\varepsilon$ is either a non-singular Euclidean ball
or a ball with a singular axis. Hence, since $E$ is
isometric to the normal cone fiber bundle of the point $p$, either $E$ is the Euclidean
space $\Bbb R^3$, or $E=\Bbb R^3(\alpha)=
\Bbb R\times(\text{open cone disk})$.    

When $\dim(S)=1$, then either $S$ is $S^1$ or an interval $[0,1]$
with silvered boundary. 

If $S=S^1$, then by convexity  either
 $S\subset \Sigma$ or $S\cap\Sigma=\emptyset$. Since $E$ is
orientable,
$\NN_{\varepsilon}(S)$ is a solid torus, possibly with a
singular core. Therefore,  $E$ is the (metrically twisted)
product of
$S^1$ with an infinite disk, possibly with a singular cone point
in the center.

If $S=[0,1]$, then  we consider the double covering $S^1\to[0,1]$
branched along the silvered boundary. The induced double branched
covering
$\tilde E\to E$ is $S^1\ltimes \Bbb R^2$, thus $E$ is the
pillow ($\Bbb R^3$ with two  silvered lines; see Figure IV.1).

When $\dim(S)=2$, we use the clasification of compact
Euclidean cone 2-man\-ifolds  with geodesic boundary and
cone angles $\leq\pi$. This clasification is easily deduced
from Gauss-Bonnet formula. In particular, either
$S=S^2(\alpha,\beta,\gamma)$ with
$\alpha+\beta+\gamma=2\pi$ or $S$ is a Euclidean orbifold
having as finite covering $\tilde S=T^2$. If
$S=S^2(\alpha,\beta,\gamma)$ or $ S=T^2$, then $S$ is two
sided and
$\NN_{\varepsilon}(S)=S\times (-\varepsilon,\varepsilon)$.
Therefore $E=S\times \Bbb R$. The remaining cases reduce
to study finite groups of isometries of $T^2$ and their
orientable isometric extension to $T^2\times \Bbb R$.
\qed \enddemo

\head 2. Totally s-convex subsets in Euclidean cone 3-manifolds.
\endhead

In this section we give some basic facts about totally s-convex subsets
in Euclidean  cone 3-manifolds,
which are used in the proof of the Soul Theorem
(in Section 3). 
Lemma 2.1 shows that totally s-convex subsets appear naturally
as level sets of continuous convex functions.

\definition{Definition} A continuous function $f\!:\! E\to\Bbb R$
on a cone 3-manifold $E$
is \it convex \rm if  $f\circ\gamma$ is convex for every
geodesic path $\gamma\!:\! [0,l]\to E$.
\enddefinition

\proclaim{Lemma 2.1} If $f\!:\! E\to\Bbb R$ is a continuous convex
function then, for every {\rm s-geodesic} path
$\gamma\!:\! [0,l]\to E$, $f\circ\gamma$ is convex. In
particular, the subset $\{x\in E\mid f(x)\leq 0\}$ is
totally s-convex.
\endproclaim

\demo{Proof}
Every s-geodesic path is locally the limit of geodesic
paths, arbitrarily close to the singular set but disjoint
from it. It follows by continuity that the inequalities
defining  convexity are satisfied locally for every
s-geodesic path. \qed 
\enddemo

\definition{Definition} Let $A\subset E$ be  a smooth submanifold
without boundary. We say that $A$ is \it totally geodesic \rm if
either $\dim A=3$ or for every  $x\in A$ the following hold:
\roster
\item"-" if $x\not\in \Sigma$, then the second fundamental form of $A\subset E$ at $x$ is trivial;
\item"-" if $x\in\Sigma$ and $\dim A=2$, then $A$ and $\Sigma$ are orthogonal at $x$;
\item"-" if $x\in\Sigma$ and $\dim A=1$, 
then
there is a neighborhood $U\subset E$ of $x$ such that $\Sigma\cap U= A\cap U$.
\endroster
For non-singular points, this definition  coincides with the usual definition 
in Riemannian geometry. We also remark that this is a local notion that does not require $A$ to be
complete.
\enddefinition

\proclaim{Proposition 2.2} Let   $E$ be a Euclidean cone 3-manifold and
$A\subset E$  a non-empty 
closed totally
s-convex subset.
 Then  $A$ is an embedded manifold, possibly with boundary, whose interior 
$A-\partial A$ is totally geodesic. 
\endproclaim

Before proving this proposition we need the following lemma, which
describes the local structure of $A$ at  the singular points.

\proclaim{Lemma 2.3} Let $x\in\Sigma$  and let
$D^2(x,\varepsilon)$ be the geodesic singular disk
transverse to $\Sigma$ with center $x$ and radius
$\varepsilon>0$. If $A$ is a totally s-convex
subset such that $A\cap D^2(x,\varepsilon)\neq\emptyset$,
 then one of
the following  possibilities happens:
\roster
\item"(i)"  $A\cap D^2(x,\varepsilon)=\{x\}$,
\item"(ii)" $A\cap D^2(x,\varepsilon)$ contains a smaller
disk $D^2(x,\delta)$, with $0<\delta<\varepsilon$,
\item"(iii)"  $x$ is silvered and $A\cap D^2(x,\varepsilon)$
 is a segment orthogonal to $\Sigma$ at
$x$.
\endroster
\endproclaim

\demo{Proof} Choose $\varepsilon>0$ so that
$D^2(x,\varepsilon)$ is a disk contained in a standard
ball. We prove first that if the cone angle $\alpha$ at $x$ is less
than $\pi$ then the convex hull of a point $y\in A\cap(
D^2(x,\varepsilon)-\{x\})$ contains a disk
$D^2(x,\delta)$, with
$0<\delta<\varepsilon$. We view $D^2(x,\varepsilon)$ as
the quotient of an angular sector $S_{\alpha}$ whose faces
are identified by an isometric rotation, and such that $y$ is obtained
by identifying two points $\tilde y_1,\tilde y_2\in
S_{\alpha}$, one in each face of $S_{\alpha}$.  Consider
the geodesic path $\tilde \sigma\!:\! [0,l]\to S_{\alpha}$
minimizing the distance between $\tilde y_1$ and $\tilde
y_2$ (see Figure IV.4), it projects to a geodesic loop
$\sigma\!:\! [0,l]\to D^2(x,\varepsilon)$  based at $y$. The convex
hull of $y$ contains $\sigma$ and we have
$$
d(\sigma,x)= d(y,x)  \cos (\alpha/2).
$$
 By using this formula and the fact that
$0<\cos (\alpha/2)<1$, we can construct a sequence of
concentric geodesic loops converging to $x$, hence the
convex hull of
$y$ contains a disk
$D^2(x,\delta)$, with $\delta>0$. This proves the lemma
when $\alpha<\pi$.

Assume now that $x$ is a silvered point. Let  $y\in
A\cap(D^2(x,\varepsilon)-\{x\})$. The set $A$ contains the
minimizing path $\sigma$ from $y$ to $x$, because
$\sigma*\sigma^{-1}$ is a s-geodesic loop based at $y$.
Moreover, if $A\cap D^2(x,\varepsilon)$ contains two such
segments, then these segments divide $D^2(x,\varepsilon)$
into two sectors of angles less than $\pi$, therefore 
$D^2(x,\delta)\subset A$ with $\delta>0$. \qed
\enddemo

\midinsert
 \centerline{\BoxedEPSF{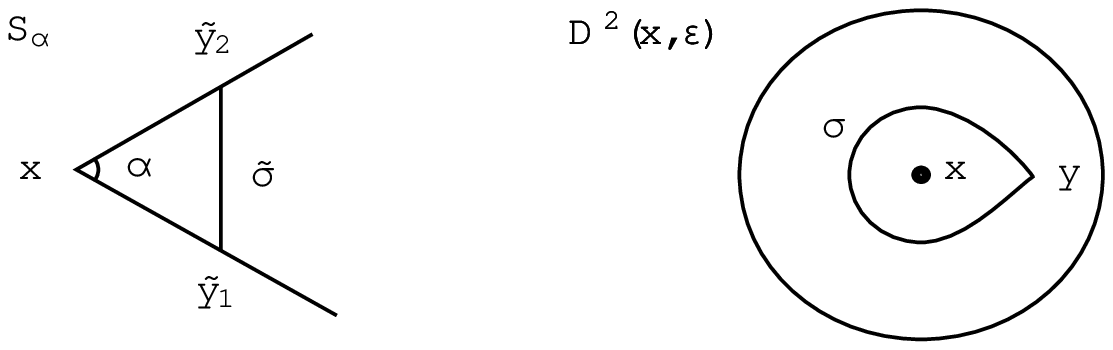 scaled 750}}
    \botcaption{Figure IV.4} The sector $S_{\alpha}$ and
the disk $D^2(x,\varepsilon)$.
    \endcaption
     \endinsert

\demo{Proof of Proposition 2.2} It suffices to prove the result locally:  
 every point  $x\in A$ has a neighborhood $U$ such that
 $\overline U\cap A$ is an embeded submanifold, possibly with boundary,
whose interior $\overline U\cap A-\partial (\overline U\cap A)$ is totally
geodesic.
If $x$ is non-singular then it is just a well known
result for locally convex subsets in $\Bbb R^3$. Hence we suppose that
$x\in\Sigma$.

We choose a neighborhood $U$ of $x$ that has a product structure. More
precisely $U$ is isometric to $D^2(0,\delta)\times(-\varepsilon,\varepsilon)$
for some $\varepsilon,\delta>0$, where $D^2(0,\delta)$ is a singular 2-disk
of radius $\delta$, with a singularity in its center $0$. For this
product structure we have that $U\cap\Sigma=\{0\}\times (-\varepsilon,\varepsilon)$
and $x= (0,0)$.

Since the intersection $A\cap U\cap \Sigma$ is a connected subset of 
$U\cap\Sigma$  containing $x$, there are the following three  possibilities:
\roster
\item"(a)"   $A\cap U\cap\Sigma=  \{ x\}$;
\item"(b)"   $A\cap U\cap\Sigma= \{0\}\times[0,\varepsilon)$, i.e. a subinterval of $\Sigma\cap U$
   having $x$ as end-point;
\item"(c)"   $A\cap U\cap\Sigma=\{0\}\times(-\varepsilon,\varepsilon) =\Sigma\cap U$.
\endroster

By using Lemma 2.3 we can describe explicitly  all the possibilities for
$A\cap \overline U$  in each case.
In case (a), when $A\cap U\cap \Sigma =\{x\}=\{(0,0)\}$, there are 3 subcases:
\roster
\item"(a1)" $A\cap U=\{x\}=\{(0,0)\}$.
\item"(a2)" $A\cap \overline U=\overline V\times\{0\}$, where $\overline V$ is a convex
neighborhood of $0$ in $D^2(0,\delta)$. In particular $\Int(A\cap U)$ is a totally geodesic
2-submanifold transverse to $\Sigma$
\item"(a3)"The point $x$ is silvered and $A\cap U$ is a segment orthogonal to $\Sigma$ at $x$.
\endroster
It folows from this explicit description that the proposition holds in the three subcases
(a1), (a2) and (a3).

In case (b), when
$A\cap U\cap\Sigma=\{0\}\times[0,\varepsilon) $, again there are  three subcases:
\roster
\item"(b1)" $A\cap U= A\cap U\cap\Sigma=\{0\}\times[0,\varepsilon) $, (i.e. $A\cap U$
is a subinterval of $\Sigma$).
\item"(b2)" For some $t_0\in [0,\varepsilon)$, $A\cap (D^2(0,\delta) \times \{t_0\})$ contains a
singular 2-disk with positive radius. 
\item"(b3)" $x$ is silvered and, for some $t_0\in[0,\varepsilon)$,
$A\cap (D^2(0,\delta) \times \{t_0\})$ is a segment perpendicular to $\Sigma$ at $x$.
\endroster
In subcase (b1) we have an explicit description of $A\cap U$ and we may conclude that 
the proposition
holds. To give an explicit description in the other cases we need further work.

In subcase (b2),
we claim that   $A\cap  U$ is a 3-manifold with boundary and that $x\in\partial(A\cap U)$. 
First we remark that, for every $t\in(0,\varepsilon)$, the intersection
$A\cap (D^2(0,\delta) \times \{t\})$ 
 contains a
singular 2-disk with positive radius, because $A\cap U$
contains the convex hull of the union of $A\cap (D^2(0,\delta) \times \{t_0\})$ and  
$A\cap U\cap\Sigma=\{0\}\times[0,\varepsilon) $.

We parametrize  $U$ by cylindrical coordinates $(r,\theta,t)\in[0,\delta)\times
[0,\alpha]\times (-\varepsilon,\varepsilon)$, where $r$ is the distance 
to $\Sigma$, $\theta\in [0,\alpha]$ is the angle parameter, $\alpha$ is the singular angle
and $t$ is the height parameter. Thus we identifie $(r,0,t)$  to $(r,\alpha,t)$, and
$(0,\theta,t)$ to $(0,\theta',t)$, for every $\theta,\theta'\in [0,\alpha]$.

By Lemma 2.3, if a point belongs to $A\cap U$, 
then so it does its projection to $\Sigma$. Therefore
there exists a function $f\!:\![0,\alpha]\times [0,\varepsilon]\to[0,\delta]$
 such that 
$$
A\cap\overline U=\{ (r,\theta,t)\in[0,\delta]\times
[0,\alpha]\times [0,\varepsilon] \mid r\leq f(\theta,t)  \}.
$$
We remark that  for every $t\in(0,\varepsilon)$ and every $\theta\in [0,\alpha]$,
$f(\theta,t)>0$
because $A\cap (D^2(0,\delta) \times \{t\})$  contains a
singular 2-disk with positive radius.

Next we show that   $f$ is continuous. The function $f$ is upper semi-continuous because
$A$ is closed. Moreover the lower semi-continuity of $f$ at a point $(\theta,t)$ can be proved by
considering the convex hull of the union of $A\cap (D^2(0,\delta) \times \{t\})$ and $A\cap\Sigma$,
because this convex hull has dimension 3. 

Since $f$ is continuous,  $A\cap\overline
U$ is a 3-manifold with boundary whose interior is 
$$\Int(A\cap\overline U)=\{ (r,\theta,t) 
\in [0,\delta]\times
[0,\alpha]\times (0,\varepsilon)
\mid r< f(\theta,t)  \}.$$
 Hence the proposition holds in   subcase (b2).

In subcase (b3) we claim that $A\cap \overline U$ is a 2-manifold that is 
 the union of a family of parallel segments perpendicular to $\Sigma$.
First we remark that, for every $t\in(0,\varepsilon)$, the intersection
$A\cap (D^2(0,\delta) \times \{t\})$ is a 
segment orthogonal to $\Sigma$, because $A\cap U$
contains the convex hull of the  union of the segment $A\cap (D^2(0,\delta) \times
\{t_0\})$ and
$A\cap U\cap\Sigma=\{0\}\times[0,\varepsilon) $. Moreover the segments 
$A\cap (D^2(0,\delta) \times \{t\})$
are parallel,
because if not the  convex hull of their union would have dimension 3 and we would be in subcase
(b2).
 
Again  we parametrize  $U$ by cylindrical coordinates $(r,\theta,t)$. By 
the same argument as in subcase (b2) we conclude that there exists a continuous function
$f\!:\!(-\varepsilon,\varepsilon)\to [0,\delta] $ such that
$$
A\cap\overline U=\{ (r,\theta,t)\in [0,\delta]\times
[0,\alpha]\times [0,\varepsilon] \mid \theta=0, \, r\leq f(\ t)  \}.
$$
Moreover, for $t>0$, $f(t)>0$. Hence $A\cap U$ is a 2-dimensional submanifold
with boundary  and with totally geodesic interior
$\{ (r,\theta,t) \mid \theta=0,\,  0<r< f(\ t)  \}$.
Thus the proposition follows from explicit description also in this case.

Finally, in case (c), when $A\cap U\cap\Sigma=\{0\}\times(-\varepsilon,\varepsilon)$,
there are again three subcases that can be treated with the same method as subcases (b1),
(b2) and (b3). These subcases are:
\roster
\item"(c1)" $A\cap U=A\cap U\cap\Sigma=\{0\}\times(-\varepsilon,\varepsilon)$.
\item"(c2)" $A\cap \overline U$ is a 3-submanifold with boundary that contains $U\cap\Sigma=
\{0\}\times(-\varepsilon,\varepsilon)$ in its interior.
\item"(c3)" $x$ is a silvered point and $A\cap\overline U$ is a 2-submanifold with boundary,
which is the
union of parallel segments orthogonal to $U\cap\Sigma$. In particular the interior of $A\cap\overline
U$ is totally geodesic and $U\cap\Sigma= 
\{0\}\times(-\varepsilon,\varepsilon)$ is contained in the boundary of $A\cap\overline U$.
\qed
\endroster
\enddemo

\remark{Remark} 
In the proof of Proposition 2.2, the cases (a3), (b3) and (c3) deal with silvered points.
Let $p\!:\! \tilde U\to U$ denote the double cover
branched along $U\cap\Sigma$, so that $\tilde U$ is non-singular. In the three cases (a3), (b3) and
(c3),
$\overline {p^{-1}(A\cap U)}$ is a manifold with boundary, of dimension 1 or 2, whose interior is
totally geodesic in $\tilde U$. Moreover 
 $x\in\partial A$, but in cases (a3) and (c3) $p^{-1}(x)$ is an interior point of $p^{-1}(A\cap U)$.
This motivates the following definition.
\endremark

\definition{Definition} Let $A\subset E$ be a closed totally s-convex subset. The {\it
silvered boundary} $\partial_SA$ is the set of points $x\in\partial A\cap\Sigma$ that are silvered
and the following holds:
if $U\subset E$ is a neighborhood of $x$ and $p\!:\!  \tilde U\to U$ is the double cover branched
along $\Sigma\cap U$, then $p^{-1}(x)$ is an interior point of $p^{-1}(A\cap U)$.
We also define the \it non-silvered boundary $\partial_{NS}A$ \rm to be the set of points
in $\partial A$ that are not in the silvered boundary: $\partial_{NS}A=\partial A-\partial_SA$.
\enddefinition

\proclaim{Lemma 2.4} Let $A\subset  E$ be a non-empty   closed totally s-convex subset
in a Euclidean cone 3-manifold. Then the following hold:
\roster
\item"(i)" The non-silvered boundary $\partial_{NS}A$ is a closed subset of $E$.
\item"(ii)" If $\dim A=0$ or $3$,
               then $\partial_SA=\emptyset$ and $\partial A=\partial_{NS}A$.
\endroster
\endproclaim

\demo{Proof} Let $x\in A\cap \Sigma$ and let $U\subset E$ be a neighborhood of $x$.
By using the explicit description of  $U\cap A$ in the proof of  Proposition 2.2, 
we have that 
$x\in\partial_S A$ if and only if $x$ is a silvered point such that, either $U\cap A$ is a segment
orthogonal to
$\Sigma$ (case (a3)),
or $U\cap A$ has dimension 2 and $U\cap\Sigma$ is a piece of $\partial_S A\subseteq\partial
A$ (case (c3)).
This description of points in $\partial_S A$ implies that  $\partial_S A$ is open in $\partial A$,
hence assertion (i) is proved. We also deduce from the description that if $\partial_S 
A\neq\emptyset$ then
$\dim A=1$ or $2$, which is equivalent to assertion (ii).\qed
\enddemo

\proclaim{Proposition 2.5} Let $A\subset  E$ be a non-empty  closed totally s-convex subset
in a Euclidean cone 3-manifold. 
If $\dim A<3$ then every point  in $\Int (A)\cup\partial_SA=A-\partial_{NS}A$ has a neighborhood
$U\subset E$ isometric to the normal cone fibre bundle over $A\cap U$. More precisely:
\roster
\item"-" if $x\in\Int (A)$ then $U$ is isometric to the product $(A\cap U)\times B^c
(0,\varepsilon)$, where $ B^c
(0,\varepsilon)$ is a ball of radius $\varepsilon>0$ and dimension $c=\operatorname{codim}(A)$,
maybe with a singularity in its center.
\item"-" if $x\in\partial_SA$ and $p\!:\! \tilde U\to U$ is the double cover branched along
$\Sigma\cap U$, then $\tilde U$ is isometric to $p^{-1}(A\cap U)\times B^c
(0,\varepsilon)$, where $ B^c
(0,\varepsilon)$ is a non-singular ball of radius $\varepsilon>0$ and dimension
$c=\operatorname{codim}(A)$.
\endroster 
\endproclaim

\demo{Proof} If $x\in\Int(A)$ and $x\not\in\Sigma$ then the proposition is clear because  $\Int(A)$
is totally geodesic.  

If $x\in\Int(A)\cap\Sigma$  then, by the description given in 
the proof of Proposition 2.2, 
$A\cap U$ is either in case (a2) or in case (c1). In case (a2), $U\cap A$ is a
totally geodesic 2-dimensional disk perpendicular to $\Sigma$,
 and  $U$ is isometric to the product of
$U\cap A$ with an interval.
In case (c1), $U\cap A$ is a subinterval of
$\Sigma$ and $U$ is isometric to the product of $U\cap A$ 
with a singular 2-disk. Hence the
proposition holds in both cases.

When $x\in\partial_SA$,  $U\cap A$ is either in case (a3) or (c3). In both cases $p^{-1}(U\cap A)$
is a  totally geodesic submanifold of $p^{-1}(A)$ and the proposition follows. \qed 
\enddemo

The following proposition shows that
$A$ has a local supporting half-space at every point of
$\partial_{NS}A$.

\proclaim{Proposition 2.6} Let $A\subset E$ be  totally s-convex,
$x\in A$, $y\in\partial_{NS} A$ and  $\gamma\!:\! [0,l]\to A$  be a path
from $x$ to $y$ that realizes the distance $d(x,\partial_{NS} A)$. 
\roster
\item"(i)" If $y\in\Sigma$, then
$\gamma([0,l])\subset\Sigma$.
\item"(ii)" Let $B(y,{\varepsilon})$ be a standard ball of radius $\varepsilon>0$ and 
let $H\subset B( y,{\varepsilon})$ be the half-ball
bounded by the (maybe singular) totally geodesic disk orthogonal to
$\gamma$ at $y$. Then  $A\cap  B(y,{\varepsilon})\subseteq
H$.
 \endroster
\endproclaim

\demo{Proof} Let $y\in\Sigma\cap\partial_{NS}(A)$.
 We choose a neighborhood $U\subset E$ of $y$.
 By using the   description
in the proof of Proposition 2.2, the intersection $U\cap A$ is either
in case (b1), (b2) or (b3). 
It follows from this description  
that   if a path $\gamma\!:\! [0,l]\to A$
from $x$ to $y$ realizes $d(x,\partial_{NS} A)$, then it also
realizes $d(x ,D^2(y,\delta))$, where
$D^2(y,\delta)$ is the totally geodesic 2-disk of
radius
$\delta>0$ transverse to $\Sigma$.
 In particular, $\gamma$ is
orthogonal to $D^2(y,\delta)$ and   assertions (i) and (ii)
hold in the singular case.

If $y\in\partial_{NS}A$ is non-singular, then  assertion
(ii) can be proved by using  developping maps to reduce the
proof to the case of locally convex subsets in  $\Bbb R^3$.\qed
\enddemo

\proclaim{Lemma 2.7}
Let $E$ be a Euclidean cone 3-manifold and let $A\subset E$ be
totally s-convex. The function
$$\matrix\format\r&\,\c\,&\l\\
  \Phi\!:\! A & \to & \Bbb R\\
       x & \mapsto & d(x,\partial_{NS}A)
\endmatrix
$$
is concave (i.e. $-\Phi$ is convex). Moreover, if for some
geodesic path $\gamma\!:\! [0,l]\to A$ $\Phi\circ\gamma$ is
constant, then for every $t\in[0,l]$ there exists a
geodesic path orthogonal to $\gamma$ that realizes the
distance from
$\gamma(t)$ to
$\partial_{NS}A$.
\endproclaim

\demo{Proof} Let $\gamma\!:\! [0,l]\to A$ be a geodesic path, we
want to prove that $\Phi\circ\gamma$ is concave. For
$t\in[0,\l]$, let $\sigma_t\!:\! [0,\lambda]\to A$ be a
minimizing path from $\gamma(t)$ to $\partial A$.  Let
$\theta\in[0,\pi]$ be the angle between $\gamma'(t)$ and
$\sigma_t'(0)$. We claim that there exists some uniform $\varepsilon>0$
such that, if $\vert s\vert<\varepsilon$ and $t+s\in[0,l]$,
then
$$
  \Phi(\gamma(t+s))\leq \Phi(\gamma(t))-s \cos(\theta).
$$
This inequality shows that $\Phi\circ\gamma$ may be
represented locally as the infimum of linear functions,
and therefore $\Phi\circ\gamma$ is concave.

To prove this inequality, we consider first the case where 
$\sigma_t([0,\lambda])\cap \Sigma=\emptyset$. Let
$D^2(\sigma_t(\lambda),\delta)$ be a totally geodesic
disk with center $\sigma_t(\lambda)$ and radius $\delta>0$ 
that is orthogonal to $\sigma_t$. By Proposition 2.6 (ii), the
disk $D^2(\sigma_t(\lambda),\delta)$ bounds a 
locally supporting
half-ball for $A$. In particular there exists $\varepsilon>0$
such that, for
 $\vert
s\vert<\varepsilon$ 
$$
\Phi(\gamma(t+s))=d(\gamma(t+s),\partial_{NS}A)
\leq d(\gamma(t+s),D^2(\sigma_t(\lambda),\delta)). 
\tag*
$$
Moreover, by considering developping maps, we can use
elementary trigonometric formulas of Euclidean space to
conclude that for $\vert
s\vert<\varepsilon$
$$
d(\gamma(t+s),D^2(\sigma_t(\lambda),\delta))
=d(\gamma(t),\sigma_t(\lambda))- s\, \cos(\theta),
$$
where  $\varepsilon>0$ is small enough, so that the tubular
neighborhood
$\NN_{\varepsilon}(\gamma_t([0,l]))$ embeds in the Euclidean
space via developping maps (see Figure IV.5).
Furthermore,  the
parallel translation of $\sigma_t$ along $\gamma$ gives a
family of geodesic paths $\{\sigma_{t+s}\mid \vert
s\vert<\varepsilon\}$,  such that $\sigma_{t+s}$ has length
$\lambda -s\, \cos(\theta)$ and
minimizes the distance of
$\gamma(t+s)$ to
$D^2(\sigma_t(\lambda),\delta)$. Therefore, when
$
\Phi\circ\gamma$ is constant
we have
equality in \ttag*, $\theta=\pi/2$, and 
$\{\sigma_{t+s}\mid \vert
s\vert<\varepsilon\}$ is a family of geodesics of constant
length, orthogonal to $\gamma$ and which minimize the distance to
$\partial_{NS}A$.

\midinsert
 \centerline{\BoxedEPSF{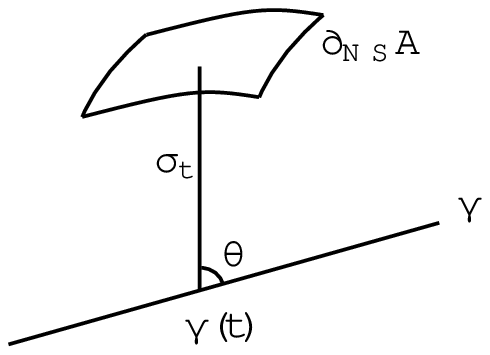 scaled 750}}
    \botcaption{Figure IV.5} 
    \endcaption
     \endinsert

When 
$\sigma_t([0,\lambda])\cap \Sigma\neq\emptyset$, since
$\sigma_t$ is minimizing, either
$\sigma_t([0,\lambda])\cap \Sigma=\{\gamma(t)\}$ or
$\sigma_t([0,\lambda])\subset \Sigma$  by Proposition 2.6 (i).
In particular, $\gamma(t)\in\Sigma$ and either 
$\gamma([0,l])\subset\Sigma$ or $\gamma(t)$ is an 
end-point of $\gamma$.
 Then the argument in the
non-singular case goes through in the singular case, by just
taking care when we use developping maps close to the singular
set.

Finally, note that a compactness argument allows to chose a
uniform
$\varepsilon>0$.
\qed
\enddemo

\head 3. Proof of the Soul Theorem. \endhead

\demo{Proof of Theorem 1.2 (Soul Theorem for Euclidean Cone
3-Manifolds)}
We start by considering Busemann
functions. We recall that a \it ray emanating from a point $p\in
C$ \rm is a continuous map $\gamma\!:\! [0,+\infty)\to E$ such
that the restriction on every compact subinterval is
minimizing. We assume that the rays are parametrized by
arc-length. The
\it Busemann function associated to $\gamma$ \rm is:
$$
b_{\gamma}(x)=\lim_{t\to+\infty} (t-d(x,\gamma(t))).
$$

By construction, Busemann functions are Lipschitz, with Lipschitz constant $1$.
In particular they are continuous.

A Euclidean cone 3-manifold with cone
angles less than $2\pi$ may be viewed as an Aleksandrov space
of curvature $\geq 0$. Next lemma is proved 
in \cite{\Yam, Proposition  6.2} for those Aleksandrov
spaces.

\proclaim{Lemma 3.1} Busemann functions on Euclidean cone
3-manifolds with cone angles
less than $2\pi$
are convex. \qed
\endproclaim

We consider a Euclidean cone 3-manifold $E$ and we fix a point
$x_0\in  E$. Following Cheeger and Gromoll \cite{\CGl}
or Sakai \cite{\Sak, Section V.3}, for $t\geq0$, we define
$$
A_t=\{ x\in E\mid b_{\gamma}(x)\leq t \text{ for every ray }
\gamma \text{ emanating from } x_0\}.
$$

\proclaim{Lemma 3.2}
For $t\geq0$, $A_t$ is a compact totally s-convex subset of
$E$, satisfying:
\roster
\item"1)" If $t_1\geq t_2\geq 0$, then $A_{t_2}\subseteq
A_{t_1}$ and $A_{t_2}=\{x\in A_{t_1}\mid d(x,\partial_{NS}
A_{t_1})\geq t_1-t_2\}$. In particular, for $t_2>0$, 
$\partial A_{t_2}=\partial_{NS}A_{t_2}=\{x\in A_{t_1}\mid d(x,\partial_{NS}
A_{t_1})= t_1-t_2\}$.
\item"2)" $E=\underset {t\geq 0} \to\bigcup A_t$.
\item"3)" $A_t$ intersects all connected components of $\Sigma$.
\endroster
\endproclaim

\demo{Proof}
The set $A_t$ is totally s-convex by Lemmas 2.1 and 3.1.
In order to prove the compactness, we suppose that there is a sequence of points
$x_n$ in $A_t$ going to infinity, and we will derive a
contradiction. For every $n\in\N$, consider the minimizing
path
$\gamma_n$ between $x_0$ and $x_n$, which is contained in $A_t$ by
convexity. Since the unit tangent bundle at $x_0$ is compact, the
sequence
$\gamma_n$ has a convergent subsequence to a ray $\gamma$
emanating from $x_0$ and contained in $A_t$; that
contradicts the definition of $A_t$.

We recall the following classical inequalities for Busemann functions,
which can be proved from triangle inequality. For every ray $\gamma$
emanating from $x_0$, every point $x\in E$ and every  real $t\geq 0$,
$$
  d(x,x_0)\geq b_{\gamma}(x)\geq t- d(x,\gamma(t)).
$$
In particular  $B(x_0,t)=\{x\in E\mid d(x,x_0)<t\}\subseteq
A_t$ for every
$t\geq 0$. Thus  assertion 2) is clear.
Moreover, for
$t>0$,
$\dim A_t=3$. Therefore, by Lemma 2.4,
$\partial A_t=\partial_{NS} A_t$.
 
To show  assertion 1) we prove first the inclusion $A_{t_2}\subseteq\{x\in A_{t_1}\mid
d(x,\partial_{NS} A_{t_1})\geq t_1-t_2\}$. Given $x\in A_{t_2}$ and
$y\not\in \Int (A_{t_1})$, we claim that $d(x,y)\geq t_1-t_2$. 
By hypothesis,
there exists a ray $\gamma$ emanating from $x_0$ such that 
$b_{\gamma}(x)\leq t_2$ and $b_{\gamma}(y)\geq t_1$.
For every 
$t>0$,
we have $d(x,y)\geq d(x,\gamma(t))-t+t-d(y,\gamma(t))$. By taking
the limit when $t\to+\infty$, we deduce that 
$
d(x,y)\geq b_{\gamma}(y)-b_{\gamma}(x)\geq t_1-t_2$,
as claimed.

To prove the reverse inclusion, we take a point $x\in A_{t_1}$ such that
$d(x,\partial_{NS} A_{t_1})\geq t_1-t_2$. We claim that for every ray $\gamma$
emanating from $x_0$, $b_{\gamma}(x)\leq t_2$. Note first that, for every
$t\geq t_1$, $d(x,\gamma(t))\geq t-t_1$, because $t_1\geq
b_{\gamma}(x)\geq t-d(x,\gamma(t))$. Let $z$ be the point in a
minimizing path between $x$ and $\gamma(t)$ such that
$d(z,\gamma(t))=t-t_1$. Then $z\not\in\Int A_{t_1}$, because 
$b_{\gamma}(z)\geq t-d(z,\gamma(t))=t_1$. It follows that $d(x,z)\geq
d(x,\partial A_{t_1})\geq 
t_1-t_2$ and
$$
 t-d(x,\gamma(t))=t-d(x,z)-d(z,\gamma(t))\leq t-(t_1-t_2)-(t-t_1)=t_2.
$$
By taking the limit when $t\to+\infty$, $b_{\gamma}(x)\leq t_2$. This proves assertion 1).

Finally assertion 3) follows from assertion 1) and the following
lemma:

\proclaim{Lemma 3.3} Let $A\subset E$ be totally s-convex
and let $A^r=\{x\in A\mid d(x,\partial_{NS}A)\geq r\}$.
If $A^r\neq\emptyset$ for some $r>0$ and $\Sigma_0$ is a component of
$\Sigma$ such that $A\cap \Sigma_0\neq\emptyset$, then 
$A^r\cap \Sigma_0\neq\emptyset$.
\endproclaim

\demo{Proof} By Lemma 2.7, $A^r$ is totally s-convex. Therefore $A^r$
intersects all components of $\Sigma$ having cone angle $\pi$. 
In general, 
let $\Sigma_0$ be a component of $\Sigma$ that intersects
$A$ and has  cone angle 
less than $\pi$. We distinguish two cases, according
to whether $\Sigma_0$
is a compact or not.

When $\Sigma_0\cong S^1$, for every $r>0$ either
$A^r\cap\Sigma_0=\emptyset$ or $\Sigma_0\subset A^r$,
because $\Sigma_0$ is a closed geodesic path and $A^r$ is
totally s-convex. Therefore, the distance to $\partial_{NS} A$
is constant on $\Sigma_0$. Let $r_0=d(\Sigma_0,\partial_{NS}
A)$, then $\Sigma_0\subseteq A^{r_0}$, because  $\Sigma_0$ intersects $A$. We claim that in fact
$A^{r_0}=\Sigma_0$.  Since $A^{r_0}$ is connected, we prove that $A^{r_0}=\Sigma_0$
by showing that there are no points in 
$A^{r_0}-\Sigma_0$ in a neighborhood of $\Sigma_0$.
Seeking a contradiction we suppose that 
 there is a point $y\in
A^{r_0}-\Sigma$ in a sufficiently small neighborhood of $\Sigma_0$.
Then, by Lemma 2.3 there is
a disk
$D^2(x,{\delta})\subset A^{r_0}$ of radius $\delta>0$, centered
at a point $x\in\Sigma_0$ and transverse to $\Sigma$. The convex
hull
of $\Sigma_0\cup D^2(x,{\delta})$ gives an open neighborhood
of $x$ in $A^{r_0}$, contradicting the fact that 
$d(x,\partial_{NS} A)=r_0$. This proves that $A^{r_0}=\Sigma_0$.
It follows that
$A^r=\emptyset$ for $r>r_0$, and $\Sigma_0\subseteq A^r$ for
$r\leq r_0$.

When $\Sigma_0\cong\Bbb R$, consider $r_0=\sup\{
d(x,\partial_{NS} A)\mid x\in\Sigma_0\cap A\}$. If $r_0=\infty$ then
there is nothing to prove, hence we can assume $r_0<\infty$. The
intersection $A^{r_0}\cap\Sigma_0$ is either a point or a
segment in $\Sigma_0$. If it is a segment, the argument
above for the closed case shows that $A^{r_0}\cap\Sigma_0=
A^{r_0}$. If the intersection is  a point then
$\dim(A^{r_0})=0$ or $2$. In any case $\dim A^{r_0}<\dim A$.
Thus $A^r=\emptyset$ for
$r>r_0$, and $A^{r}\cap\Sigma_0\neq\emptyset$ for
$r\leq r_0$.

This finishes the proofs of  Lemmas 3.2 and 3.3.
\qed
\enddemo
\enddemo

We fix  a value $t>0$ and we set $A=A_t$, where $A_t$ is defined as in Lemma
3.2. The subset $A\subset E$ is compact, totally s-convex
and $\dim A=3$. In particular
$\partial_{NS} A=\partial A\neq\emptyset$.

 For $r>0$, we
consider
$$
A^r=\{x\in A\mid d(x,\partial_{NS}A)\geq r\}
\quad\text{and}\quad A^{\max}=\bigcap\{ A^r\mid
A^r\neq\emptyset\}.
$$
If $A^r\neq\emptyset$, then $A^r$ is totally s-convex by
Lemma 2.7. Let $r_0=\max\{ d(x,\partial_{NS}A)\mid x\in
A\}$, then
$
A^{max}=A^{r_0}=\{x\in A\mid d(x,\partial_{NS}A)= r_0\}
$.
By Lemma 2.7, every geodesic in $A^{max}$ is perpendicular
to a geodesic minimizing the distance to $\partial_{NS} A$,
hence $\dim A^{max}<\dim A=3$. 

We set $A(1)=A^{max}$. If $\partial_{NS} A(1)=\emptyset$
or $\dim A(1)=0$, then we take $S=A(1)$. Otherwise, we construct
$A(2)=A(1)^{max}$ and so on. Since $\dim A(i+1)<\dim A(i)$,
this process stops and we obtain $S=A(i)$ for either $i=1$,
$2$ or $3$. Thus $S$ is a compact totally s-convex subset
 with $\partial_{NS}S=\emptyset$ and $\dim S<3$.

Next we prove that $E$ is isometric to the normal cone fiber bundle of $S$. The
key point of the proof is  the following lemma:

\proclaim{Lemma 3.4} Each point of $E-S$ has a unique minimizing
geodesic path to $S$. Moreover for singular points, this path is
contained in $\Sigma$.
\endproclaim

\demo{Proof}
We consider the following subset of $E-S$
 $$
X=\left\{x\in E-S\,\left\vert 
\matrix\format \l\\
 x \text{ has more than one minimizing path to }
S, \text{ or }\\ x\in \Sigma \text{ has a minimizing path to }
S 
\text{ not contained in }\Sigma
\endmatrix \right.
\right\}
$$

\proclaim{Claim 3.5} $X$ is a closed subset of $E-S$ and $d(X,S)>0$.
\endproclaim

\demo{Proof of the claim} By Lemma 1.3 of Chapter III, each point has a  finite
number of minimizing paths to the totally s-convex submanifold $S$.
Moreover, given a point $x\in E-S$, in a sufficiently small neighborhood of $x$ the
minimizing paths to $S$ are obtained by perturbation of those of $x$.
It follows that the property of having a unique minimizing path to $S$,
 contained in $\Sigma$ for singular points, is an open
property in $E-S$, and thus $X\subset E-S$ is closed.

Since $S$ is compact, totally s-geodesic and 
$\partial_{NS}(S)=\emptyset$, by Proposition 2.5
$S$ has a metric tubular neighborhood. Thus
$d(X,S)>0$.
\qed
\enddemo

We come back to the proof of Lemma 3.4. Seeking
a contradiction, we assume $X\neq\emptyset$.

Let $x_0\in X$ be such that
$d(x_0,S)=d(X,S)$. Either $x_0$ has two minimizing paths to $S$ or
$x_0\in\Sigma$ has a minimizing path  to S not in $\Sigma$.

We assume first that there are  two
minimizing paths $\gamma_1$ and $\gamma_2$ from $x_0$ to $S$. If at $x_0$ the angle
$\angle(\gamma_1,\gamma_2)<\pi$, then $\gamma_1$ and
$\gamma_2$ can be deformed to shorter paths, with common
origine, that minimize the distance to $S$, contradicting the definition of $x_0$. Therefore
$\angle(\gamma_1,\gamma_2)=\pi$ and $\gamma_1^{-1}*\gamma_2$ is a
geodesic with end-points in $S$, contradicting the fact that
$S$ is totally s-geodesic.

Now we assume that $x_0\in\Sigma$ has a minimizing path $\gamma$ to 
$S$ not contained in $\Sigma$. If at $x_0$ the angle 
$\angle(\gamma,\Sigma)<\pi/2$, then we can perturb
$\gamma$ to a shorter path, not contained in $\Sigma$ and minimizing the distance of a point in
$\Sigma$ to $S$. This contradicts the definition of $x_0$, hence we can
 assume that the angle $\angle(\gamma,\Sigma)=\pi/2$. 
We remark that the point $x_0$ is not silvered, otherwise $x_0\in S$,
because $\gamma^{-1}*\gamma$ would be a s-geodesic path with end-points in $S$.
By
construction there exists a totally s-convex subset $A\subset E$ such
that $S\subset \Int(A)$ and $x_0\in\partial_{NS}A$ (this is one of the sets
constructed in Lemma 3.2 and 3.3, i.e. $A=A_t$ for some $t>0$ or $A=A(i)^r$ for
some
$r>0$ and $i=1$, $2$ or $3$).
Since the angle
$\angle(\gamma,\Sigma)=\pi/2$, Proposition 2.6 says that $\gamma$ is
tangent to the boundary of a supporting half-space for $A$ at the point
$x$. In particular $\gamma$ does not go to the interior of $A$. This
contradicts  the fact that $S\subset \Int (A)$, because
$d(S,\partial_{NS} A)>0$.

This finishes the proof of Lemma 3.4.
\qed \enddemo

We have shown that every point in $E-S$ has a unique minimizing path to
$S$, and that for singular points this path is  contained in
$\Sigma$. 
By proposition 2.5, since $S$ is compact, totally s-convex and 
$\partial_{NS}(S)=\emptyset$, for some 
 $\varepsilon>0$ the tubular
neighborhood $\NN_{\varepsilon}(S)$ 
is isometric to the normal cone fiber bundle of $S$.
The fact that every point in $E-S$ has a unique minimizing
geodesic to $S$ implies that the radius $\varepsilon$ 
of the tubular neighborhood
can be taken
arbitrarily large, and  thus $E$ is isometric to the normal cone fiber bundle of $S$.
\qed

\enddemo

\head 4. Proof of the Local Soul Theorem.
\endhead

\demo{Proof}
Seeking a contradicition, we assume 
that there exist some
$\varepsilon_0>0$ and some
$D_0>1$ such that for every $\delta>0$ and every $R>D_0$
there are a hyperbolic cone 3-manifold $C_{\delta,R}$ with
cone angles in $[\omega,\pi]$ and a point
$x\in C_{\delta,R}$ with $\Inj(x)<\delta$ that does not verify the
statement of the Local Soul  Theorem with parameters
$\varepsilon_0$, $D_0$ and $R$. By
taking
$\delta=\frac1n$ and
$R=n$, we obtain a sequence of pointed hyperbolic cone
3-manifolds
$(C_n, x_n)_{n\in\N}$ such that $\Inj(x_n)<\frac1n$ and
$x_n$ does not verify the Local Soul Theorem with
parameters $\varepsilon_0$, $D_0$ and $R=n$.

We apply the Compactness Theorem (Chapter III) to the  sequence of rescaled
pointed cone 3-manifolds
$(\bar C_n, \bar x_n)_{n\in\N}=(\frac1{\Inj(x_n)}C_n,
x_n)_{n\in\N}$. Then a subsequence of $(\bar C_n,\bar x_n)$
converges to a pointed Euclidean cone 3-manifold
$(C_{\infty},x_{\infty})$.

If the limit $C_{\infty}$ is compact, then the geometric
convergence implies that for some integer 
$n_0$  there exists a $(1+\varepsilon_0)$-bilipschitz
homeomorphism $f\!:\! C_{\infty}\to \bar C_{n_0}$.
We can also choose $n_0$ such that 
$$
n_0>\Diam (C_{\infty}). 
$$
Then the rescaled Euclidean cone  3-manifold $E=\Inj (x_{n_0})
C_{\infty}$ is $(1+\varepsilon_0)$-bilipschitz homeomorphic
 to $C_{n_0}$ and we have
$$
\Diam (E)< n_0\Inj(x_{n_0})=R\,\Inj(x_{n_0}).
$$
Hence $x_{n_0}$ satisfies the statement of the Local Soul
Theorem in the compact case and we get a contradiction.

If the limit $C_{\infty}$ is not compact, then by the Soul Theorem 1.2, 
$C_{\infty}$ has a soul
$S_{\infty}$ and is the normal cone bundle of $S_{\infty}$.
Since
$\Inj(x_{\infty})\leq 1$, the soul $S_{\infty}$ has
dimension $1$ or $2$. We choose a real number $\nu_{\infty}$
verifying
$$
\nu_{\infty}> D_0\,\max\big(\Diam (S_{\infty}),
d (x_{\infty},S_{\infty})+
(1+\varepsilon_0)\varepsilon_0, 1\big).
$$
For $n_0$ sufficiently large, the geometric convergence implies  
the existence of a
$(1+\varepsilon_0)$-bilipschitz embedding
$\bar g\!:\! \NN_{\nu_{\infty}}(S_{\infty})\to \bar C_{n_0}$ such
that
$d(\bar g(x_{\infty}),x_{n_0})<\varepsilon_0$ and $n_0\geq
\nu_{\infty}$. The image $U=\bar
g(\NN_{\nu_{\infty}}(S_{\infty}))$ is an open neighborhood 
of $x_{n_0}$, because
$$
\multline
d(\partial U,\bar g(x_{\infty}))\geq\frac1{1+\varepsilon_0}
d(\partial\NN_{\nu_{\infty}}(S_{\infty}),x_{\infty}) \geq
\\
\frac1{1+\varepsilon_0} 
(\nu_{\infty}-d(x_{\infty},S_{\infty})  )
> 
 \varepsilon_0
>
d(\bar g(x_{\infty}),x_{n_0}).
\endmultline
$$

 As in the compact case, we
consider the rescaled 3-manifold
$E=\Inj(x_{n_0})\,  C_{\infty}$ with soul 
$S=\Inj(x_{n_0})\, S_{\infty}$. By taking
$\nu=\Inj(x_{n_0})\,\nu_{\infty}$, $\bar g^{-1}$ induces a
$(1+\varepsilon_0)$-bilipschitz homeomorphism
$f\!:\! U\to \NN_{\nu}(S)$. Moreover, the constants have been
chosen so that $\nu\leq n_0\Inj(x_{n_0})\leq 1$ and
$$
\max\big(\Inj(x_{n_0}),d(f(x_{n_0}),S),\Diam(S)\big)\leq
\nu/D_0.
$$
Thus $x_{n_0}$ verifies the statement of the Local Soul Theorem in the non-compact
 case and we obtain a contradiction  again. This finishes the proof of the Local Soul Theorem.
\qed
\enddemo

\proclaim{Corollary 4.1} 
Let $(C_n)_{n\in\N}$ be a sequence of hyperbolic cone
3-manifolds such that $\sup\{\Inj(x)\mid x\in C_n\}$ converges
to zero when
$n\to\infty$. Then for every $\varepsilon>0$ and $D>1$, 
\roster
\item"--" either
there exists $n_0= n_0(\varepsilon,D)$ such that for $n\geq
n_0$ the Local Soul Theorem with parameters $\varepsilon,D$
applies to every point of $C_n$ and the local models are
non-com\-pact;
\item"--" or, after rescaling, a
subsequence of  $(C_n)_{n\in\N}$  converges to a
closed Euclidean cone 3-manifold.
\endroster
\endproclaim

\demo{Proof} Since the supremum on $C_n$ of the injectivity
radius goes to zero, there exists an integer $n_0>0$ such that, for $n>n_0$,
the Local Soul Theorem
applies to every point of $C_n$ with compact or non-compact
models. 
 We assume the existence of a sequence of points
$(x_{k})_{k\in\N}$ such that $x_k\in C_{n_k}$, the local
model of $x_k$ is compact, and $n_k\to\infty$ as
$k\to\infty$. In particular, for every $k\in\N$ there exists a
closed Euclidean cone 3-manifold $E_k$ and a
$(1+\varepsilon)$-bilipschitz homeomorphism $f_k\!:\!  E_k\to
C_{n_k}$ such that $\Diam(E_k)\leq R\,\Inj(x_k)$. Therefore,
$$
\Diam(C_{n_k})\leq (1+\varepsilon) \Diam(E_k)
\leq (1+\varepsilon) R\,\Inj(x_k).
$$
Hence, the diameter of the rescaled cone 3-manifold 
$\bar C_{n_k}=\frac1{\Inj({x_k})}C_{n_k}$ is uniformly
bounded above. By the Compactness Theorem (Chapter III)  $(\bar
C_{n_k},x_{n_k})_{k\in\N}$ has a subsequence converging to a
pointed Euclidean cone 3-manifold $(E,x_{\infty})$. Moreover, the limit
$E$ is compact, because the diameter of $\bar C_{n_k}$ has a
uniform upper bound.
\qed
\enddemo

\newpage
\rightheadtext{ }
\leftheadtext{V \qquad    Sequences with cone angles less than $\pi$}

\

\centerline{\smc chapter \  v} 

\

\centerline{\chapt SEQUENCES  \,  OF \,  CLOSED \,  HYPERBOLIC \,  CONE}

\vglue.2cm

\centerline{\chapt  3-MANIFOLDS \, WITH \,  CONE \,  ANGLES \,  LESS \,  THAN \,   $\pi$}

\

\

This chapter is devoted to the proof of the following theorem.

\proclaim{Theorem A} 
Let $(C_n)_{n\in\N}$ be a sequence of closed hyperbolic  cone 3-manifolds
with fixed topological type $(C,\Sigma)$
such that the cone angles increase and are contained in  
$[\omega_0,\omega_1]$, with
$0<\omega_0<\omega_1<\pi$.
Then there exists a subsequence $(C_{n_k})_{k\in\N}$ such that one of 
the following occurs: 
\roster
  \item"1)" The  sequence $(C_{n_k})_{k\in\N}$ converges 
geometrically to a hyperbolic cone 3-man\-ifold with topological type
$(C,\Sigma)$ whose cone angles are the limit of the cone angles of
$C_{n_k}$.
   \item"2)" For every $k$, $C_{n_k}$ contains an  embedded $2$-sphere $S_{n_k}^2\subset
C_{n_k}$
that intersects $\Sigma$ in three points, and the sum of the three cone
angles at $S_{n_k}^2\cap\Sigma$ converges to $2\pi$.
   \item"3)" There is a sequence of positive  reals $\lambda_k$ approaching
$0$ such that the subsequence of rescaled cone 3-manifolds
$(\lambda_k^{-1}C_{n_k})_{k\in\N}$ converges geometrically to a 
Euclidean cone 3-manifold of topological type $(C,\Sigma)$ and whose
cone angles are the limit of the cone angles of $C_{n_k}$.
\endroster
\endproclaim

The proof of Theorem A  splits into two cases, according to whether the
sequence
 $(C_n)_{n\in\N}$ collapses or not.

\definition{Definition} We  say that a  sequence $(C_n)_{n\in\N}$ of
cone 3-manifolds  \it collapses  \rm if the sequence $(\sup\{\Inj(x)\mid
x\in C_n\})_{n\in\N}$ goes to zero.
\enddefinition

\head 1. The non-collapsing case
 \endhead

The following proposition (see \cite{\Zhou}, \cite{\SOK} and
\cite{\Hod}) proves Theorem A when the sequence $(C_n)_{n\in\N}$ does
not collapse.

\proclaim{Proposition 1.1}  Let $(C_n)_{n\in\N}$ be a sequence of
hyperbolic cone 3-manifolds satisfying the hypothesis of Theorem A. If
the sequence $(C_n)_{n\in\N}$ does not collapse, then there is a
subsequence $(C_{n_k})_{k\in\N}$ that verifies assertion {\rm 1)} or
{\rm 2)} of Theorem A.
\endproclaim

\demo{Proof} Since the sequence $C_n$ does  not collapse, after passing to a subsequence
if necessary, there is a positive real number $a>0$ and, for every $n\in\N$,  a
point $x_n\in C_n$ such that
$\Inj(x_n)\geq a$. Thus the sequence  $(C_n,x_n)_{n\in\N}$ is
contained in  $\CC_{[\omega_0,\omega_1],a}$,
the space of pointed cone 3-manifolds $(C,x)$ with
constant curvature in  $[-1,0]$,  cone angles in $[w_0,w_1]$, and
 such that $\Inj(x)\geq a>0$.
 By the Compactness Theorem (Chapter III), the sequence
$(C_n,x_n)_{n\in\N}$ has a convergent subsequence, which we denote
again by $(C_n,x_n)_{n\in\N}$.
 Hence, we can assume that the sequence
 $(C_n,x_n)_{n\in\N}$  converges geometrically to a pointed hyperbolic
cone 3-manifold
$(C_{\infty},x_{\infty})$, which may be compact or not.

If the limit cone 3-manifold $C_{\infty}$ is compact, then the  geometric
convergence implies that $C_{\infty}$ has the same topological type $(C,\Sigma)$
as the cone 3-manifolds of the sequence $C_n$. Moreover the  cone angles of
$C_{\infty}$ are the limit of the  cone angles of $C_n$. This shows
that in this case the assertion 1) of Theorem A holds. If the limit cone
3-manifold is not compact, then the next proposition shows that we get
the assertion 2) of Theorem A.
 
 \proclaim{Proposition 1.2} If the limit cone 3-manifold $C_{\infty}$ is
not compact, then for $n$ sufficiently large, $C_n$ contains an embedded
$2$-sphere $S_{n}\subset C_n$ that intersects $\Sigma$ in three points, and the sum of the three
cone angles at $S_{n}\cap\Sigma$ converges to $2\pi$. 
\endproclaim

\rightheadtext{V \qquad    Sequences with cone angles less than $\pi$}

We start with the following lemma.

 \proclaim{Lemma 1.3} The limit cone 3-manifold $C_{\infty}$ has finite
volume.
\endproclaim

\demo{Proof of Lemma 1.3} 
 Since $\Vol(C_{\infty})=\lim\limits_{R\to\infty}\Vol\big(
B(x_{\infty},R)\big)$, it suffices to bound
$\Vol\big(B(x_{\infty},R)\big)$ independently of $R$. From the
geometric convergence, for every $R>0$ there is $n_0$ so that, for
$n>n_0$, there exists a
$(1+\varepsilon_n)$-bilipschitz embedding $f_n:B(x_{\infty},R) \to (C_n,x_n)$,
with 
$\varepsilon_n\to 0$.
Hence, for $R>0$ and $n>n_0$,  $\Vol(B(x_{\infty},R))\leq
(1+\varepsilon_n)^3\Vol(C_n)$, and we get the  bound
$\Vol(B(x_{\infty},R))\leq 2^3\Vol(C_n)$. According to Schl\"afli's formula for
cone 3-manifolds (cf. \cite{\Hod} or \cite{\PoTwo}), since the cone angles of the 3-manifolds $C_n$
increase, the sequence $(\Vol(C_n))_{n\in\N}$ decreases, hence it is  
bounded above. This fact follows from \cite{\PoTwo, Prop. 4.2}, since the sequence
$(\rho_n)$ of  holonomy representations of the cone 3-manifolds $C_n$ belongs to
a piecewise analytical path in the variety of representations.
  \qed
  \enddemo

Proposition 1.2 follows from the next one.

\proclaim{Proposition 1.4} If the limit cone 3-manifold $C_{\infty}$ is not compact, then
neither is its singular set $\Sigma_{\infty}$.
\endproclaim

 \demo{Proof of Proposition 1.2 from proposition 1.4}
 From 1.4, there is a connected component $\Sigma^0_{\infty}$ of
$\Sigma_{\infty}$ which is not compact.
 Since $\Vol(C_{\infty})$ is finite by Lemma 1.3, the cone-injectivity radius
along 
$\Sigma^0_{\infty}$ is not bounded away from zero.

 By the Local soul theorem (Chapter IV), there is a point
$y\in\Sigma^0_{\infty}$ having a neighborhood
$(1+\varepsilon)$-bilipschitz homeomorphic to a product
$S^2(\alpha,\beta,\gamma)\times(-\nu,\nu)$, where $\nu>0$ and 
$S^2(\alpha,\beta,\gamma)$ is a two-dimensional Euclidean cone 3-manifold
with underlying space the sphere $S^2$ and singular set three cone
points with singular angles 
$\alpha$, $\beta$ and $\gamma$ such that $\alpha+\beta+\gamma=2\pi$.

 Since $(C_n,x_n)$ converges geometrically to
$(C_{\infty},x_{\infty})$, there is an
 $(1+\varepsilon_n)$-bilip\-schitz embedding
$f_n:S^2(\alpha,\beta,\gamma)\times\{0\}\to C_n$, with
$\lim\varepsilon_n=0$. The image $S^2_n=f_n(S^2(\alpha,\beta,\gamma))$ is
a 2-sphere  embedded in $C_n$ that intersects the
singular set in three points and the sum of the three cone
angles $\alpha_n+\beta_n+\gamma_n$  at $S^2_{n_k}\cap\Sigma$ converges to $2\pi$. 
This proves Proposition 1.2.\qed
 \enddemo

The remaining of this section is devoted to the proof of Proposition 1.4.

\demo{Proof of Proposition 1.4}  Seeking a contradiction, we
 suppose that the limit cone
3-manifold $C_{\infty}$ is not compact, but that its
singular set $\Sigma_{\infty}$ is compact.  We use the following lemma
of \cite{\KoOne} (see also \cite{\Zhou}).

\proclaim{Lemma 1.5} Let $C_{\infty}$ be a hyperbolic  cone 3-manifold
of finite volume whose singular set $\Sigma_{\infty}$ is compact. Then
$C_{\infty} - \Sigma_{\infty}$ admits a complete hyperbolic structure
of finite volume.
\endproclaim

\demo{Proof} The proof consists in deforming the metric in 
$\NN_{\varepsilon}(\Sigma_{\infty})-\Sigma_{\infty}$, where
$\NN_{\varepsilon}(\Sigma_{\infty})$ is a tubular neighborhood
of radius $\varepsilon>0$, so that $C_{\infty} - \Sigma_{\infty}$ 
admits a complete metric of (non-constant) sectional curvature 
$K\leq -a^2<0$. 
With this complete metric $C_{\infty}-\Sigma_{\infty}$ has a finite volume,
therefore by \cite{\Ebe, Thm. 3.1} it has only finitely many ends and each end is parabolic.
In particular $C_{\infty}-\Sigma_{\infty}$ is the interior of a compact
manifold with toral boundary. Since strictly negative curvature forbides
essential spheres and tori as well as Seifert fibrations,  
Thurston's Hyperbolization Theorem (for Haken manifolds) provides a complete
hyperbolic structure on $C_{\infty} - \Sigma_{\infty}$.
See \cite{\KoOne} for the details of the deformation. \qed
\enddemo

\remark{Remark} If $\Sigma_{\infty}$ is compact, then the ends of $C_{\infty}$
are cusps \cite{\Ebe, Thm. 3.1}. In particular the ends are topologycally $T^2\times\Bbb [0,\infty)$.
Moreover, if $\rho_{\infty}:\pi_1(C_{\infty} - \Sigma_{\infty})\to PSL_2(\Bbb C)$
is the holonomy of $C_{\infty}$, then  the restriction of $\rho_{\infty}$
to $\pi_1(T^2\times\{0\})$ is parabolic and faithfull. This is a consequence
of the fact that, in the proof of Lemma 1.5, the metric has not been changed on
the ends.
\endremark

Let $N_{\infty}\subset C_{\infty} - \Sigma_{\infty}$ be  a compact
core containing the base point $x_{\infty}$. If $\Sigma_{\infty}$ is
compact, then the boundary $\partial N_{\infty}$ is a collection of
tori $T_1,\ldots,T_p$ and:
$$
 C_{\infty} -\Sigma_{\infty}=
N_{\infty} \underset{\partial N_{\infty}}\to\cup
        \bigsqcup_{i=1}^p T_i\times [0,\infty).
$$
We set $X=C-\NN(\Sigma)=C_n-\NN(\Sigma_n)$, where $\NN$ denotes an open
tubular neighborhood. From the geometric convergence (for
$n$ sufficiently large) there is an $(1+\varepsilon_n)$-bilipschitz embedding
$f_n:N_{\infty}\to C_n $,  with $\varepsilon_n\to 0$, such that 
$$
N_n=f_n(N_{\infty})\subset
C_n-\NN(\Sigma_n)=X.$$

\proclaim{Claim 1.6} For $n$ sufficiently large, every connected
 component of $X-\Int(N_n)$ is
either a solid torus $S^1\times D^2$ or a product $T^2\times[0,1]$.
\endproclaim

\demo{Proof} First we show that $X-\Int(N_n)$ is irreducible for $n$
 sufficiently large. Otherwise, after
 passing to a subsequence, we can assume that $X-\Int(N_n)$ is
reducible  for every $n$. This implies that there is a ball $B_n\subset X$ such that $N_n\subset
B_n$, because $X$ is irreducible. Let
  $ \rho_n:\pi_1(X,x_n)\to PSL_2(\Bbb C)$ and
  $ \rho_{\infty}:\pi_1(N_{\infty},x_{\infty})\to PSL_2(\Bbb C)$
denote the holonomy representations of $C_n$ and $C_{\infty}$ respectively. The
geometric convergence implies the algebraic convergence of the holonomies (Proposition
5.4, Chapter III). This means that for every $\gamma\in\pi_1(N_{\infty},x_{\infty})$,
$\rho_n(f_{n*}(\gamma))$ converges to $\rho_{\infty}(\gamma)$. Since
$N_n$ is contained in a ball,  $f_{n*}(\gamma)=1$, so $\rho_{\infty}(\gamma)=1$ for
every $\gamma\in\pi_1(N_{\infty},x_{\infty})$.
 Since the holonomy representation of
$C_{\infty}$ is non-trivial, we get a contradiction. This proves the irreducibility
of  $X-\Int(N_n)$. 

Since $X-\Int(N_n)$ is irreducible and $\partial N_n$ is a collection of tori, the claim 
follows easily from the
fact that $X$ is irreducible and atoroidal. 
\qed
 \enddemo

In order to get a contradiction with the hypothesis that $\Sigma_{\infty}$ is compact we
 need in addition the following claim.

\proclaim{Claim 1.7} For $n$ sufficiently large, at least one component $X-\Int(N_n)$ is
 a solid torus.
\endproclaim

\demo{Proof} We assume that the claim is false and look for a contradiction. Thus,
after passing  to  a subsequence if necessary, we can assume that all the components of $\partial
N_n$ are parallel to the boundary of $\partial X$; this means that $f_n:N_{\infty}\to
X$ is a homotopy equivalence. 

 If $T^2\subset\partial N_{\infty}$ is a component corresponding to an end of
$C_{\infty}$, then the image $\rho_{\infty}(\pi_1(T^2,x_{\infty}))$ is a
parabolic subgroup of $PSL_2(\Bbb C)$ by the remark following Lemma 1.5.
Furthermore, since $C_n$ converges geometrically to
$C_{\infty}$, for every $\gamma\in \pi_1(T^2,x_{\infty})$,
$
    \rho_{\infty}(\gamma)=\lim\limits_{n\to\infty}\rho_n(f_{n*}(\gamma)).
$

Since $X$ has a complete hyperbolic structure, by Mostow's rigidity theorem \cite
{\Mos}  and Whaldhausen's theorem \cite {\Wal} the group
$\pi_0(\operatorname{Diff}(X))$ is finite (see also \cite {\Joh}). Hence, after passing to a
subsequence, we can choose $\gamma\in
 \pi_1(T,x_{\infty})$ such that, for every $n$ , $f_{n*}(\gamma)$ is
conjugate to
a meridian $\mu_0$  of a fixed component $\Sigma_0$ of $\Sigma$. Since
$\mu_0$ is elliptic, $\Trace(\rho_nf_{n*}(\gamma))=\pm 2
\cos(\alpha_n/2)$, where
 $\alpha_n\in[\omega_0,\omega_1]$ is the cone angle of the manifold
$C_n$ at the component
$\Sigma_0$.
Since $0<\omega_0<\omega_1\leq\pi$ the sequence 
 $\vert\Trace(\rho_n(f_{n*}(\gamma)))\vert$ is bounded away from $2$.
As $\rho_{\infty}(\gamma)$ is parabolic,
 $\vert\Trace(\rho_{\infty}(\gamma))\vert=2$, and we obtain
a contradiction with the convergence of $\rho_n(f_{n*}(\gamma))$ to
$\rho_{\infty}(\gamma)$.
\qed
\enddemo

 From  the Claim 1.7, there is a collection $T_1,\ldots,T_q$  of components
of $\partial N_{\infty}$  such that, for $n$ sufficiently large, $f_n(T_i)$ bounds
a solid torus $V^i_n\subset X$, for $i=1,\ldots,q$.

 Let $\lambda^i_n\subset f_n(T_i)$ be the boundary of a meridian of
the solid torus
$V^i_n$, for $i=1,\ldots,q$. The inverse images 
$\tilde\lambda^1_n=
f_n^{-1}(\lambda^1_n),\ldots,\tilde\lambda^q_n=f_n^{-1}(\lambda^q_n)$ are the
meridians of the Dehn fillings of $N_{\infty}=f_n^{-1}(N_n)$ which give $X$.
More precisely,
$$
X=N_{\infty}\bigcup_{\phi_{i,n}}\bigsqcup_{i=1}^nS^1\times D^2_i,
$$
 where, for $i=1,\ldots,q$, the gluing maps
$\phi_{i,n}:S^1\times \partial D^2_i\cong T_i\subset\partial N_{\infty}$
satisfy $\phi_{i,n}(\{*\}\times\partial D_i^2)=\tilde
\lambda^i_n$.

We have now the following claim:

 \proclaim{Claim 1.8} For every $i=1,2,\ldots,q$, the sequence of
simple closed curves
 $(\tilde\lambda^i_n)_{n\geq n_0}$ represents infinitely many distinct elements in
$H_1(T_i)$. Hence, after passing to a subsequence, the lenght of $\tilde\lambda^i_n$
goes to infinity with $n$.
 \endproclaim

 \demo{Proof} If the claim is not true, then, after passing to  a
subsequence, we can assume that there is an index $i\in\{
1,2,\ldots,q\}$ such that the curves
 $\tilde\lambda^i_n$ are all homotopic to a fixed curve
$\tilde\lambda^i$, for every $n$.

 Since the sequence $(C_n,x_n)$ converges geometrically to $(C_{\infty},x_{\infty})$ we
have: 
$$
    \rho_{\infty}(\tilde\lambda^i)=
          \lim_{n\to\infty}\rho_n(f_{n*}(\tilde\lambda^i))=\pm\operatorname{Id}
$$
because, for $n$ sufficiently  large, $f_{n}(\tilde\lambda^i)=\lambda^i_n$ bounds a meridian disk of
a solid torus $V^i_n$. 
Since $T_i$ corresponds to a cusp of $C_{\infty}$, the holonomy 
$\rho_{\infty}(\tilde\lambda^i)$ is not trivial and we get a 
contradiction.
\qed
\enddemo

 We are now ready to  contradict the hypothesis that
$\Sigma_{\infty}$ is compact.
 If $\Sigma_{\infty}$ is compact, then, by Claims 1.6, 1.7 and 1.8, 
we have, for $i=1,\ldots,q$, a sequence of 
curves $(\tilde\lambda^i_n)_{n\geq n_0}$ in $T_i\subset\partial N_{\infty}$ whose
lengths go to infinity with n, and so that the 3-manifold obtained by Dehn filling
with meridians $\{\tilde\lambda^1_n,\ldots,\tilde\lambda^q_n\}$
is always
$X=C-\NN(\Sigma)$. 
According to Thurston's Hyperbolic Dehn filling Theorem \cite{\ThuNotes}, (cf.
Theorem 1.1, Chapter II), almost all these Dehn fillings are hyperbolic. Furthermore,
by Schl\"afli's formula, almost all of them have different hyperbolic volumes. Thus
we get a contradiction, because our Dehn fillings give always the same 3-manifold
$X$. This finishes the proof of Propositions 1.4 and  1.2.
\qed
 \enddemo
\enddemo

\head 2. The collapsing case
 \endhead

The next proposition proves Theorem A when the sequence of hyperbolic cone
3-man\-ifolds
$C_n$ collapses.

\proclaim{Proposition 2.1}   Let $(C_n)_{n\in\N}$ be a sequence of
hyperbolic cone 3-manifolds with the same hypothesis as in Theorem A. If
the sequence $(C_n)_{n\in\N}$ collapses, then there is a subsequence
$(C_{n_k})_{k\in\N}$ that satisfies assertions {\rm 2)} or {\rm 3}) of
Theorem A.
\endproclaim

The proof uses Gromov simplicial volume of a compact oriented 3-manifold $M$ and
the dual notion of real bounded cohomology of $M$, both introduced by M. Gromov
\cite{\Gro} (see also \cite{\Iva}).

The simplicial volume $\norm M$ of a compact, orientable, 3-manifold
$M$, with boundary $\partial M$ (possibly empty) is defined as follows
 $$
\norm M =\inf\left\{  \sum_{i=1}^n \vert\lambda_i\vert \ \Bigg\vert\ 
\matrix \format\l\\
     \sum\limits_{i=1}^n \lambda_i\sigma_i \text{ is a cycle representing a fundamental}\\
     \text{class in } H_3(M,\partial M;\Bbb R),\text{ where }\sigma_i:\Delta^3\to M \\
     \text{is a singular simplex and } \lambda_i\in\Bbb R,\ i=1,\ldots,n.
\endmatrix
\right\}
$$
In particular, when $C$ is a closed and oriented 3-manifold and $\Sigma\subset C$ is a
link, we define the simplicial volume  $\norm{C-\Sigma}=\norm{C-\Int (\NN(\Sigma))}$,
where $\NN(\Sigma)$ is a tubular neighborhood of $\Sigma$ in $C$.

 We are starting now to prove Proposition 2.1.

\demo{Proof of Proposition 2.1} We are going to show that if Assertions 2) and 3) of Theorem A
do not hold, then the simplicial volume $\norm{C-\Sigma}$ is zero, and this would contradict
the hyperbolicity of $C-\Sigma$ (Lemma 1.5). To show that the simplicial volume
vanishes, we need for a subset of $C$ the notion of abelianity in $C-\Sigma$.

\definition{Definition} We say that a subset $U\subset C$ is \it abelian in \rm
$C-\Sigma$ if the image $i_*(\pi_1(U-\Sigma))$ is an abelian subgroup of
$\pi_1(C-\Sigma)$, where $i_*$ is the morphism induced by the inclusion
$i:(U-\Sigma)\to (C-\Sigma)$.
  \enddefinition

\definition{Definition} Let $C$ be a hyperbolic cone 3-manifold,  $x\in C$, and
$\varepsilon, D>0$. An \it $(\varepsilon,D)$-Mar\-gu\-lis' neighborhood  of abelian
type \rm  of $x$ is a neighborhood $U_x$ $(1+\varepsilon)$-bilips\-chitz homeomorphic to
the normal cone fiber bundle
$\NN_{\nu}(S)$, of radius $\nu\leq 1$, of the soul $S$ of one of the
following non-compact orientable Euclidean cone 3-manifolds: 
$$
   T\times\Bbb R,\qquad S^1\ltimes\Bbb R^2,\qquad S^1\ltimes(\text{cone
disk}),
$$
where $\ltimes$ denotes the metrically twisted product. Moreover, the 
$(1+\varepsilon)$-bilipschitz homeomorphism $f:U_x\to N_{\nu}(S)$ satisfies:
$$
\max(\Inj(x),d(f(x),S),\Diam(S))\leq\nu/D.
$$
  \enddefinition

Note that a $(\varepsilon,D)$-Margulis' neighborhood  of abelian type is
abelian in $C-\Sigma$. This definition is motivated by the following lemma, which is
the first step in the proof of Proposition 2.1.

\proclaim{Lemma 2.2}  Let $(C_n)_{n\in\N}$ be a sequence of hyperbolic cone
3-manifolds which collapses and
satisfies the hypothesis of Theorem A.
 If  both assertions 2) and 3) of Theorem A fail to hold, then, for every
$\varepsilon,D>0$, there exists $n_0$ such that, for $n\geq n_0$, every $x\in C_n$
has a  $(\varepsilon,D)$-Margulis' neighborhood  of abelian type. 
\endproclaim

\demo{Proof} Since the sequence collapses,  we can apply the Local Soul Theorem
(Chapter IV) and we show that the only possible local models are the three ones of
abelian type.

More precisely, since the supremum of the cone-injectivity radius 
converge to zero when $n$ goes to infinity,
 given $\varepsilon,D>0$ there exists $n_0$ such that
for $n\geq n_0$ the Local Soul Theorem applies to every point 
$x\in C_n$,
Since we assume that assertion 3) of Theorem A does not hold,
 by Corollary 4.1 of Chapter IV, the
compact models are excluded.
Hence we have to consider only the non-compact local models.  

From the hypothesis that assertion 2) of theorem A
does not hold, we get rid of the product model
$S^2(\alpha,\beta,\gamma)\times\Bbb R$, where $S^2(\alpha,\beta,\gamma)$
is a  Euclidean cone 2-manifold with underlying space the
sphere $S^2$ and singular set three points at which the sum of the three
cone angles is
$\alpha+\beta+\gamma=2\pi$.
 
Since the cone angles belong to $[\omega_0,\omega_1]$, with $0<\omega_0<\omega_1<\pi$, the local
models with a cone angle equal to $\pi$ cannot occur.

Finally, the last model to be eliminated is the one corresponding to the
normal bundle of the soul of a twisted fiber bundle over the Klein bottle $
K^2\tilde\times \Bbb R$. A neighborhood $U_x$ $(1+\varepsilon)$-bilips\-chitz
homeomorphic to this model does not intersect the singular set $\Sigma$. If this
local  model occurred, there would be a  Klein bottle $K^2\times\{0\}$ embedded in
$C-\Sigma$. This would contradict the fact that $C-\Sigma$ admits a complete
hyperbolic structure (Lemma 1.5). \qed
 \enddemo

Since $C-\Sigma$ admits a complete hyperbolic  structure (Lemma 1.5),
by Gromov \cite{\Gro} and Thurston \cite{\ThuNotes, Ch. 6},
$\norm{C-\Sigma}=\Vol(C-\Sigma)/v_3$,
where $v_3>0$ is a constant depending only on the dimension.  In particular,
$\norm{C-\Sigma}\neq 0$. Then the proof of Proposition 2.1 follows
from Lemma 2.2 and the next proposition:

\proclaim{Proposition 2.3}  There exists a universal constant $D_0>0$
such that, if $C$ is a closed hyperbolic cone 3-manifold  where every
point has an 
$(\varepsilon,D)$-Mar\-gu\-lis' neighborhood   of abelian type, with
$\varepsilon<1/2$ and $D>D_0$, then the simplicial volume
$\norm{C-\Sigma}$ is zero.
\endproclaim

We prove this proposition in Sections 3 and 4. In order to show that $\norm{C-\Sigma}$
vanishes, we adapt a construction of Gromov \cite{\Gro, Sec. 3.4} to the relative case.
This construction gives  a covering of $C$ by open sets that are
abelian in $C-\Sigma$, and the dimension of the covering is
$2$ in $C$ and $0$ in $\Sigma$. In fact, Proposition 2.3 can be seen as a version
of Gromov's Isolation Theorem \cite{\Gro, Sec. 3.4} for cone 3-manifolds.
\enddemo

\head 3. Coverings \`a la Gromov \endhead

\definition{Definition} For $\eta>0$, a covering $(V_i)_{i\in I}$ of a hyperbolic cone
3-manifold $C$ by open subsets is said to be a \it $\eta$-covering \`a la Gromov \rm
if it satisfies: \roster
  \item"1)" for every $i\in I$, there exists a metric ball $B(x_i,r_i)$ of radius 
          $r_i\leq 1$
          that contains $V_i$;
   \item"2)" if $B(x_i,r_i)\cap B(x_j,r_j)\neq\emptyset$, then
          $3/4\leq{r_i}/{r_j}\leq 4/3$;
   \item"3)" for $i\neq j$, $B(x_i,{r_i}/4)\cap B(x_j,{r_j}/4)=\emptyset$;
   \item"4)" every $x\in C$ belongs to an open set $V_i$ such that $d(x,\partial
V_i)
           \geq{r_i}/3$;
   \item"5)" for every $i\in I$, $\Vol(V_i)\leq\eta\, r_i^3$.
\endroster\enddefinition

\remark{Remark} Every $\eta$-covering  \`a la Gromov of a
closed hyperbolic cone 3-manifold is finite, because properties 2) and
3) forbide acumulating
sequences. \endremark

Our interest in $\eta$-coverings \`a la Gromov comes from the following proposition. Our proof of
this proposition follows closely Gromov's proof \cite{\Gro, Sec. 3.4} for the
Riemannian (non-singular) case.

\proclaim{Proposition 3.1} There exists a universal constant $\eta_0>0$ such that, for any
closed hyperbolic cone 3-manifold $C$ admitting a $\eta$-covering \`a la Gromov
$(V_i)_{i\in I}$
 with
$\eta\leq\eta_0$, there exists a continuous map  from $C$ to a simplicial
$2$-complex $f:C\to K^{(2)}$ satisfying:
\roster
\item"i)" for every vertex $v$ of $K^{(2)}$, there is $i(v)\in I$ such that 
$$
f^{-1}(\Star v) \subset \bigcup\limits_{V_j\cap V_{i(v)}\neq\emptyset} V_j
$$
\item"ii)" for every $x\in C$ that belongs to only one open set of the covering, $f(x)$
is a vertex of $K^2$.
 \endroster
 \endproclaim

\demo{Proof}
The proof consist of a sequence of lemmas, like in Gromov's proof  \cite{\Gro, Sec.
3.4}. We recall that  a covering has dimension $n$ if every point belongs to at most $n+1$ open sets
of the covering.

\proclaim{Lemma 3.2}  There is a universal integer $N>0$ such that, for
every closed hyperbolic cone 3-manifold $C$ and for every $\eta>0$, the dimension of any
$\eta$-covering \`a la Gromov of $C$  is at most $N$. 
 \endproclaim

\demo{Proof of Lemma 3.2} We shall bound the number $N_i$ of balls
$B(x_j,r_j)$ that intersect  a given ball $B(x_i,r_i)$. From property 
2) of the definition of a $\eta$-covering \`a la Gromov, if 
$B(x_j,r_j)\cap B(x_i,r_i)\neq\emptyset$, then  $3/4\leq
r_j/r_i\leq4/3$ and we have:
$$
\align
   &B(x_j,r_j)\subset B(x_i,r_i+2 r_j)\subset B(x_i, 4 r_i)\\
   \text{and }\quad &B(x_i, 4 r_i)\subset B(x_j, 5 r_i+r_j)\subset 
    B(x_j, 8 r_j).
\endalign
$$
By using these inclusions and the fact that the balls 
 $(B(x_j,r_j/4))_{j\in I}$ are pairwise disjoint, it follows that the
number $N_i$ of balls that intersect  a given $B(x_i,r_i)$ is bounded
above by:
$$
\align
  N_i&\,\leq\sup\left\{\frac{\Vol(B(x_i,4 r_i))}{\Vol(B(x_j, r_j/4))} 
   \mid 
   B(x_i,r_i)\cap B(x_j,r_j)\neq\emptyset \right\}
\\&
\leq\sup_{j\in I}\left\{\frac{\Vol(B(x_j,8 r_j))}{\Vol(B(x_j, r_j/4))}\right\}.
\endalign
$$
Now, the uniform upper bound for $N_i$ comes from the Bishop-Gromov 
Inequality (Proposition 1.6, Chapter III), which shows that:
$$
    \frac{\Vol(B(x_j,8 r_j))}{\Vol(B(x_j, r_j/4))}\leq 
       \frac{\V_{-1}( 8 r_j)}{\V_{-1}( r_j/4)},
$$
where $\V_{-1}(r)=\pi(\sinh(2r)-2r)$ is the volume of the ball  of
radius $r$ in the hyperbolic space $\Bbb H^3_{-1}$. Since the function
$
r\mapsto {\V_{-1}( 8 r)}/{\V_{-1}( r/4)}
$
is continuous, it is bounded on $[0,1]$. Hence, for any integer $N$  bounding
above this function on $[0,1]$, we get $N_i\leq N$ and the lemma is proved.
\qed
\enddemo

Given a $\eta$-covering \`a la Gromov,  its nerve $K$ is a simplicial
complex and, according to Lemma 3.2, the dimension of $K$ is at most
$N$, where $N$ is a uniform constant. Since we work with a
compact cone 3-manifold $C$, every
$\eta$-covering \`a la Gromov is finite and its nerve $K$ is compact.
We  canonically embeds $K$ in $\Bbb R^p$, where $p$ is the number of
vertices of $K$, which equals the number of open sets of this
covering. More precisely, every vertex of $K$ corresponds to a vector of the form
 $(0,\ldots,1,\ldots,0)$ and the simplices of positive dimension
 are defined by
linear extension.

 The proof of Proposition 3.1 goes as follows. We start 
in Lemma 3.3 by
constructing  a Lipschitz map from $C$ to the nerve of the covering
$f:C\to K$ which satisfies properties i) and ii) of 
Proposition 3.1. Next, in Lemma 3.4, we deform the map
$f$ to a Lipschitz map $f_3:C\to K^{(3)}$  where $K^{(3)}$ is the
$3$-skeleton of $K$. Finally, in Lemma 3.5, we prove that for
$\eta>0$  sufficiently small we can deform  $f_3:C\to K^{(3)}$ to
the $2$-skeleton $K^{(2)}$,  keeping properties i) and ii) of
Proposition 3.1. To prove the existence of such a universal constant
$\eta_0>0$ we need uniform constants in the Lemmas, the first example being the upper
bound $N$ of the dimension of $K$.

\proclaim{Lemma 3.3} Let $C$ be a hyperbolic cone $3$-manifold equipped
with a
$\eta$-covering \`a la Gromov $(V_i)_{1\leq i\leq p}$.  Let
$K=K^{(k)}$ be the nerve of
 this covering, which has dimension $k\leq N$. Then there exists a Lipschitz
map $f_k:C\to K$ that verifies properties i) and ii) of Proposition
3.1 and in addition:
\roster
\item"iii)" there exists a uniform constant $\xi_k$,  depending only
on the dimension
$k$, such that, for $1\leq i\leq p$, $$
\forall x,y\in\bigcup_{V_j\cap V_i\neq\emptyset} V_j, \qquad
\norm{f_k(x)-f_k(y)}\leq\frac{\xi_k}{r_i} d(x,y).
$$
\endroster\endproclaim

In this lemma, $d$  denotes the hyperbolic distance on $C$ and
$\norm{\,\,}$ the Euclidean norm on $\Bbb R^p$, since we assume that
$K$ is canonically embedded in $\Bbb R^p$.

\demo{Proof} We choose a  smooth  function $\phi:\Bbb R\to [0,1]$ such
that
$\phi((-\infty,0])=0$, $\phi([1/3,+\infty))=1$,  and
$\vert\phi'(t)\vert\leq 4$ for every $t\in\Bbb R$.

\midinsert
 \centerline{\BoxedEPSF{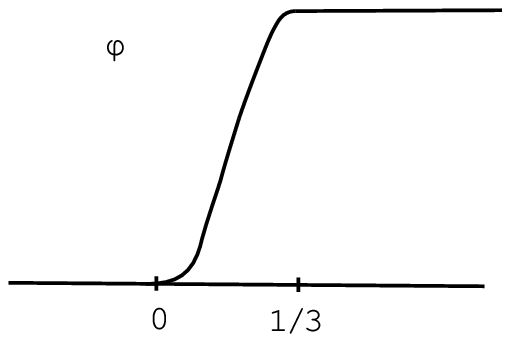 scaled 1000}}
    \botcaption{Figure V.1} The function $\phi$
    \endcaption
     \endinsert

For every $i=1,\ldots,p$, let $\phi_i:\overline{V_i}\to\Bbb R$ be 
the Lipschitz map such that
$
\forall x\in \overline{V_i}$  $\phi_i(x)=\phi({d(x,\partial
V_i)}/{r_i})$, where $\partial V_i$ is the boundary of $V_i$. 
Since $\phi_i$ vanishes on $\partial V_i$,  we extend it to the whole
manifold just by taking zero outside $V_i$. Then we have the property
$$
\forall x, y\in C, \qquad\norm{\phi_i(x)-\phi_i(y)}\leq\frac4{r_i} 
d(x,y) \tag*
$$
because the map
$$
x\mapsto \left\{\matrix \format \c\,&\l\\ d(x,\partial V_i)& \text { if } x\in 
\overline{V_i}\\
                           0 & \text{ otherwise }
              \endmatrix
\right.
$$
has Lipschitz constant $1$.

Let  $\Delta^{p-1}=\{(u_1,\ldots,u_p)\in\Bbb R^p\mid
u_1+\cdots+u_p=1\text{ and } u_i\geq 0\text{ for } i=1,\ldots,p\}
$
be the unit  simplex of $\Bbb R^p$. We define the map $f_k:C\to
\Delta^{p-1}$ to be:
$$
   \forall x\in C,\qquad
f_k(x)=\frac1{\sum\limits_{i=1}^p\phi_i(x)}(\phi_1(x),\cdots,\phi_p(x)).
 $$
This map is well defined, since $\sum\limits_{i=1}^p\phi_i(x)\geq 1$
by property 4) of a $\eta$-covering  \`a la Gromov.

The nerve $K$ of  the covering embeds canonically in $\Bbb R^p$ as a
subcomplex of
$\Delta^{p-1}$. Namely, if $V_1,\ldots,V_p$ are the open sets of the covering, then the
vertex of $K$ corresponding to $V_i$ is mapped to $i$-th vertex $(0,\ldots, 1,\ldots,0)$
of $\Delta^{p-1}$. 
By construction the image $f_k(C)$ is contained in
$K\subset\Delta^{p-1}$ and satisfies properties i) and ii) of
Proposition 3.1.

 We show now that $f_k$ satisfies property iii). 
We view $f_k$ as a  composition $f_k=g\circ (\phi_1,\cdots,\phi_p)$
where
$g(u_1,\ldots,u_p)=\frac1{u_1+\cdots+u_p}(u_1,\ldots,u_p) $. Since 
$\phi_1(x)+\cdots+\phi_p(x)\geq 1$ for every $x\in C$, and the map $g$
restricted to the region
$$
 \{(u_1,\ldots,u_p)\in\Bbb R^p\mid
u_1+\cdots+u_p\geq 1\text{ and } u_i\geq 0\text{ for } i=1,\ldots,p\}
$$
 has  Lipschtitz constant $\sqrt 2$, we have:
$$
\forall x,y \in C\qquad 
\norm{f_k(x)-f_k(y)}\leq \sqrt 2
 \norm{(\phi_1(x),\ldots,\phi_p(x))-(\phi_1(y),\ldots,\phi_p(y))}.
$$
From   inequality \ttag *, for every $i=1,\ldots,p,$ 
$$
\multline
\forall x,y\in\underset{V_j\cap V_i\neq\emptyset}\to\bigcup V_j, 
\qquad
\norm{(\phi_1(x),\ldots,\phi_p(x))-(\phi_1(y),\ldots,\phi_p(y))}^2
=\\
\sum\limits_{V_j\cap
V_i\neq\emptyset}(\phi_j(x)-\phi_j(y))^2
\leq 4^2 \sum\limits_{V_j\cap
V_i\neq\emptyset}\frac{1}{r_j^2}  d(x,y)^2
\endmultline
$$
From the property 2) of a $\eta$-covering \`a la Gromov:
$$
 \sum\limits_{V_j\cap
V_i\neq\emptyset}\frac1{r_j^2}\leq (k+1)\left(\frac{4
}{3 r_i}\right)^2,
$$
where $k$ is  the dimension of the covering.
Summarizing all the inequalities we get:
$$
\forall i=1,\ldots,p,\quad \forall x,y, \in\underset{V_j\cap
V_i\neq\emptyset}\to\bigcup V_j,  \qquad\quad
\norm{f_k(x)-f_k(y)}\leq \frac{16\sqrt{2(k+1)}\, d(x,y)}{3r_i}.
$$
Thus the lemma is proved by taking $ \xi_k=\frac{16 \sqrt{2(k+1)}}{3
}$.\qed
\enddemo

\proclaim{Lemma 3.4} With the hypothesis of Lemma 3.3, the Lipschitz
map $f_k:C\to K$ can be deformed to a Lipschitz map $f_3:C\to K^{(3)}$
into the $3$-skeleton which satisfies properties i), ii) (from
Proposition 3.2) and iii) (from Lemma 3.3).
\endproclaim

\demo{Proof} 
We start with the map $f_k:C\to K$ obtained in Lemma 3.3. If $k=\dim K=3$,
we are done. Hence we assume that $k>3$ and we prove the lemma by
induction: we show that whenever we have a map $f_k:C\to K^{(k)}$
satisfying properties i), ii) and iii) with $k>3$, then we can deform
it to a map into the $(k-1)$-skeleton $ K^{(k-1)}$ satisfying the same
properties. The key point in the argument is the following technical
claim.
\enddemo

\proclaim{Claim 3.5} Given  a Lipschitz map $f_k:C\to K^{(k)}$ satisfying properties
i), ii) and iii), there is a  uniform constant $\varepsilon_k>0$
(which depends only on $k$) such that
 every $k$-simplex $\Delta^k\subset K$ contains a point $z$ 
at distance at least $\varepsilon_k$ from the image $f_k(C)$ and the boundary
$\partial\Delta^k$.
\endproclaim

\demo{Proof of the claim}  Let $\varepsilon>0$ be such that every
point in
$\Delta^k$ is at distance at most  $\varepsilon>0$ from the union 
$f_k(C)\cup\partial\Delta^k$. We are going to find a uniform constant $\varepsilon_k>0$ 
such that $\varepsilon\geq\varepsilon_k$.

Let $\{z_1,\ldots,z_s\}\subset\Delta^k$ be a maximal  family of points
such that 
$d(z_i,\partial\Delta^k)\geq 3\varepsilon$ and  $\norm{z_i-z_j}\geq
3\varepsilon$ for 
$i\neq j$.  There exists a constant $c_1=c_1(k)>0$ depending only on
the dimension $k$ such that, for $\varepsilon$ sufficiently small, we
can find at least
$c_1\varepsilon^{-k}$ points in this family. So we can assume 
$s\geq c_1\varepsilon ^{-k}$.

By the hypothesis on  $\varepsilon>0$, we can find a family of points
 $\{y_1,\ldots,y_s\}\subset
f_k(C)\cap\Delta^k$ such that
$\norm{z_i-y_i}\leq\varepsilon$, for $i=1,\ldots,s$.  In particular,
$\norm{y_i-y_j}>\varepsilon$ (if $i\neq j$) and $d(y_i,\partial
\Delta^k)>\varepsilon$. Choose points $\{\bar y_1,\ldots,\bar
y_s\}\subset C$ such that $f_k(\bar y_i)=y_i$,
$i=1,\ldots,s$. From property i) of $f_k$, the points 
$\{\bar y_1,\ldots,\bar y_s\}$ belong to $\underset{V_j\cap V_{i(v)}\neq 0}\to \bigcup
V_j$, where the open set $V_{i(v)}$ corresponds to a vertex $v$ of $\Delta^k$.
So, from property iii) we have 
$$
  \forall i\neq j\in\{1,\cdots,s\}, \qquad d(\bar y_i,\bar y_j)\geq
 \frac{r_{i(v)}}{\xi_k}\norm{y_i-y_j} >
\frac{r_{i(v)}}{\xi_k}\varepsilon.
$$
This implies that the balls  $B(\bar y_j,\frac{r_{i(v)}}{2 \xi_k}\varepsilon)$ are
pairwise disjoint and verify 
$$
B(\bar y_j,\frac{r_{i(v)}}{2 \xi_k}\varepsilon)\subset\underset{V_j\cap V_{i(v)}\neq
0}\to \bigcup V_j \subset B(x_{i(v)},4 r_{i(v)}),
$$
where the last inclusion follows from property 2) of a $\eta$-covering \`a la
Gromov. We get the following upper bound for the number $s$ of such balls:
$$
s\leq\max_{j=1,\ldots,s}\Big( \frac{\Vol(B(x_{i(v)},4r_{i(v)}))}
{\Vol(B(\bar y_j,\frac{r_{i(v)}}{2\xi_k}\varepsilon))} \Big)
\leq\max_{j=1,\ldots,s}\Big( \frac{\Vol(B(\bar y_j,8r_{i(v)}))}
{\Vol(B(\bar y_j,\frac{r_{i(v)}}{2\xi_k}\varepsilon))} \Big),
$$
because $B(x_{i(v)},4r_{i(v)}) \subset B(\bar y_j,8r_{i(v)})$, for $j=1,\ldots,s$.
From the Bishop-Gromov inequality (Proposition 1.6, Chapter III) we obtain:
$$
s\leq \frac{\V_{-1}(8r_{i(v)})}{\V_{-1}(\frac{r_{i(v)}}{2\xi_k}\varepsilon)},
$$
where $\V_{-1}(r)=\pi(\sinh(2r)-2r)$ is the volume of the ball of radius $r$ in
the hyperbolic $3$-space. There exists a constant $a>0$ such that 
$\frac{r^3}a\leq\V_{-1}(r)\leq a r^3 $
$
\forall r\in [0,8]$.
Since $r_{i(v)}\leq 1$, we obtain the upper bound:
$$
s\leq \frac{a^2(8 r_{i(v)})^3 8 \xi_k^3}{r_{i(v)}^3\varepsilon^3}=c_2\varepsilon^{-3},
$$
where $c_2=2^{12} a^2 \xi_k^3>0$ depends only on the dimension $k$.
By combining both inequalities
$c_1\varepsilon^{-k}\leq{s}\leq{c_2}\varepsilon^{-3}$,   we conclude
that
$
\varepsilon\geq (c_1/c_2)^{\frac1{k-3}}$, with $k>3$.
This finishes the proof of  Claim 3.5.\qed
 \enddemo

\demo{End of the proof of Lemma 3.4}  We asume $k>3$ and we want to
construct
$f_{k-1}:C\to K^{(k-1)}$.  Let $\Delta^k_1,\ldots,\Delta^k_q$ be the
$k$-simplices of
$K$. From Claim 3.5, for every $k$-simplex $\Delta^k_i\subset K$ we
can choose a point
$z_i\in\Delta^k_i$ so that
$d(z_i,f_k(C)\cup\partial\Delta^k_i)>\varepsilon_k$.  We consider the
map
$R_i:K-\{z_i\}\to K$ which  is defined by the radial retraction of
$\Delta^k_i-\{z_i\}$ onto $\partial\Delta^k_i$ and the identity on $K-\Delta^k_i$. 
Since the points $\{ z_1,\ldots,z_q\}$ do not belong to the image of 
$f_k$, the composition 
$$
  f_{k-1}=R_1\circ \cdots \circ R_q\circ f_k: C\to K
$$
is well defined, and the image $f_{k-1}(C)$ lies in the $(k-1)$-skeleton $K^{(k-1)}$.
Moreover, it follows from the construction that $f_{k-1}$ satisfies properties i) and
ii) of Proposition 3.1, because the retractions $R_i$ preserve the vertices and their stars.

For $i=1,\ldots,q$, the retraction $R_i:K-\{z_i\}\to K$ is piecewise smooth. From the 
inequality $d(z_i,f_k(C)\cup\partial\Delta^k_i)>\varepsilon_k$, it follows that the
local Lipschitz constant of $R_1\circ \cdots \circ R_q$ is uniformly bounded on the
image
$f_k(C)$; moreover the bound  depends only on the dimension $k$, because the constant
$\varepsilon_k$ is uniform, depending only on the dimension $k$. Thus $f_{k-1}$ satisfies also
property iii) of Lemma 3.3. \qed
 \enddemo

Next lemma completes the proof of Proposition 3.1.

\proclaim{Lemma 3.6} There exists a universal constant $\eta_0>0$ such
that, for every
 $\eta$-cov\-er\-ing \`a la Gromov of $C$, the map $f_3:C\to K^{(3)}$ of
Lemma 3.4 can be deformed to a continuous map $f_2:C\to K^{(2)}$  into
the
$2$-skeleton which satisfies properties i) and ii) of Proposition 3.1.
\endproclaim

\demo{Proof of Lemma 3.6} 
To deform $f_3$ to $f_2$,
it suffices to prove that in every  $3$-simplex $\Delta^3\subset K$,
there is a point
$z\in\Int(\Delta^3)$ that does not belong to the image $f_3(C)$. 
Then such a deformation is constructed by composing  $f_3$ with all the
radial retractions  from $\Delta^3-\{z\}$ to $\partial\Delta^3$ as in
Lemma 3.4. The  map $f_2$ will satisfy properties i) and ii) of
Proposition 3.1 by construction. Next claim shows that 
$\Int(\Delta^3)-f_3(C)$ is non-empty whenever $\eta$ is less than a
universal constant
$\eta_0>0$. This will conclude the proof of Lemma 3.6.

\proclaim{Claim 3.7}
There exists a universal constant $\eta_0>0$ such that, if $C$ admits
a $\eta$-covering \`a la Gromov
with $\eta<\eta_0$, then for every $3$-simplex $\Delta^3\subset K^{(3)}$
$$
\Vol(\Delta^3\cap f_3(C)) < \Vol(\Delta^3).
$$
 \endproclaim

\demo{Proof of the claim} Property ii) of the map $f_3:C\to K^{(3)}$ implies the
 following inequality for every $3$-simplex $\Delta^3\subset K^{(3)}$:
$$
\Vol(\Delta^3\cap f_3(C))\leq\sum_{V_j\cap V_{i(v)}\neq\emptyset}\Vol(f_3(V_j)),
$$
where $V_{i(v)}$ is the open set corresponding to a vertex $v$ of $\Delta^3$. The map
$f_3$ is Lipschitz, and from property iii), its restriction to 
$\underset{V_j\cap V_{i(v)}\neq\emptyset}\to\bigcup V_j$ has Lipschitz constant
$\xi_3/r_{i(v)}$. Hence, according to the formula giving a bound for the volume of the image
of a Lipschitz map (see \cite{\Fed, Corollary 2.10.11}), we get:
$$
\sum_{V_j\cap V_{i(v)}\neq\emptyset}\Vol(f_3(V_j))\leq  
\sum_{V_j\cap V_{i(v)}\neq\emptyset}\Vol(V_j) \left(\frac{\xi_3}{r_{i(v)}}\right)^3.
$$
Property 5) of a $\eta$-covering \`a la  Gromov asserts that $\Vol(V_j)\leq \eta\, r_j^3$.
Furthermore, from property 2) of these coverings, we have  $ r_j\leq \frac43r_{i(v)}$
whenever 
$V_j\cap V_{i(v)}\neq\emptyset$. Thus we deduce the following inequalities:
$$
\Vol(\Delta^3\cap f_3(C))\leq\sum_{V_j\cap V_{i(v)}\neq\emptyset}\Vol(f_3(V_j))
\leq \eta \,\left(\frac43\xi_3\right)^3 (N+1),
$$ 
where $N$ is the universal upper bound of the dimension of the covering given by Lemma 3.2.
Hence it suffices to take $\eta_0 < \Vol(\Delta^3)/((N+1)(\frac43\xi_3)^3)$  to prove the
claim. \qed
\enddemo
\enddemo
\enddemo

\head  
      4. From  $(\varepsilon,D)$-Mar\-gu\-lis'
coverings  of abelian type to
$\eta$-coverings \`a la Gromov 
\endhead

The aim of this section is to prove Proposition 2.3. 
We recall the statement:

\proclaim{Proposition 2.3}  There exists a universal
constant $D_0>0$ such that, if $C$ is a closed
hyperbolic cone 3-manifold  where every point has an 
$(\varepsilon,D)$-Mar\-gu\-lis' neighborhood  of
abelian type with $\varepsilon<1/2$ and $D>D_0$, then
the simplicial volume $\norm{C-\Sigma}$ is zero.
\endproclaim

The proof follows from Proposition 3.1 and the following:

\proclaim{Proposition 4.1}  There is a universal constant $b_0>0$
such that,  if $C$ is a closed hyperbolic cone 3-manifold 
where each point $x\in C$ has an
$(\varepsilon,D)$-Margulis' neighborhood of abelian type with
$\varepsilon\leq\frac12$ and $D\geq300$, then  $C$ admits a
$\eta$-covering \`a la Gromov with $\eta<b_0/D$.

 Moreover, the open sets
$(V_i)_{i\in I}$ of the $\eta$-covering
\`a la Gromov satisfy the following additional properties:
\roster
\item"6)" $\forall i\in I$, $\underset{V_j\cap V_{i}\neq\emptyset}\to\bigcup V_j$ is
abelian in $C-\Sigma$.
\item"7)" there is a tubular neighborhood $\NN (\Sigma)$ of  $\Sigma$
such that every component of  $\NN (\Sigma)$ is contained in only one open set of the
covering.
\endroster
\endproclaim

\demo{Proof of Proposition 2.3}
 We choose $D_0=\max(b_0/\eta_0,300)$, where $\eta_0>0$ is the universal constant
of Proposition 3.1. From Propositions 3.1 and 4.1, since every point of
$C$ has an $(\varepsilon, D)$-neighborhood of abelian type, we can
construct a continuous map from $C$ to a
$2$-dimensional simplicial complex: $f:C\to K^2$. Moreover properties i) and ii) of
Proposition  3.1 together with properties 6) and 7) of Proposition 4.1 imply that $f$
satisfies:
\roster
\item"i$'$)" for every vertex $v$ of $K^2$,
$f^{-1}(\Star v)$ is abelian in
$C-\Sigma$;
\item"ii$'$)" there is an open tubular neighborhood $\NN(\Sigma)$ of $\Sigma$ such that
$f(\overline{\NN(\Sigma_i)})$ is a vertex of $K^2$ for every component
$\Sigma_i$ of $\Sigma$.
 \endroster

 Let $C(\lambda_1,\ldots,\lambda_q)$ denote the closed 3-manifold obtained by gluing $q$
solid tori to the boundary of the manifold $C-\NN(\Sigma)$, so that the
boundaries of the meridian disks are identified respectively to the simple closed
curves $\lambda_1,\ldots,\lambda_q$ in $\partial \overline{\NN(\Sigma)}$. More precisely,
$$
C(\lambda_1,\ldots,\lambda_q)=(C-\NN(\Sigma))
\bigcup_{\phi_{1},\ldots,\phi_q}\,\bigsqcup_{i=1}^qS^1\times D^2_i,
$$
where the gluing maps $\phi_i:\partial\NN(\Sigma_i)\to S^1\times\partial D^2_i$
satisfy $\phi_i(\lambda_i)=(\{ *\}\times\partial D^2_i)$, for $i=1,\ldots,q$.

From properties i$'$) and ii$'$), the continuous map $f:C\to K^2$
induces a map
$\bar f: C(\lambda_1,\ldots,\lambda_q)\to K^2$. Since abelianity is preserved by quotient,
$\bar f^{-1}(\Star v)$ is abelian in 
$C(\lambda_1,\ldots,\lambda_q)$, for every vertex $v$ of $K^2$.

It follows that the closed orientable 3-manifold 
$C(\lambda_1,\ldots,\lambda_q)$ admits an abel\-ian covering of
dimension 2. The vanishing theorem for simplicial volume (cf.
\cite{\Gro, Sec. 3.1} and \cite{\Iva})  shows that
$\norm{C(\lambda_1,\ldots,\lambda_q)}=0$. This holds for every choice
of simple closed curves on $\partial \NN(\Sigma)$. Thus, from
Thurston's hyperbolic Dehn filling theorem \cite{\ThuNotes} (cf. Theorem 1.1, Chapter II):
$$
\norm{C-\Sigma}=\lim_{\operatorname{Lenght}(\lambda_i)\to\infty}
\norm{C(\lambda_1,\ldots,\lambda_q)}=0. \qed
$$
\enddemo

The remaining of this section is devoted to the proof of   Proposition
4.1.

\demo{Proof of Proposition 4.1} Let $C$ be a closed  hyperbolic cone
3-manifold so that every point $x\in C$ admits un
$(\varepsilon,D)$-Margulis' neighborhood of abelian type, with
$\varepsilon<\frac12$ and $D>300$. It means that $x$
has a neihgbourhood $U_x\subset C$ that is bilipschitz homeomorphic
to the normal cone fiber bundle  $\NN_{\nu}(S)$, of radius $\nu<1$,  of the soul
$S$ of one of the following non-compact Euclidean cone 3-manifolds: $T^2\times\Bbb
R$, $S^1\ltimes\Bbb R^2$, $S^1\ltimes( \text{cone disk})$. Moreover,
the $(1+\varepsilon)$-bilipschitz homeomorphism $f: U_x\to\NN_{\nu}(S)$
satisfies:
\roster
\item"(a)" $\max\big(\Inj(x),d(f(x),S),\Diam(S)\big)\leq \frac{\nu}D$,
\endroster
(cf. Local Soul Theorem, Chapter IV, and Lemma 2.2, Chapter V).

For every point $x\in C$, we define the abelianity radius $\ab x$ to be:
$$
    \ab x =\sup\{ r>0\mid B(x,r) \text{ is abelian in } C-\Sigma \}.
$$
 By using the $(1+\varepsilon)$-bilipschitz homeomorphism $f: U_x\to\NN_{\nu}(S)$ and the
majoration (a), we get:
$$
  \ab x\geq\frac{\nu}{1+\varepsilon}(1-\frac1D)\geq\frac{\nu}2\geq\frac{D}2\Inj(x).
$$

For every $x\in C$ we define $r(x)=\inf(\frac{\ab x}8,1)$. Lemma 4.2 and 4.3 give the
first properties of  the balls $B(x,r(x))$.

\proclaim{Lemma 4.2}  Let $x,y\in C$. If $B(x,r(x))\cap
B(y,r(y))\neq\emptyset$, then \roster
\item"(b)" $3/4\leq{r(x)}/{r(y)}\leq 4/3$;
\item"(c)" $B(x,r(x))\subset B(y,4r(y))$.
\endroster
\endproclaim

\demo{Proof} Assume $r(x)\geq r(y)$. Either $r(y)=1$ or $r(y)=\ab
y/8$. If $r(y)=1$, then $r(x)=1$ and assertion (b) is clear. If $r(y)=\ab
y/8$, by using the inclusion
$
B(y,6 r(x))\subset
B(x,8 r(x))
$ 
and the fact that $8r(x)\leq\ab x$, it follows that
 $B(y,6r(x))$ is abelian in
$C-\Sigma$.
Hence
 $r(x)\leq \ab y/6\leq 4 r(y)/3$ and
(b) is proved. 

Assertion (c) follows easily from (b) and the
inclusion $B(x,r(x))\subset B(y,2 r(x)+r(y))$.
\qed
\enddemo

\proclaim{Lemma 4.3} 
Let $\Sigma=\Sigma_1\sqcup\cdots\sqcup\Sigma_q$ be the
singular  set of $C$. We choose   a point $x_i$ in each
connected component $\Sigma_i$.  Then, we have
the following properties: \roster
\item"(d)" for $i=1,\ldots,q$, if $\mu>0$ is 
sufficiently small, then $\NN_{\mu}(\Sigma_i)\subset
B(x_i,\frac{r(x_i)}4)$, where  $\NN_{\mu}(\Sigma_i)$ is
the tubular neighborhood of radius  $\mu$ around the
connected component $\Sigma_i$; \item"(e)"
$B(x_i,r(x_i))\cap B(x_j,r(x_j))=\emptyset$, for $i\neq
j$, $i,j\in\{1,\ldots,q\}$. \endroster
\endproclaim

\demo{Proof} Property (d) follows from  the hypothesis
that $x_i$ has an $(\varepsilon, D)$-Margulis'
neigbourhood of abelian type. Since $x_i$ is singular,
the local model is the normal cone fiber bundle
$\NN_{\nu}(S)$, of radius $\nu<1$, of the soul  $S=S^1\times\{\text{cone
point}\}$ of the Euclidean cone 3-manifold
$S^1\ltimes(\text{cone disk})$. 

Let $f:U_{x_i}\to \NN_{\nu}(S)$ be the  $(1+\varepsilon)$-bilipschitz
homeomorphism between
$U_{x_i}$ and the local model, then
 $U_{x_i}\cap\Sigma=\Sigma_i=f^{-1}(S)=f^{-1}(S^1\times\{\text{cone
point}\})
$.
Since $\varepsilon\leq 1/2$, it follows from the upper bound (a) that
$\Diam(\Sigma_i)\leq\Diam(S)(1+\varepsilon)\leq 2\frac{\nu}D$. Furthermore, since
$\nu\leq\inf(1,2\ab {x_i})$ ,
$r(x_i)=\inf(1,\frac{\ab{x_i}}8)$ and $D>300$ we get:
$$
\Diam(\Sigma_i)\leq\inf(\frac2D,\frac{4\ab{x_i}}D)< r(x_i)/9,
$$
         
Hence $\Sigma_i\subset B(x_i,\frac{r(x_i)}9)$.
By taking $\mu\leq\inf\{\frac{r(x_i)}{18}\mid i=1,\ldots,q\}$  we
obtain the inclusion
$\NN_{\mu}(\Sigma_i)\subset B(x_i,\frac{r(x_i)}4)$.

To show property (e),  we assume that there are $i\neq
j$ such that $B(x_i,r(x_i))\cap
B(x_j,r(x_j))\neq\emptyset$ and we seek a
contradiction. From property (c) of Lemma 4.2, the fact
that the balls intersect implies that
$B(x_j,r(x_j))\subset  B(x_i,4r(x_i))$. Hence, by
property (d), $\Sigma_i\cup\Sigma_j\subset B(x_i,4
r(x_i))$,  which is an abelian ball in $C-\Sigma$. This
implies that the two peripheral elements of
$\pi_1(C-\Sigma)$ represented by the meridians of
$\Sigma_i$ and $\Sigma_j$ commute. This contradicts
the fact that $C-\Sigma$ admits a complete hyperbolic structure  (Lemma 1.5). \qed
 \enddemo

We can now start the construction of   the $\eta$-covering
\`a la Gromov. First, we choose a point $x_i$ on each
connected component $\Sigma_i$ of $\Sigma$,
$i\in\{1,\ldots,q\}$. We fix the points
$\{x_1,\ldots,x_q\}$ and we consider all the possible finite
sequences of points  $\{x_1,\ldots,x_q,x_{q+1},\ldots,x_p\}$, 
starting with these fixed $q$ points and having the
property that $$
 \text{the balls } B\big(x_n,\tfrac{r(x_n)}4\big) \text{
are pairwise disjoint.}\tag* $$
Note that a sequence satisfying \ttag*  and Lemma 4.2 is necessarily
finite because $C$ is compact. The following lemma is
due to Gromov \cite{\Gro, Sec. 3.4 Lemma B}:

\proclaim{Lemma 4.4}
Let $x_1,\ldots,x_p$ be a finite sequence in $C$ as above, with the first fixed $q$ points in
the singular set. If it is maximal for property \ttag*,   then the balls
$B\big(x_1,\frac23r(x_1)\big),\ldots,B\big(x_p,\frac23r(x_p)\big)$ cover $C$.
\endproclaim

\demo{Proof} Let $x\in C$. By maximality, the ball $B\big(x,\frac{r(x)}4\big)$ intersects
$B\big(x_i,\frac{r(x_i)}4\big)$ for some $i\in\{1,\ldots,p\}$. From property (b) of Lemma 4.2,
$r(x)\leq\frac43r(x_i)$ and thus $x\in B\big(x_i,\frac{r(x_i)+r(x)}4\big)\subset
B\big(x_i,\frac23 r(x_i)\big)$.\qed
 \enddemo

Let $0<\mu\leq\inf\{\frac{r(x_i)}{18}\mid i=1,\ldots,q\}$  so that
$\NN_{\mu}(\Sigma_i)\subset  B\big(x_i,\frac{r(x_i)}4\big)$, as in Lemma
4.3 (d). Let $x_1,\ldots,x_p$ be a sequence as in Lemma 4.4, we consider the covering
$(V_i)_{i\in\{1,\ldots,p\}}$ defined by:
$$
\left\{\matrix\format\l&\qquad\l\\
  V_i= B(x_i,r(x_i)) &\text{ for } i=1,\ldots,q;  \\
   V_i= B(x_i,r(x_i))-\NN_{\mu}(\Sigma) &\text{ for }i=q+1,\ldots,p. 
\endmatrix\right. \tag{**}
$$

The following Lemma finishes the proof of Proposition 4.1.

\proclaim{Lemma 4.5} There is a universal constant $b_0>0$ such that the
above covering 
$(V_i)_{i\in\{1,\ldots,p\}}$ defined by \ttag{**} is a
$\eta$-covering \`a la Gromov with $\eta<\frac{b_0}D$ and satisfies
Properties {\rm 6)} and {\rm 7)} of Proposition {\rm 4.1}.
 \endproclaim

\demo{Proof of Lemma 4.5} We start by checking  that the covering 
satisfies properties 1) to 5) of a $\eta$-covering \`a la Gromov.
Property 1) follows from the construction by setting $r_i=r(x_i)$, for
$i=1,\ldots,p$. Property 2) follows immediately from Lemma 4.2, and
property 3) is the hypothesis
\ttag*.

\proclaim{Claim 4.6} The covering $(V_i)_{i\in\{1,\ldots,p\}}$
satisfies property 4) of a $\eta$-covering
\`a la Gromov. That is, $\forall x\in C$ there is an open set $V_i$ such
that $x\in V_i$ and $d(x,\partial V_i)>r_i/3$.
 \endproclaim

\demo{Proof} Let $x\in C$. From Lemma 4.4, $x\in B(x_i,\frac23 r_i)$
for some
$i=1,\ldots,p$. If $i\in\{1,\ldots,q\}$ (i.e. if $x_i\in\Sigma$) or if
$\NN_{\mu}(\Sigma)\cap B(x_i,r_i)=\emptyset$, then by construction
\ttag{**} 
$V_i=B(x_i,r_i)$ and we have $d(x,\partial V_i)\geq\frac{r_i}3$. 

Thus we may assume
that $\NN_{\mu}(\Sigma)\cap B(x_i,r_i)\neq\emptyset$. 
Let $j\in\{1,\ldots,q\}$ be an
index so that $\NN_{\mu}(\Sigma_j)\cap B(x_i,r_i)\neq\emptyset$; we can also assume that
$d(x,x_j)\geq\frac23r_j$.
By construction $V_i=B(x_i,r_i)-\NN_{\mu}(\Sigma)$; hence it is enough to show that the
distance of $x$ to the component $\NN_{\mu}(\Sigma_j)$ is at least $\frac13r_i$
whenever $\NN_{\mu}(\Sigma_j)\cap B(x_i,r_i)\neq\emptyset$.

Since  $d(x,\Sigma_j)\geq d(x,x_j)-\Diam(\Sigma_j)$, $ d(x,x_j)\geq\frac23r_j$   and
$\Diam(\Sigma_j)\leq\frac{r_j}9$  (by the proof of Lemma 4.3), we get
$d(x,\Sigma_j)\geq\frac59r_j$. Moreover, from property (b) of Lemma 4.2,
$ \frac43r_j\geq r_i$, hence:
$$
d(x,\Sigma_j)\geq(\tfrac49+\tfrac19)r_j\geq\tfrac13r_i+\tfrac19r_j.
$$
By the choice of $\mu\leq\inf\{\frac1{18} r_j\mid j=1,\ldots,q\}$ we can conclude that 
$
d(x,\NN_{\mu}(\Sigma_j))>\frac13r_i.
$
Hence $d(x,\partial V_i)>\frac13r_i$ and the claim is proved.
\qed
\enddemo

Property 5) of a $\eta$-covering \`a la Gromov given by the following claim:

\proclaim{Claim 4.7}  There is a universal constant
$b_0>0$ such that
$$
\Vol(V_i)\leq\Vol(B(x_i,r_i))
\leq\frac{b_0}Dr^3_i,\qquad\text{ for }i=1,\ldots,p. 
$$
\endproclaim

\demo{Proof} For $i=1,\ldots,p$, $x_i$ has  an
$(\varepsilon,D)$-Margulis' neighborhood of abelian type
$U_x$ which is $(1+\varepsilon)$-bilipschitz homeomorphic
to the normal cone fiber bundle  $\NN_{\nu}(S)$,  of radius $\nu<1$,  of the soul
$S$ of one of the following non-compact  Euclidean cone 3-manifolds: $T^2\times\Bbb
R$, $S^1\ltimes\Bbb R^2$, $S^1\ltimes( \text{cone disk})$.
The $(1+\varepsilon)$-bilipschitz homeomorphism $f:U_{x_i}\to
\NN_{\nu}(S)$ satisfies  $\varepsilon<1/2$ and
$\max(\Inj(x_i),d(f(x_i),S),\Diam(S))\leq\nu/D$. From
these inequalities we deduce that $\ab{x_i}\geq\nu/2$; hence
$r_i\geq\nu/16$.

 From the Bishop-Gromov inequality (Proposition 1.6
Chapter III) we get:  $$ \Vol(B(x_i,r_i))\leq
\Vol\big(B(x_i,\tfrac{\nu}{16})\big)\frac{\V_{-1}(r_i)}{\V_{-1}(\frac{\nu}{16})}. 
$$ Let $a>0$ be a constant so that
${t^3}/a\leq{\V_{-1}(t)}\leq a{t^3}$ for every $t\in[0,1]$. Since
$\nu\leq 1$ and $r_i\leq 1$, we get:
$$
\Vol(B(x_i,r_i))\leq
\Vol\big(B(x_i,\tfrac{\nu}{16})\big)2^{12}a^2\frac{r_i^3}{\nu^3}.
$$
Since $d(f(x_i),S)\leq\frac{\nu}{D}\leq\frac{\nu}{300}$, we have the inclusion
$f(B(x_i,\frac{\nu}{16}))\subset\NN_{\nu}(S)$. Thus:
$$
\Vol\big(B(x_i,\tfrac{\nu}{16})\big)\leq (1+\varepsilon)^3\Vol(\NN_{\nu}(S))
\leq 2^3\Vol(\NN_{\nu}(S)),
$$
because $f$ is $(1+\varepsilon)$-bilipschitz with $\varepsilon<1/2$.

By using the upper bound $\Diam(S)\leq\nu/D$ and the fact that S is of dimension
$1$ or $2$, a simple computation of Euclidean volumes gives the
upper bound:
 $$
\Vol(\NN_{\nu}(S))\leq\frac{\pi}D\nu^3.
$$
Thus:
$$
\Vol(V_i)\leq\Vol(B(x_i,r_i))\leq\frac{b_0}Dr^3_i,  \qquad\text{ where }
b_0=2^{15}\pi a^2.\qed
 $$
 \enddemo

Property 6) of Proposition 4.1 follows  from property
c) of Lemma 4.2, because $$
\forall i=1,\ldots,q,  \qquad \bigcup_{V_j\cap
V_i\neq\emptyset} V_j
\subset
\bigcup_{V_j\cap 
V_i\neq\emptyset}
B(x_j,r_j)
     \subset  B(x_i,4r_i)
$$
and by construction 
the ball $B(x_i,4 r_i)$ is abelian in
$C-\Sigma$.

Finally, property 7) of  
Proposition 4.1 follows immediatly from  property d) of
Lemma 4.3 and the construction of the covering. 

This finishes the proof of Proposition 4.1,  and thus
of Theorem A. \qed
\enddemo
 \enddemo

\leftheadtext{ }
\newpage
\rightheadtext{VI \qquad    Orbifolds and sequences of cone manifolds}
\leftheadtext{VI \qquad    Orbifolds and sequences of cone manifolds}

\

\centerline{\smc chapter  \ vi} 

\

\centerline{\chapt ORBIFOLDS  \, AND \,  SEQUENCES \, OF}

\vglue.2cm

\centerline{\chapt  HYPERBOLIC  \, CONE \,  3-MANIFOLDS}

\

\

This chapter is devoted to the proof of  Theorem B.

\proclaim{Theorem B} Let $\OO$ be a closed orientable connected 
irreducible very good
$3$-orbifold with topological type $(\Mod\OO,\Sigma)$ and
ramification indices $n_1,\ldots,n_k$. Assume that there exists a
sequence of hyperbolic cone 3-manifolds 
$(C_n)_{n\in\N}$  with the same topological type 
$(\Mod\OO,\Sigma)$ and
such that, for each component of $\Sigma$,  the cone angles form an
increasing sequence that converges   to $2\pi/n_i$ when $n$ approaches
$\infty$.

Then $\OO$ contains a non-empty compact
essential $3$-suborbifold
$\OO'\subseteq\OO$, which is not a product and which is either complete 
hyperbolic of finite volume, Euclidean,  Seifert fibered or $Sol$.
\endproclaim

We recall that a  compact
orientable $3$-suborbifold $\OO'$ is essential in a $3$-orbifold $\OO$
if the 2-suborbifold $\partial\OO'$
is either empty or incompressible in $\OO$.

The suborbifold  $\OO'$ of the theorem is not necessarily proper, it can be
$\OO'=\OO$, but it is non-empty. By saying that $\OO'$ is complete hyperbolic of
finite volume we mean that its interior has a complete hyperbolic structure
of finite volume. In particular, $\partial\OO'$ is a
collection of Euclidean $2$-suborbifolds.

The proof of Theorem B splits into two cases, according to
whether  the  sequence of cone 3-manifolds
$(C_n)_{n\in\N}$ collapses or not, as in Theorem A.

\head 1. The non-collapsing case
\endhead

Next proposition proves Theorem B when the  sequence of cone
3-manifolds $(C_n)_{n\in\N}$ does not collapse (i.e. $\sup\{\Inj(x)\mid x\in
C_n\}$ does not converge to zero).

\proclaim{Proposition 1.1}
Let $\OO$ and $(C_n)_{n\in\N}$  satisfy  the hypothesis of Theorem B.
 If the sequence $(C_n)_{n\in\N}$ does not collapse, then $\OO$ contains a non-empty compact
essential $3$-suborbifold that is complete hyperbolic of finite volume.
\endproclaim

\demo{Proof} Since the sequence $(C_n)_{n\in\N}$ does not collapse,
 after passing to a subsequence if
necessary, there  is a constant
$a>0$ and, for every $n\in\N$, 
there is a
point $x_n\in C_n$ such that $\Inj(x_n)\geq
a$. Thus, the sequence of pointed cone 3-manifolds $(C_n,x_n)$ is contained in
$\CC_{[\omega_0,\pi],a}$, for some $\omega_0>0$, because the cone angles of
$C_n$ converge to angles of the form $2\pi/n_i$. 

Since $(C_n,x_n)\in\CC_{[\omega_0,\pi],a}$, by the Compactness
Theorem (Chapter III), after passing to a subsequence, we can assume that $(C_n,x_n)_{n\in\N}$
converges geometrically to a pointed hyperbolic cone 3-manifold $(C_{\infty},
x_{\infty})$. By hypothesis, the cone angles of $C_{\infty}$ are of the
form $2\pi/m$ with $m\in\N$, hence $C_{\infty}$ is an orbifold.
We distinguish two cases, according to whether the limit 3-orbifold $C_{\infty}$ is
compact or not.

If the limit 3-orbifold $C_{\infty}$ is compact, then the geometric
convergence implies that $C_{\infty}$ has the same topological type
$(C,\Sigma)$ as the orbifold $\OO$. Moreover, the branching indices of
$C_{\infty}$ agree with the ones of $\OO$. Therefore
as an orbifold $C_{\infty}=\OO$ and $\OO$ is a closed hyperbolic orbifold. Thus
Proposition 1.1 is proved in this case.

If the limit $C_{\infty}$ is not compact, then we need further work, as in
Chapter V. The first step is the following lemma.

\proclaim{Lemma 1.2} If the limit 3-orbifold $C_{\infty}$ is not compact, then
\roster
\item"i)" $C_{\infty}$ has a finite volume,
\item"ii)" the singular set $\Sigma_{\infty}$ of $C_{\infty}$ is not
compact.
\endroster
\endproclaim

\demo{Proof} Assertion i) is Lemma 1.3 of Chapter V and assertion ii) is Proposition
1.4 of Chapter V, whose proofs do not require the cone angles to be strictly
less  than $\pi$ but only less than or equal to $\pi$. \qed
\enddemo

Hence, the non-compact orbifold $C_{\infty}$ is hyperbolic with finite
volume. Let $N_{\infty}\subset C_{\infty}$ be a compact core corresponding
to the thick part of the orbifold. The thin part $C_{\infty}-N_{\infty}$ is
a union of cusps of the form $F\times (0,+\infty)$, where $F$ is an
orientable closed $2$-dimensional Euclidean orbifold. Moreover, since
$\Sigma_{\infty}$ is not compact, at least one of the cusps is singular.

Propostion 1.1, in the case where $C_{\infty}$ is not compact, follows from the
following one.

\proclaim{Proposition 1.3} 
Let $N_{\infty}\subset C_{\infty}$ be the compact core of the hyperbolic
3-orbifold $C_{\infty}$. Then $N_{\infty}$ embeds in $\OO$ as an essential
3-suborbifold.
\endproclaim

\demo{Proof} The geometric convergence implies that, for $n$ sufficiently large, there is a
$(1+\varepsilon_n)$\-bilipschitz embedding
$f_n:(N_{\infty},\Sigma_{\infty}\cap N_{\infty})\to (C_n,\Sigma)$ with
$\varepsilon_n\to 0$. Since the $3$-orbifold $\OO$ and the cone $3$-manifolds $C_n$
have the same topological type, we view the image $f_n(N_{\infty})$ as a
suborbifold of $\OO$, which we denote by $N_n\subset\OO$.
 The  orbifold $N_n$ is homeomorphic to $N_{\infty}$, thus
$N_n$ is an orbifold whose interior is hyperbolic of finite volume.

In Lemma 1.5 we are going to prove that  $\partial N_n$ is incompressible
in
$\OO$, but before we need the
following lemma.

\proclaim{Lemma 1.4} For $n$ sufficiently large, the orbifold $\OO-\Int(N_n)$ is
irreducible. 
\endproclaim

\demo{Proof}
Seeking a contradiction, we assume that $\OO-\Int(N_n)$  contains a
 spherical $2$-suborbifold
$F^2$ which is essential.
Since $\OO$ is irreducible, $F^2$ bounds a spherical $3$-orbifold $\Delta^3$,
which is the quotient of a standard $3$-ball $B^3$ by the ortogonal action
of a finite subgroup of $SO(3)$. Since the singular set of $\OO$ is a link,
the only possibility is that $\Delta^3$ is either a non-singular ball $B^3$ or its
quotient by a finite cyclic group. Therefore the topological type of $\Delta^3$ is 
$(B^3,A)$, where $A=\Delta^3\cap\Sigma$ is either empty or an unknotted 
proper arc.
Moreover, by hypothesis, $N_n\subset \Delta^3$, for $n$ sufficiently large.

Let $\rho_n:\pi_1(C_n-\Sigma, x_n)\to PSL_2(\Bbb C)$ be the holonomy
representation of $C_n$ and let 
$f_n:(N_{\infty},\Sigma_{\infty}\cap N_{\infty})\to (C_n,\Sigma)$ be the
$(1+\varepsilon_n)$-bilipschitz embedding such that $N_n=f_n(N_{\infty})$. 
For $n$ sufficiently large, the representation $\rho_n\circ
f_{n*}:\pi_1(N_{\infty}-\Sigma_{\infty},x_\infty)\to PSL_2(\Bbb C)$ is
either cyclic or
trivial, since $N_n\subset \Delta^3$. Hence, the holonomy of
$C_{\infty}$ is abelian, because  the geometric convergence implies the
convergence of the holonomies (Proposition 5.4, Chapter III). This contradicts the fact
that $C_{\infty}$ is a complete hyperbolic orbifold of finite volume.
\qed
\enddemo

\proclaim{Lemma 1.5} For $n$ sufficiently large, the boundary $\partial N_n$ is
incompressible in $\OO$.
\endproclaim

\demo{Proof}
Seeking a contradiction, we suppose that the lemma is not true. So, after passing to a
subsequence if necessary, we can assume that
$\partial N_n$ is compressible in $\OO$ and furthermore that  $\OO-\Int(N_n)$ is
irreducible (by Lemma 1.4). Let $F_1,\ldots,F_p$ be the components of
$\partial N_{\infty}$. By  passing again to a subsequence if necessary, we can assume moreover
that the embedded components $f_n(F_1),\ldots,f_n(F_q)$ are precisely the 
compressible ones, with $p\geq q$, where
$f_n:(N_{\infty},\Sigma_{\infty}\cap N_{\infty})\to (C_n,\Sigma)$ is the
$(1+\varepsilon_n)$\-bilipschitz embedding defining $N_n$. 

 For $i=1,\ldots,q$, let $\lambda^i_n$ be an essential curve on
$f_n(F_i)$ which bounds a  properly embedded disk in $\OO-\Int(N-n)$,
intersecting $\Sigma$ in at most
one point. Consider the simple closed essential curves 
$\tilde\lambda^i_n=f^{-1}_n(\lambda^i_n)\subset F_i$, for 
$i=1,\ldots,q$.

\proclaim{Claim 1.6} For each $i=1,\ldots,q$, the sequence of simple
closed essential curves $(\tilde\lambda^i_n)_{n\geq n_0}$ represents infinitely
many different homotopy classes in the  fundamental orbifold group
$\pi_1^o(F_i)$.
\endproclaim

\demo{Proof}
If the claim is false, then, by passing to a subsequence and changing the
indices of the $F_i$, we can suppose that the curves $\tilde\lambda^1_n$
represent a fixed class $\tilde\lambda^1\in\pi_1^o(F_1)$ which does not
depend on $n$. Let 
$\rho_n:\pi_1(C_n-\Sigma, x_n)\to PSL_2(\Bbb C)$ be the holonomy
representation of $C_n$ and 
$\rho_{\infty}:\pi_1(C_{\infty}-\Sigma_{\infty}, x_{\infty})
=
\pi_1(N_\infty-\Sigma_\infty,x_{\infty})\to PSL_2(\Bbb
C)$ be the holonomy  of $C_{\infty}$. From the geometric convergence (Proposition 5.4, Chapter
III):
$$
\rho_{\infty}(\tilde\lambda^1)
=\lim_{n\to\infty}\rho_n(f_n(\tilde\lambda^1_n))
=\lim_{n\to\infty}\rho_n(\lambda^1_n).
$$
Since the curves $\lambda^1_n$ are compressible in $\OO$, their holonomies 
$\rho_n(\lambda^1_n)$ are either elliptic or trivial. Moreover, since 
$\pi_1^o(F_1)$ is parabolic and $\tilde\lambda^1$ is essential in $F_1$,
the holonomy $\rho_{\infty}(\tilde\lambda^1)$ is non-trivial and
parabolic.
 We remark
that possibly  some non-trivial elements in $\pi^o_1(F_1)$ have an elliptic
holonomy, but they are not represented by an essential simple closed curve
in $F_1$. Thus we obtain a contradiction comparing
$\rho_{\infty}(\tilde\lambda^1)$ with the limit of $\rho_n(\lambda^1_n)$.
\qed
\enddemo

We come back to the proof of Lemma 1.5.

Since $\OO$ is a very good orbifold, it has a finite covering $M$ which is
a (closed oriented) 3-manifold. The 3-suborbifolds
$N_n\subset\OO$ lift to compact 3-submanifolds $\bar N_n\subset M$, all of
them  being homeomorphic to a 3-manifold $\bar N_{\infty}$ covering $N_{\infty}$.
Hence,
$\Int(\bar N_n)\cong\Int(\bar N_{\infty})$ is complete hyperbolic with finite volume
and $\partial\bar N_n$ is a union of tori. 
Let $T_1,\ldots, T_r$ be the
components of $\partial \bar N_{\infty}$ which cover the components  $f_n(F_1),\ldots,f_n(F_q)$ 
of $\partial N_{\infty}$; they are compressible for the lifted
embeddings $\bar f_n:N_{\infty}\to M$, such that $\bar N_n=\bar f_n(\bar
N_{\infty})$.

Since the orbifold $\OO-\Int(N_n)$ is irreducible, the Equivariant Sphere
Theorem (\cite{\DD}, \cite{\MYOne,2}, \cite{\JR}) implies that
$M-\Int(N_n)$ is also irreducible. Hence, for $i=1,\ldots,r$ , the tori
$\bar f_n(T_i)\subset \partial\bar N_n$ bound solid tori in $M$. Let
$\mu^i_n\subset \bar f_n(T_i)$ be the boundary of a meridian disk of
this solid torus, for $i=1,\ldots,r$;
the curve $\mu^i_n$ is a lift to $M$ of the compressing curve
$\lambda^j_n\subset f_n(F_j)\subset\partial N_k$ (where $F_j$ is the $\partial$-component covered
by $T_i$).
 
Claim 1.6 implies that, for each $i=1,\ldots,r$, the curves
$\tilde\mu^i_n=\bar f_n^{-1}(\mu^i_n)\subset T_i$ represent infinitely many
different homotopy classes in $\pi_1(T_i)$.  Since $\bar N_{\infty}$ is
hyperbolic, Thurston's hyperbolic Dehn filling theorem \cite{\ThuNotes} implies that, for $n$
sufficiently large, the 3-manifold
$$
\bar N_{\infty}(\tilde \mu^1_n,\ldots,\tilde\mu^r_n)=
\bar N_{\infty}\cup\bigsqcup_{i=1}^r S^1\times D^2_i
$$
obtained by Dehn filling along the curves $\tilde
\mu^1_n,\ldots,\tilde\mu^r_n$ is hyperbolic.  In particular, 
$\bar N_{\infty}(\tilde \mu^1_n,\ldots,\tilde\mu^r_n)$ has incompressible
boundary.

Since the curves $\tilde\mu^i_n=\bar f_n^{-1}(\mu^i_n)\subset T_i$
 represent infinitely many
different homotopy classes in $\pi_1(T_i)$,
it follows from volume estimations that the sequence  $(\bar
N_{\infty}(\tilde \mu^1_n,\ldots,\tilde\mu^r_n))_{n\in\N}$ contains infinitely many 
non-homeomorphic 3-manifolds.
We shall obtain a contradiction by showing that in fact all these manifolds are 
homeomorphic to finitely many ones.
For $n$ large, the boundary $\partial\bar
N_{\infty}(\tilde \mu^1_n,\ldots,\tilde\mu^r_n)$ is incompressible in $M$,
because $M-\Int(\bar
N_{\infty}(\tilde \mu^1_n,\ldots,\tilde\mu^r_n))$ is irreducible with
incompressible boundary and $\bar
N_{\infty}(\tilde \mu^1_n,\ldots,\tilde\mu^r_n)$ is hyperbolic. Hence, this
3-submanifold is a piece of the Jaco-Shalen \cite{\JS} and Johannson
\cite{\Joh} splitting
 of the 3-manifold M. Unicity of this splitting implies that the 3-manifolds 
$\bar N_{\infty}(\tilde \mu^1_n,\ldots,\tilde\mu^r_n)$ are only finitely many. Hence we get
the contradiction that proves Lemma 1.5 and therefore Proposition 1.3.\qed
\enddemo
\enddemo
\enddemo

\head
2. The collapsing case
\endhead

Next proposition proves Theorem B in the collapsing case.

\proclaim{Proposition 2.1}
Let $\OO$ and $(C_n)_{n\in\N}$ satisfy the hypothesis of Theorem B.
 If the sequence $(C_n)_{n\in\N}$ 
 collapses, then $\OO$  contains a non-empty compact essential
$3$-suborbifold,  which is not a product and which is either Euclidean, Seifert or Sol.
\endproclaim

\demo{Proof} Since the sequence  $(C_n)_{n\in\N}$  collapses,
by Corollary 4.1 of the Local Soul theorem (Chapter IV), either
 there is a subsequence that, after
rescaling, converges to a closed Euclidean cone 3-manifold, or the Local Soul theorem, with any
parameters $\varepsilon>0$, $D>0$ and non-compact local models, applies to every point  $x\in
C_n$. Thus we distinguish again two cases, according to whether we obtain a compact limit or not.

In the first case, for every $n\in\Bbb N$ there is $x_n\in C_n$ such that the sequence of rescaled
cone 3-manifolds
$(\bar C_n,x_n)=(\tfrac1{\Inj(x_n)} C_n,x_n)$ has a subsequence that converges geometrically to
a compact cone 3-manifold $(\bar C_{\infty}, x_{\infty})$. The geometric convergence
implies that $\bar C_{\infty}$ is a closed Euclidean 3-orbifold  with the same
topological type and the same branching
indices as  $\OO$. Therefore, as an orbifold $C_{\infty}=\OO$ and so $\OO$ is Euclidean. Thus 
Proposition 2.1  holds in this case.

The second case, when we  cannot find such a compact limit,
is the difficult case to which the remaining of this chapter is devoted.
Hence, from now on
 we suppose that the Local Soul theorem, with any
parameters $\varepsilon>0$, $D>0$ and non-compact local models, applies to every point  $x\in
C_n$, for n sufficiently large.

\proclaim{Lemma 2.2} Under the hypothesis of the second case, for any
$\varepsilon>0$ and
$D>0$, there exists $n_0>0$ such that, for $n\geq n_0$, every $x\in C_n$ has an
open neighborhood $U_x$ $(1+\varepsilon)$-bilipschitz
homeomorphic to  the normal fiber bundle $\NN_{\nu}(S)$, with
radius $\nu<1$, of the soul $S$ of one of the following non-compact orientable Euclidean
orbifolds:
\roster
\item"a)" $T^2\times \Bbb R$;  $S^1\ltimes\Bbb E^2$;
$S^1\ltimes D^2(\frac{2\pi}p)$;
\item"b)" 
 $S^2(\frac{2\pi}{p_1},\frac{2\pi}{p_2},\frac{2\pi}{p_3})\times\Bbb R$,
 with 
 $\frac1{p_1}+\frac1{p_2}+\frac1{p_3}=1$; $S^2(\pi,\pi,\pi,\pi)\times\Bbb R$; the
pillow;  
\item"c)" $\Bbb {P}^2(\pi,\pi)\tilde\times\Bbb R$, which is the twisted orientable line bundle over
$\Bbb P^2(\pi,\pi)$; and the quotient of  $S^2(\pi,\pi,\pi,\pi)\times\Bbb R$ by an involution
 that
gives the orientable  bundle over $D^2(\pi,\pi)$, with silvered boundary (cf. Figure VI.1). 
\endroster
Moreover,  if $f:U_x\to\NN_{\nu}(S)$ is the 
$(1+\varepsilon)$-bilipschitz homeomorphism, then 
$$
\max(\Inj(x),d(f(x),S),\Diam(S))\leq\nu/D.
$$
\endproclaim

\midinsert
 \centerline{\BoxedEPSF{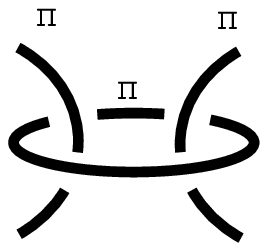 scaled 750}}
    \botcaption{Figure VI.1}
    \endcaption
     \endinsert

Recall that the pillow is the orbifold with  underlying space $\Bbb R^3$
and branching set two straigh lines of branching order 2 (cf. Figure IV.1). It is the
quotient of
$S^1\ltimes \Bbb R^2$ by an involution.

\demo{Proof}  From the hypothesis of the second case, for every $\varepsilon>0$ and $D>0$, there
exists an $n_0$ such that for $n>n_0$  we can apply  the Local Soul Theorem, with 
parameters $\varepsilon>0$, $D>0$ and non-compact local models, to every point  $x\in
C_n$. Moreover from the hypothesis about the cone angles, the local
models are Euclidean non-compact 3-orbifolds.
Now it remains to eliminate the
Euclidean 3-orbifolds model that are no listed in Lemma 2.2  (as in Lemma 2.2 of Chapter V).
Hence, by the Local Soul theorem, we only have to get rid of the
twisted line bundle over the Klein bottle $K^2\tilde\times\Bbb R^1$ and the two
models of Figure VI.2, which correspond  to an orientable  bundle over either an annulus or a
M\"obius strip, with silvered boundary in both cases.

\midinsert
 \centerline{\BoxedEPSF{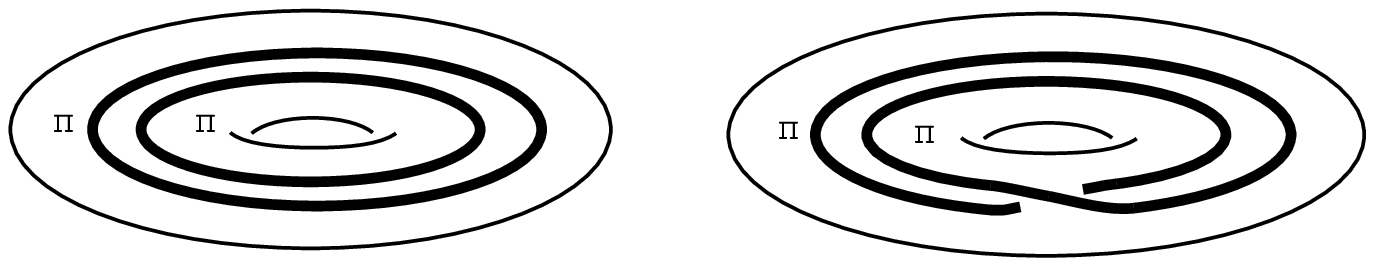 scaled 750}}
   \botcaption{Figure VI.2}
    \endcaption
     \endinsert

Let $S$ be the soul of one of these three Euclidean non-compact orbifolds, and let
$\NN_r(S)$ denote its normal fiber bundle of radius $r$. Then $\partial\NN_r(S)$
is an incompressible torus in $\NN_r(S)-\Sigma$, for every $r>0$. Therefore the appearance of one of
these models  would contradict the fact that $C-\Sigma$ is topologically atoroidal and not Seifert
fibered.\qed

\enddemo

As in the previous chapter, the neighborhoods given by Lemma 2.2 are called
{\it $(\varepsilon,D)$-Margulis' neighborhoods.}
 We remark that the $(\varepsilon,D)$-Margulis' neighborhoods of abelian
type correspond to the local models listed in a). The local models listed in c) are
Seifert fibered and different from a product. If neighborhoods corresponding to these local models
appear, the following lemma proves proposition 2.1.

\proclaim{Lemma 2.3} If for every $n\geq n_0$ there is a point $x\in
C_n$ having a $(\varepsilon,D)$-Margu\-lis' neighborhood of  type c)
in Lemma 2.2, then the orbifold $\OO$ contains a non-empty compact essential 
orientable $3$-suborbifold $\OO'$ which is Seifert fibered and 
different from a product.
\endproclaim

\demo{Proof of Lemma 2.3}
Let $S$ be the soul of one of the Euclidean local models given in c). This soul is
either a projective plane with two cone points $\Bbb P^2(\pi,\pi)$ or a disk
with two cone points and mirror boundary $\bar D^2(\pi,\pi)$. In both cases
a regular neighborhood $\NN(S)$ of $S$ embeds as a compact suborbifold of
$\OO$. This suborbifold $\OO'=\NN(S)$ is Seifert fibered, it is not a
product and its boundary $\partial\OO'=S^2(\pi,\pi,\pi,\pi)$ is
incompressible in $\OO'$. It remains to show that either it is also
incompressible in $\OO$ or $\OO$ is Seifert fibered itself. 

First note that $\OO-\Int(\OO')$ is irreducible, because the
soul $S=\Bbb P^2(\pi,\pi)$ or  $\bar D^2(\pi,\pi)$ cannot be contained in
the quotient of a ball by a cyclic action. If $\partial\OO'$ is compressible
in $\OO$, then the orbifold 
$\OO-\Int(\OO')$ is irreducible  with compressible boundary
$\partial (\OO-\Int(\OO'))=S^2(\pi,\pi,\pi,\pi)$. Therefore 
$\OO-\Int(\OO')$ is a pillow and $\OO$ is Seifert fibered.
\qed
\enddemo

Lemma 2.3 shows that in applying Lemma 2.2 we only need to
consider local models of types a) and b). The
following Lemma shows that we must consider local models of type b). 

\proclaim{Lemma 2.4} There is a constant $D_0>0$ such that  
if every point of a closed hyperbolic cone 3-manifold $C$ has a
$(\varepsilon,D)$-Margulis' neighborhood of type a) or b), 
with $\varepsilon<1/2$ and
$D>D_0$,
then at least one of the neighborhoods is of type b).
\endproclaim 

\demo{Proof} 
By Proposition 2.3 of Chapter V, there exists a uniform constant $D_0>0$ such that,
if every point of a closed hyperbolic cone 3-manifold 
 $C$ has a $(\varepsilon,D)$-Mar\-gu\-lis' neighborhood of type a)
(abelian type), whith $\varepsilon<1/2$ and $D>D_0$, then the simplicial
volume $\norm{C-\Sigma}=0$. Therefore, there must be a point whose local model 
is of type b), because the fact that $C-\Sigma$ admits a complete hyperbolic structure (Chapter V,
Lemma 1.5) implies that $\norm{C-\Sigma}\neq 0$.
\qed
\enddemo

By using Lemmas 2.2, 2.3 and 2.4, Proposition 2.1 follows from Proposition 2.5 below:

\definition{Definition} We say that a compact oriented 3-orbifold
$\OO$ is a \it graph orbifold \rm if there exists a family of orientable
Euclidean closed $2$-suborbifolds that decompose $\OO$ into Seifert fibered
pieces. In particular, Seifert fibered 3-orbifolds are graph orbifolds. 
\enddefinition

\proclaim{Proposition 2.5} 
Let 
$(C_n)_{n\in\Bbb N}$ and $\OO$ satisfy the hypothesis of Theorem B. There is a universal
constant $D_1>0$ such that, if for some $n$ every point of $C_n$
admits a $(\varepsilon,D)$-Margulis' neighborhood of type a) or
b), with  $\varepsilon<1/2$ and $D>D_1$, then $\OO$ is a graph
orbifold.
 \endproclaim

\demo{Proof of Proposition 2.5} Let $0<\varepsilon<1/2$ and $D>D_0$,
where $D_0$ is the constant of Lemma 2.4. Assume that, for some $n$ fixed, 
every point of the hyperbolic cone 3-manifold $C_n$ has a 
$(\varepsilon,D)$-Margulis' neighborhood of type a) or b).

We choose a
point $x_0\in C_n$ having a local model of type b). It means that $x_0$ has a
neighborhood $U_{x_0}\subset C_n$ with a $(1+\varepsilon)$-bilipschitz
homeomorphism $f_0:U_{x_0}\to \NN_{\nu}(S)$, where $\NN_{\nu}(S)$ is the normal fiber bundle,
with radius $\nu<1$, of the soul $S$ of a non-compact Euclidean $3$-orbifold of the family b).

Let $W_0=f_0^{-1}\big(\overline{\NN_{\nu_0/D}(S)}\big)\subset
U_{x_0}$ be the inverse image  of the closed normal fiber bundle of the
soul $S$, with radius $\nu_0/D$. We view $W_0$ as a suborbifold of $\OO$.

Since $\OO$ is very good, there is a regular finite covering $p:M\to\OO$ such
that $M$ is a closed 3-manifold. We need the following  proposition, that we shall prove in the
next section. 
We recall that $D_0$ is the universal constant of Lemma 2.4; we can suppose $D_0>10^4$.

\proclaim{Proposition 2.6} 
Let $\OO$ and
$(C_n)_{n\in\Bbb N}$ satisfy the hypothesis of Theorem B. There is a universal
constant $b_1>0$ such that, if for some $n$ every point of $C_n$
admits a $(\varepsilon,D)$-Margulis' neighborhood of type a) or
b), with  $\varepsilon<1/2$ and $D>D_0>10^4$, then 
$C_n$ admits a $\eta$-covering \`a la Gromov $(V_i)_{i\in I}$ with
$\eta<b_1/D$. 

Moreover, there is a choice of $x_0$ and $W_0\subset \OO$ such that the covering
$(V_i)_{i\in I}$  satisfies the additional properties:
 \roster
\item"6)" for every $i\in I$, 
$p^{-1}\Big(   \bigcup\limits_{V_j\cap
V_i\neq\emptyset} V_j    \Big)$
is abelian in
$M-p^{-1}(W_0)$;
\item"7)" $W_0$ intersects only one open set $V_i$ of the
covering.
 \endroster
 \endproclaim

We recall that  $U\subset M$ is \it abelian in
$M-p^{-1}(W_0)$  \rm if the homomorphism induced by  the inclusion
$$
i_*:\pi_1(U-p^{-1}(W_0))\to
\pi_1(M-\Int(p^{-1}(W_0)))
$$ 
has  an abelian image.

\demo{Proof of Proposition 2.5 assuming Proposition 2.6} 
Set $\widetilde W_0=p^{-1}(W_0)\subset M$.
Let
$\eta_0>0$ be the univeral constant of Proposition 3.1 of Chapter V. We choose
$D_1=\sup(b_1/\eta_0,10^4)$. Proposition 2.6 of this chapter and  Proposition 3.1 of Chapter V 
imply the
existence of a continuous map $g:C\to K^2$, from $C$ to a simplicial $2$-complex  $ K^2$, such
that:
\roster
\item"i)" for every vertex $v$ of $K^2$, $p^{-1}(g^{-1}(\Star{v}))$ is
abelian in $M-\widetilde W_0$;
\item"ii)" $g(W_0)$ is a vertex of $K^2$.
\endroster

By composing $g$ with the projection of the covering map $p:M\to\OO$, we have a
continuous map $f=g\circ p:M\to K^2$ with the following properties:
\roster
\item"i)" for every vertex $v$ of $K^2$, $f^{-1}(\Star{v})$ is
abelian in $M-\widetilde W_0$;
\item"ii)" $f(\widetilde W_0)$ is a vertex $v_0$ of $K^2$.
\endroster

Now we use the map  $f$ to show that all Dehn fillings along the boundary of any connected component
of  $M-\Int(\widetilde W_0)$ have simplicial volume zero. Let $N$ be a
connected component of $M-\Int(\widetilde W_0)$. Its boundary $\partial N$ is a union of tori. Let
$$
\bar N=N\underset{\partial N}\to\cup\bigsqcup_{i=1}^p D^2\times S^1
$$
be any closed Dehn filling of $N$ along $\partial N$. Since $f(\widetilde W_0)$ is a vertex
$v_0$ of $K^2$, the map
$f:M\to K^2$ induces a map $\bar f:\bar N\to K^2$  that 
coincides with $f$ in $N$ and maps each filling solid torus 
$D^2\times S^1$ to the vertex $v_0$.
By property i), for every vertex $v\in K^2$, $\bar f^{-1}(\Star v)$ is
abelian in $\bar N$. 
Hence, the closed oriented 3-manifold $\bar N$
admits an abelian covering of dimension $2$. By Gromov's Vanishing
Theorem \cite{\Gro, Sec. 3.1} (cf \cite{\Iva})  the simplicial
volume $\norm{\bar N}=0$, as claimed.

Next step is the following lemma, whose proof is postponed to the end of the
section.

\proclaim{Lemma 2.7} The 3-manifold $M-\Int(\widetilde W_0)$ is irreducible.
\endproclaim

Assuming this lemma, any connected component $N$ of $ M-\Int(\widetilde
W_0)$ is irreducible. If the boundary
$\partial N$ is incompressible in $N$, then, by \cite{\BDV}, $\norm N=0$. If
$\partial N$ is compressible in $N$, then $N$ is a solid torus because it is
irreducible; in particular $\norm N=0$. Therefore $\norm {M-\Int(\widetilde
W_0)}=0$.

We prove now that $M-\Int(\widetilde W_0)$ is a graph manifold.
Since $M-\Int(\widetilde W_0)$ is compact, irreducible and
with toral boundary, according  to Jaco-Shalen \cite{\JS} and Johannson
\cite{\Joh}  $M-\Int(\widetilde W_0)$ splits along
incompressible tori into  Seifert and simple pieces. By Thurston's
Hyperbolization Theorem the simple pieces are hyperbolic, hence they have
non-zero simplicial volume (\cite{\Gro} and \cite{\Thu, Ch. 6}). Since the
simplicial volume is additive under gluing along incompressible tori (\cite{\Gro} and
\cite{\Som}), the fact that $\norm{M-\Int(\widetilde W_0)}=0$ implies that all the pieces in the
JSJ-splitting are Seifert; therefore $M-\Int(\widetilde W_0)$ is a graph manifold.

We deduce that $\OO-W_0$ is
a graph orbifold, 
because $p:M-\Int(\widetilde W_0)\to \OO-W_0$ is a regular covering 
 and the results of Meeks and Scott \cite{MS} provide a graph structure on
$M-\Int(\widetilde W_0)$ invariant by the action of the deck transformations
group of the covering.

The $3$-suborbifold $W_0\subset\OO$ is either $
S^2(\frac{2\pi}{p_1},\frac{2\pi}{p_2},\frac{2\pi}{p_3})\times [0,1]$ (with
$\frac1{p_1}+\frac1{p_2}+\frac1{p_3}=1$),
$S^2(\pi,\pi,\pi,\pi)\times[0,1]$ or the pillow. Therefore $W_0$ is
Seifert fibered and $\OO$ admits a graph structure.

This proves  Proposition 2.5 from Proposition 2.6 and Lemma 2.7. The
proof of Proposition 2.6 is given in the next section and the proof of Lemma
2.7 comes now.

\demo{Proof of Lemma 2.7}
By the Equivariant Sphere Theorem
(\cite{\DD}, \cite{\MYOne,2}, \cite{\JR}), to
show that  the 3-manifold $M-\Int(\widetilde W_0)$ is irreducible it suffices
to prove the irreducibility of the orbifold $\OO-W_0$.

 Seeking a
contradiction, we assume that $\OO-W_0$ is reducible. It means that there exists an
essential spherical $2$-suborbifold $F^2\subset\OO$. Since $\OO$ is
irreducible, $F^2$ bounds a discal 3-suborbifold  $D^3$ in $\OO$ and
$D^3$ contains $W_0$. Since the branching locus $\Sigma\subset\OO$ is a link,
the  discal suborbifold $D^3$ is the quotient of a ball by a finite cyclic orthogonal
action. Hence the topological type of $D^3$ is $(B^3,A)$, where
$A=B^3\cap\Sigma$ is a proper unknotted arc in the ball $B^3$.

Since $W_0\subset \Int(D^3)$, we already have a contradiction 
in the case where $W_0=\NN
(S^2(\frac{\pi}{p_1},\frac{\pi}{p_2},\frac{\pi}{p_3}))$, because there
is no way to embedd a $2$-sphere in $B^3$ that intersects $A$ in 3
points. 

Hence we assume that $W_0=\NN(S^2(\pi,\pi,\pi,\pi))$  or $W_0$ is
the pillow. In both cases, $\Sigma\cap W_0$ is not  connected and we
find a contradiction using a
Dirichlet polyhedron and the fact that $A=B^3\cap\Sigma$ is connected.
More precisely, these local models imply that there is a metric ball
$B(x,r)\subset W_0\subset C_n$ such that $B(x,r)\cap\Sigma$ is not
connected. We consider  the Dirichlet polyhedron $P_x$ of $C_n$ 
centered at $x$. This polyhedron is convex, because the cone
angles of $C_n$ are equal to or less than $\pi$.  By convexity,
different connected components of  $B(x,r)\cap\Sigma$ give different
edges of $\partial P_x$ that belong to different geodesics of $\Bbb
H^3$. In particular,  the holonomy of the meridians of different
components of $B(x,r)\cap\Sigma$ are not contained in a cyclic group.
This contradicts the inclusion $(W_0,\Sigma\cap W_0)\subset
(D^3,\Sigma\cap D^3)$, because $\pi_1(D^3-\Sigma)\cong\pi_1(B^3-A)$ is
cyclic. Thus we get a contradiction and the lemma is proved.  \qed
 \enddemo

 \enddemo

\enddemo

\enddemo

\head  
      3. From $(\varepsilon,D)$-Mar\-gu\-lis' coverings  of 
type a) and b) to
$\eta$-coverings \`a la Gromov 
\endhead

This section is devoted to the proof of Proposition 2.6, which constructs
the  required $\eta$-covering \`a la Gromov.

We recall that we had applied the Local Soul theorem (Chapter IV) to the hyperbolic cone 3-manifold
$C_n$ with parameters $(\varepsilon,D)$. For any point $x_0\in C_n$
whose local model is of type b), there is a neighborhood $U_{x_0}$ and 
a $(1+\varepsilon)$-bilipschitz
homeomorphism $f_0:U_{x_0}\to \NN_{\nu_0}(S)$ 
where $\NN_{\nu_0}(S)$ is a normal fiber bundle, with radius $\nu_0 \leq 1$, of the
soul $S$ of a non-compact Euclidean cone 3-manifold of type b). We have defined
 $W_0=f_0^{-1}\big(\overline{\NN_{\nu_0/D}(S)}\big)\subset
U_{x_0}$.

Let $p:M\to\OO$ be a finite regular covering of $\OO$ that is a closed
$3$-manifold. The proposition we want to prove is the following.

\proclaim{Proposition 2.6} 
Let $\OO$ and
$(C_n)_{n\in\Bbb N}$ satisfy the hypothesis of Theorem B.
 There is a universal
constant $b_1>0$ such that, if for some $n$ every point of $C_n$
admits a $(\varepsilon,D)$-Margulis' neighborhood of type a) or
b), with  $\varepsilon<1/2$ and $D>D_0>10^4$, then 
$C_n$ admits a $\eta$-covering \`a la Gromov $(V_i)_{i\in I}$ with
$\eta<b_1/D$. 

Moreover, there is a choice of $x_0$ and $W_0\subset \OO$ such that the covering
$(V_i)_{i\in I}$  satisfies the additional properties:
 \roster
\item"6)" for every $i\in I$, 
$p^{-1}\Big(   \bigcup\limits_{V_j\cap
V_i\neq\emptyset} V_j    \Big)$
is abelian in
$M-p^{-1}(W_0)$;
\item"7)" $W_0$ intersects only one open set $V_i$ of the
covering.
 \endroster
 \endproclaim

\demo{Proof of Proposition 2.6}
In the proof we set $C_n=C$ to simplify notation.

First we  describe the choices of $x_0\in C$ and $W_0$.
Given $\varepsilon>1/2$ and $D>D_0>10^4$, we consider
$$
 T_{(\varepsilon,D)}=\left\{ x\in C\ \Big\vert \matrix \format\l\\
x \text{ admits an } (\varepsilon, D) \text{-Margulis'}\\
\text{neighborhood of type b)}
\endmatrix
\right\}
$$
Since $D>D_0$, Lemma 2.4 implies that $T_{(\varepsilon,D)}\neq\emptyset$.
For $x\in T_{(\varepsilon,D)}$, let $U_x$ denote the
${(\varepsilon,D)}$-Margulis' neighborhood of type b) and let
$f:U_x\to \NN_{\nu(x)}(S)$ be the $(1+\varepsilon)$-bilipschitz homeomorphism
between $U_x$ and the normal fiber bundle, with radius $\nu(x)\leq 1$, of the
compact soul $S$ of a local model of type b). We choose a
point $x_0\in T_{(\varepsilon,D)}$ such that 
$$
\nu(x_0)=\nu_0\geq\frac1{1+\varepsilon}\sup 
\{\nu(x)\mid x\in T_{(\varepsilon,D)}\}.
 $$
Let $W_0=f_0^{-1}\big(\overline{\NN_{\nu_0/D}(S)}\big)\subset
U_{x_0}$ be the inverse image  of a closed normal neighborhood of the
soul $S$ of radius $\nu_0/D$, where 
$f_0:U_{x_0}\to \NN_{\nu_0}(S)$ is the $(1+\varepsilon)$-bilipschitz
homeomorphism.

For every $x\in C$ we define the abelianity radius (relative to
$\widetilde W_0=p^{-1}(W_0)\subset M$):
$$
\ab x=\sup\{r\in\Bbb R\mid p^{-1}(B(x,r))\text{ is abelian in }
M-\widetilde W_0\}.
 $$
We set $r(x)=\inf\{ 1,\frac{\ab x}8 \}$.
 
This definition is analogous to the one given in Section V.4. For
instance, the following lemma has the same proof as Lemma 4.2 of Chapter V:

\proclaim{Lemma 3.1}   Let $x,y\in C$. If $B(x,r(x))\cap
B(y,r(y))\neq\emptyset$, then \roster
\item"(a)" $3/4\leq{r(x)}/{r(y)}\leq 4/3$;
\item"(b)" $B(x,r(x))\subset B(y,4r(y))$. \qed
\endroster
\endproclaim

\proclaim{Lemma 3.2} For every $x_0\in W_0$, $W_0\subset
B(x_0,\frac{r(x_0)}9)$.
\endproclaim

\demo{Proof}
This lemma is equivalent to the inequality 
$$
\Diam(W_0)< r(x_0)/9.
$$
Since  $W_0=f_0^{-1}\big(\overline{\NN_{\nu_0/D}(S)}\big)\subset
U_{x_0}$, where $f_0:U_{x_0}\to \NN_{\nu_0}(S)$ 
is a $(1+\varepsilon)$-bilipschitz homeomorphism,
and $\NN_{\nu_0}(S)$ is a normal fiber bundle, with radius $\nu_0$, of 
the soul $S$ of a non-compact Euclidean cone 3-manifold of type b), we have:
 $$
\Diam (W_0)\leq (1 +\varepsilon)\Diam\big( \NN_{\frac{\nu_0}D}(S)\big)
\leq (1+\varepsilon) (\Diam (S)+ 2 \,\frac{\nu_0}D)\leq 6\,\frac{\nu_0}D,
$$
because 
$\Diam(S)\leq\frac{\nu_0}D$ and $\varepsilon\leq 1/2$.
By definition 
$
\ab{x_0}\geq\frac1{1+\varepsilon}(\nu_0-\frac{\nu_0}D)
\geq\frac{\nu_0}2,$ moreover 
 $\nu_0\leq 1$ and $D>10^4$, thus we obtain the following inequalities: 
$$
\Diam(W_0)\leq 6\,\frac{\nu_0}D \leq\inf\{
\frac6D,\frac{12\ab{x_0}}D\}<\frac{r(x_0)}9.\qed  
$$
\enddemo

Now we give the construction of the $\eta$-covering
\`a la Gromov. We fix a point $x_0\in W_0$, we consider then all the
possible finite sequences $\{x_0,x_1,\ldots,x_p\}$, starting with $x_0$, such that:
$$
\text{the balls }
B\Big(x_0,\frac{r(x_0)}4\Big),\ldots,B\Big(x_p,\frac{r(x_p)}4\Big)
\text{ are pairwise disjoint}. \tag*
  $$
A sequence satisfying \ttag* and Lemma 3.1 is finite by compactness. 
Moreover we have the following property, proved in Chapter V, Lemma 4.4.

\proclaim{Lemma 3.3} If the sequence $\{x_0,x_1,\ldots,x_p\}$ is maximal
for property \ttag*, then the balls $B(x_0,\frac23 r(x_0)),\ldots,
B(x_p,\frac23 r(x_p))$ cover $C$.\qed
\endproclaim

Given a sequence $\{x_0,x_1,\ldots,x_p\}$,  maximal
for property \ttag* and starting with $x_0\in W_0$ , we consider the  covering
of $C$ by the following open sets :
$$
\left\{\matrix \format\l&\qquad\l\\ V_0=B(x_0,r(x_0))&\\
V_i=B(x_i,r(x_i))-W_0, &\text{ for }
i=1,\ldots,p.
\endmatrix
\right.
$$ 

Next lemma concludes the proof of Proposition 2.6.

\proclaim{Lemma 3.4} There is a universal constant $b_1>0$ such that, for
$\varepsilon<1/2$ and $D>10^4$, the  open sets $V_0,\ldots,V_p$ define a
$\eta$-covering \`a la Gromov of $C$, with $\eta<b_1/D$. Moreover this covering
satisfies properties 6) and 7) of Proposition 2.6. 
\endproclaim

\demo{Proof} Lemmas 3.2  and 3.3 garantee that the open  sets
$V_0,\ldots,V_p$ cover $C$. Then by setting $r_i=r(x_i)$ for $i=1,\ldots,p$,
properties 1), 2) and 3) of a $\eta$-covering \`a la Gromov follow from the
construction and Lemma 3.1.

 Next claim shows that the covering
 $(V_i)_{i\in\{0,\ldots,p\}}$ satisfies also property 4).
 
\proclaim{Claim 3.5} For every $ x\in C$ there is an open set $V_i$,
with $i\in\{0,\ldots,p\}$, such
that $x\in V_i$ and $d(x,\partial V_i)>r_i/3$
\endproclaim

\demo{Proof of Claim 3.5} Let $x\in C$, then by Lemma 3.3  $x\in
B(x_i,\frac23 r_i)$ for some $i\in\{0,\ldots,p\}$; we fix this index
$i$.   If $i=0$ or if $B(x_i,r_i)\cap W_0=\emptyset$, then
$V_i=B(x_i,r_i)$ and the lemma holds. Hence  we may assume that $i>0$ and
$B(x_i,r_i)\cap W_0\neq\emptyset$. Moreover, we can suppose
$d(x,x_0)>\frac23r_0$. In this case $V_i=B(x_i,r_i)-W_0$ and we
claim that $d(x,W_0)>\frac13r_i$. 

To prove this claim, we use the inequality:
$$
d(x,W_0)\geq d(x,x_0)-\Diam (W_0)\geq \tfrac23r_0-\tfrac19r_0=\tfrac59r_0,
$$
which is true because by assumption $d(x,x_0)>\frac23r_0$ and by Lemma 3.2 
 $\Diam(W_0)\leq\frac19r_0$.  Since $B(x_0,r_0)\cap
B(x_i,r_i)\neq\emptyset$, Lemma 3.1 implies that $r_0\geq \frac34 r_i$.
Hence $d(x,W_0)\geq \frac5{12}r_i>\frac13r_i$ and the claim is proved.
\qed
 \enddemo

Before proving property 5) of a $\eta$-covering \`a la Gromov, we point out that
property 6) of Proposition 2.6 follows from Lemma 3.1 and the fact that the balls $B(x_i,4
r_i)$ are abelian in $M-\Int( \widetilde W_0)$.
Moreover property 7) of Proposition 2.6 is satisfied by construction and Lemma 3.2.

Next claim proves property 5) of a $\eta$-covering \`a la Gromov and completes the proof of
Proposition 2.6.

\proclaim{Claim 3.6} There exists a universal constant $b_1>0$ such that  
$$
\Vol(V_i)\leq\Vol\big(B(x_i,r_i)\big)\leq\frac{b_1}Dr_i^3,
\qquad \forall i=0,\ldots,p
$$
\endproclaim

\demo{Proof of Claim 3.6} To estimate the volume of $B(x_i,r_i)$ we use
the same method as in Claim 4.7 of Chapter V.
To fix the notations, for $i=0,\ldots,p$, let
$f_i:U_{x_i}\to\NN_{\nu_i}(S_i)$ be the $(1+\varepsilon)$-bilipschitz
homeomorphism given by the Local Soul Theorem (Chapter IV).

 We need the following technical claim, whose proof is postponed  to the end of the section.

\proclaim{Claim 3.7} For $i=0,\ldots,p$, let $\nu_i$  denote the
radius of the normal fiber bundle 
 of the soul of the
Euclidean local model given by the Local Soul theorem. Then
$r_i>{\nu_i}/{2^{11}}.$
 \endproclaim

Assuming that Claim 3.7 is true, we can compare the volumes of the balls $B(x_i,r_i)$
and  $B(x_i,\nu_i/2^{11})$. Since $r_i>{\nu_i}/{2^{11}}$, by the Bishop-Gromov inequality
(Pro\-po\-sition
 1.6 of Chapter III)  we get:
$$
\Vol\big(B(x_i,r_i)\big)\leq
\Vol\big(B(x_i,\frac{\nu_i}{2^{11}})\big)
\frac{\V_{-1}(r_i)}{\V_{-1}(\nu_i/2^{11})},
$$
where $\V_{-1}(t)=\pi(\sinh(2t)-2t)$.

As in Claim 4.7 of Chapter V, let $a>0$ be a constant such that 
$t^3/a\leq{\V_{-1}(t)}\leq a\, t^3$ for  every
$t\in[0,1]$. Since $\nu_i\leq 1$ and $r_i\leq 1$, we get:
$$
\Vol\big(B(x_i,r_i)\big)\leq
\Vol\big(B(x_i,\frac{\nu_i}{2^{11}})\big)
a^2 2^{33}\frac{r_i^3}{\nu_i^3}.
$$

Since $d(f_i(x_i),S_i)\leq\nu_i/D<\nu_i\,{10^{-4}}$, we have that
 $f_i\big(B(x_i,\nu_i/2^{11})\big)\subset\NN_{\nu_i}(S_i)$. Thus
$$
\Vol\big(B(x_i,\frac{\nu_i}{2^{11}})\big)\leq
(1+\varepsilon)^3
\Vol(\NN_{\nu_i}(S_i))
\leq 2^3\Vol(\NN_{\nu_i}(S_i)),
 $$
because $f_i$ is $(1+\varepsilon)$-bilipschitz, with $\varepsilon\leq1/2$.

Using the bound $\Diam(S_i)\leq\nu_i/D$ and the fact that
the dimension of the soul $S_i$ is $1$ or $2$, we easily get the upper bound
$\Vol(\NN_{\nu_i}(S_i))\leq (\pi/D) \nu_i^3$, as in Claim 4.7 of Chapter V. Hence
$$ \Vol\big(B(x_i,r_i)\big)\leq 
\frac{b_1}D r_i^3,\quad \text{ with }
b_1=2^{36}a^2\pi. \qed
$$
 \enddemo

Finally the proof of Claim 3.7  concludes the proof of Proposition 2.6.

\demo{Proof of Claim 3.7} For $i=0,\ldots,p$, let
$f_i:U_{x_i}\to\NN_{\nu_i}(S_i)$ be the $(1+\varepsilon)$-bilipschitz
homeomorphism given by the Local Soul theorem (Chapter IV). We recall the upper bound
$$
\max(\Inj (x), d(f_i(x_i),S_i),\Diam(S_i))\leq\frac{\nu_i}D.
$$

If $i=0$, then it is clear that $\ab {x_0}\geq\nu_0/2$; thus
$r_0\geq\nu_0/16$, because $\nu_0\leq 1$. 

If $i\geq 1$, then
$\ab{x_i}\geq
\inf\big(\frac1{1+\varepsilon}\nu_i(1-\frac1D),d(x_i,W_0)\big)$,
because $W_0$ can intersect the neigbourhood $U_{x_i}$. Since
$\varepsilon<1/2$ and $D>10^4$, this inequality becomes 
$$
\ab{x_i}\geq\inf\big(\frac{\nu_i}2,d(x_i,W_0)\big).
$$
 Now we want to
find a  lower bound for $d(x_i,W_0)$.

Since  $d(x_i,x_0)>\frac{r_0}4$ by the choice of the sequence
$x_0,\ldots,x_p$ (property \ttag* above) and since $\Diam(W_0)\leq\frac{6\nu_0}D$ by the
proof of Lemma 3.2, first we  get the following lower bound:
$$
d(x_i,W_0)\geq d(x_i,x_0)-\Diam(W_0)>\frac{r_0}4-\frac{6\nu_0}D,
$$

Moreover, $d(x_i,W_0)>\nu_0(\frac1{64}-\frac6D)>\frac{\nu_0}{128}$, because $r_0\geq\nu_0/16$.
Therefore, since $\nu_0$ and $\nu_i\leq 1$, we obtain: 
$$
r_i\geq  \frac18\ab{ x_i} \geq\inf\Big(\frac{\nu_i}{2^4},\frac{\nu_0}{2^{10}}\Big).
$$
To compare $\nu_0$ and $\nu_i$ we distinguish two cases,
according to whether the local model for $U_{x_i}$ is of type a) or b).

If the local model for  $U_{x_i}$ is of type b), then by the choice of 
$x_0$, we have $\nu_0\geq\nu_i/2$, hence $r_i\geq\nu_i/2^{11}$.

When  the local model for $U_{x_i}$ is  of type a), again we
distinguish  two cases according to whether the intersection 
 $f_i^{-1}(\NN_{\nu_i/8}(S_i))\cap W_0$ is empty or not.

If $f_i^{-1}(\NN_{\nu_i/8}(S_i))\cap W_0=\emptyset$, then
$p^{-1}f_i^{-1}(\NN_{\nu_i/8}(S_i))$ is abelian in $M-\Int(\widetilde W_0)$ and
 we  have
$$
\ab{x_i}\geq d(x_i,\partial f_i^{-1}(\NN_{\nu_i/8}(S_i)))\geq
\frac1{1+\varepsilon}\big(\frac{\nu_i}8-d(f_i(x_i),S_i)\big)\geq
\frac1{1+\varepsilon}\big(\frac{\nu_i}8-\frac{\nu_i}D\big)>\frac{\nu_i}{16}  
$$
 and
we conclude that $r_i>\nu_i/128$. 

If  $f_i^{-1}(\NN_{\nu_i/8}(S_i))\cap W_0\neq\emptyset$, then there exists $y_0\in
W_0$ such that
$d(y_0,f_i^{-1}(S_i))\leq(1+\varepsilon){\nu_i}/8<{\nu_i}/4$.
Hence, 
$$
\text{for every }x\in W_0 , \quad d(x,f_i^{-1}(S_i))\leq 
    d(y_0,f_i^{-1}(S_i))+\Diam(W_0)\leq \frac{\nu_i}4+\frac{6\nu_0}D.
$$
Since $W_0$ corresponds to a $(\varepsilon,D)$-Margulis neighborhood of
type b), it cannot be contained in a
$(\varepsilon,D)$-Margulis neighborhood of type a). In particular,
$W_0$ cannot be contained in $f_i^{-1}(\NN_{\nu_i}(S_i))$ and we have:
$$
\frac{\nu_i}4+\frac{6\nu_0}D>\frac{\nu_i}{1+\varepsilon}>\frac{\nu_i}2.
$$
 We deduce that $\nu_0\geq {D}\, {\nu_i}/{24}>32\,\nu_i$, because
$D>10^4$. Thus 
$$
r_i\geq\inf\Big(\frac{\nu_i}{2^4},\frac{\nu_0}{2^{10}}\Big)  \geq\nu_i/32
$$ and the claim is proved.

This also concludes the proof of Proposition 2.6.
\qed

\enddemo

\enddemo
\enddemo

\newpage
\rightheadtext{ }
\leftheadtext{ }

\

\centerline{\chapt REFERENCES}

\

\

\widestnumber\key{1234567}

\ref
 \key{\BDV}
  \by  M. Boileau, S. Druck,  E. Vogt
   \paper Vanishing of Gromov volume and circle foliations of Epstein length 1 on open 3-manifolds
    \jour Preprint 1996
     \endref
\ref 
 \key{\BSOne} 
  \by F. Bonahon, L. C. Siebenmann 
   \paper The characteristic tori splitting of irreducible
     compact three-orbifolds \yr 1987 \vol 278
    \jour  Math. Ann. \pages 441-479
     \endref
\ref 
 \key{\BSTwo} 
  \by F. Bonahon, L. C. Siebenmann 
   \paper The classification of Seifert fibered three- manifolds   
    \ed R. Fenn   
     \inbook Low Dimensional Topology
      \bookinfo London Math. Soc. Lecture Notes Ser.  
       \publ Cambridge Univ. Press \publaddr Cambridge
        \pages 19-85 
         \vol 95 \yr 1985   
          \endref

\ref 
 \key{\BSThree} 
  \by F. Bonahon, L. C. Siebenmann 
   \book  Geometric splittings of classical knots and the algebraic knots of Conway
     \endref
\ref
 \key{\Bou}
  \by J.-P. Bourgignon
   \paper L'\'equation de la chaleur associ\'ee \`a la courbure
     de Ricci (d'apr\`es R. S. Hamilton)
     \inbook S\'eminaire Bourbaki 38\`eme ann\'ee, n. 653
      \bookinfo Ast\'e\-risque
       \vol 145-146
        \yr 1987
         \endref
\ref
 \key{\BrS}
  \by M. R. Bridson, G. A. Swarup
   \paper On Hausdorff-Gromov convergence and a theorem of Paulin
    \jour L'Ens. Math. \vol 40 \yr 1994
     \pages  267-289
      \endref
\ref
 \key{\CEG}
  \by R. D. Canary, D. B. A. Epstein and P. Green
   \paper Notes on notes of Thurston
    \yr 1987 \vol 111
     \ed D. B. A. Epstein
      \inbook Analytical and Geometric Aspects of Hyperbolic Space
       \bookinfo London Math. Soc. Lecture Notes Ser.  
        \publ Cambridge Univ. Press \publaddr Cambridge
         \pages 3-92
          \endref
\ref
 \key{\CG}
  \by J. Cheeger, M. Gromov
   \paper Collapsing Riemannian manifolds while keeping their
     curvature bounded
     \jour J. Diff. Geom. \yr 1990 \vol 32 \pages 269-298
      \endref
\ref
 \key{\CGl}
  \by J. Cheeger, D. Gromoll
   \paper On the structure of complete manifolds of nonnegative
     curvature
     \jour Ann. of Math. \yr 1986 \vol 72 \pages 413-443
      \endref
\ref 
 \key{\CHK}
  \by D. Cooper, C. Hodgson, S. Kerchkoff
   \book
    \bookinfo In preparation
     \yr 1998
      \endref
\ref 
 \key{\Cul} 
  \by M. Culler 
   \paper Lifting representations to covering groups \yr 1986 \vol 59 
    \jour Adv. Math.\pages 64-70 
     \endref
\ref 
 \key{\CS} 
  \by M. Culler, P. B. Shalen 
   \paper Varieties of group representations and splittings of
    3-manifolds 
     \yr 1983 \vol 117  \jour Ann. of Math.\pages 109-146
      \endref
\ref
 \key{\DaM}
  \by M. W. Davis, J. W. Morgan
   \paper Finite group actions on homotopy 3-spheres
    \pages 181-225
     \ed  H. Bass, J. W. Morgan
      \inbook  The Smith Conjecture
       \publ Academic Press \yr 1984
     \endref
\ref
 \key{\DD} \by W\. Dicks, M\. J\. Dunwoody 
  \book Groups acting on graphs
   \publ Cambridge University Press
    \publaddr Cambridge, 1989
     \endref
\ref 
 \key{\Dun} 
  \by W. D. Dunbar 
   \paper Geometric Orbifolds 
     \yr 1988 \vol 1  \jour Rev. Mat. Univ. Comp. Madrid \pages 67-99
      \endref
\ref 
 \key{\DuM} 
  \by W. D. Dunbar, R. G. Meyerhoff
   \paper Volumes of Hyperbolic 3-Orbifolds 
     \yr 1994 \vol 43  \jour Indiana Univ. Math. J. \pages 611-637
      \endref
\ref 
 \key{\Ebe} 
  \by P\. Eberlein 
   \paper Lattices in spaces of nonpositive curvature
     \yr 1980 \vol 111  \jour Ann. of Math.\pages 435-476
      \endref
\ref
 \key{\Fed}
  \by H\. Federer
   \book Geometric Measure Theory
    \yr 1969  \publ Springer Verlag \publaddr Berlin
     \endref
\ref
 \key{\GM}
  \by Gonzalez-Acu\~na, J\. M\. Montesinos 
   \paper On the character variety of group representations in $SL(2,\Bbb C)$
     et $PSL(2,\Bbb C)$ \jour Math. Z. \vol 214 \yr 1983 \pages 627-652
    \endref
\ref
 \key{\Gro}
  \by M\. Gromov
   \paper Volume and bounded cohomology
    \yr 1983 \vol 56 \jour I. H. E. S. Publ. Math.
     \pages 5-99
      \endref 
\ref
 \key{GLP}
  \by M\. Gromov, J\. Lafontaine, P\. Pansu
   \book Structures m\'etriques pour les vari\'et\'es riemanniennes
    \yr 1981 \publ Cedic/Fernand Nathan \publaddr Paris
     \endref
\ref
 \key{\GT}
  \by M. Gromov, W. Thurston
   \paper  Pinching constants for hyperbolic manifolds
    \yr 1987 \vol 89 \jour Invent. Math.
     \pages 1-12
      \endref 
\ref
 \key{\HamOne}
  \by R. S. Hamilton 
   \paper Three-manifolds with positive Ricci curvature
    \yr 1982 \vol 17 \jour J. Diff. Geom.
     \pages 255-306
      \endref 
\ref
 \key{\HamTwo}
  \by R. S. Hamilton 
   \paper Four-manifolds with positive curvature operator
    \yr 1986 \vol 24 \jour J. Diff. Geom.
     \pages 153-179
      \endref 
\ref
 \key{\Hod}
  \by C. Hodgson
   \paper Geometric structures on  3-dimensional orbifolds: Notes on
     Thurston's proof
     \jour Preprint
      \endref
\ref 
 \key{\HodTwo}
  \by C. Hodgson 
   \book Degeneration and Regeneration of Hyperbolic Structures on
     Three-Manifolds 
     \bookinfo Thesis, Princeton University  \yr 1986
      \endref
\ref 
 \key{\HK}
  \by C. Hodgson, S. Kerchkoff
   \paper A rigidity theorem for hyperbolic cone manifolds 
    \yr 1993\jour Preprint
     \endref
\ref 
 \key{\HT}
  \by C. Hodgson, J. Tysk
   \paper Eigenvalue estimates and isoperimetric inequalities
     for cone manifolds
     \jour  Bull. Austral. Math. Soc.
      \vol 47 \yr 1993 \pages 127-144
       \endref 
\ref 
 \key{\Iva}
  \by N. V. Ivanov
   \paper Foundations of the theory of bounded cohomology
    \jour J. Sov. Math.  \vol 37 \yr 1987 \pages 1090-1115
     \endref
\ref 
 \key{\JR}
  \by W. H. Jaco, J.H. Rubinstein
   \paper P-L minimal surfaces in 3-manifolds
    \jour J. Diff. Geom. \vol 27 \yr 1988 \pages 493-524 
     \endref
\ref 
 \key{\JS}
  \by W. H. Jaco,  P. B. Shalen
   \paper Seifert fibred spaces in 3-manifolds
    \jour Mem. Amer. Math. Soc.  \vol 220 \yr 1979
     \endref
\ref
 \key{\Joh}
  \by K. Johannson
   \book Homotopy equivalences of 3-manifolds with boundary
    \bookinfo  Lecture Notes in Mathematics \vol 761
     \publ Springer Verlag \publaddr Berlin \yr 1979
      \endref
\ref 
 \key{\Jon}
  \by K. Jonnes
   \endref
\ref
 \key{\Kap}
  \by M. Kapovich
   \book Hyperbolic Manifolds and Discrete Groups: Notes on 
    Thurs\-ton's Hyperbolization
    \bookinfo University of Utah lecture Notes, 1993-1994
     \endref
\ref
 \key{\Kir}
  \by R. Kirby
   \paper Problems in Low-Dimensional Topology
    \inbook Geometric Topology \vol 2 \yr 1997   \ed W.H. Kazez
     \bookinfo Amer. Math. Soc. Internat. Press
      \endref
\ref
 \key{\KoOne}
  \by S. Kojima  
   \paper Deformations of hyperbolic 3-cone-manifolds
    \jour Preprint \yr 1997
     \endref
\ref
 \key{\KoTwo}
  \by S. Kojima
   \paper Finiteness of Symmetries on 3-manifolds
    \jour Preprint
     \endref
\ref 
 \key{\McMOne}
  \by C. McMullen
   \paper Iteration on Teichm\"uller space
    \jour  Invent. Math. \vol 99 \yr 1990 \pages 425-454 
     \endref
\ref 
 \key{\McMTwo}
  \by C. McMullen
   \book Renormalization and 3-manifolds which Fiber over the Circle 
    \bookinfo Annals of Math. Studies
     \publ Princeton University Press
      \vol 144 \yr 1996
       \endref
\ref 
 \key{\MS}
  \by W. H. Meeks, P. Scott
   \paper  Finite group actions on 3-manifolds
    \jour   Invent. math. 
     \vol 86  \yr 1986 \pages 287-346
      \endref
\ref 
 \key{\MYOne}
  \by W. H. Meeks, S.T. Yau
   \paper Topology of three-dimensional manifolds and the embedding problems
    in minimal surface theory
     \jour  Ann. of Math.
      \vol 112  \yr 1980 \pages 441-484
       \endref 

\ref 
 \key{\MYTwo}
  \by W. H. Meeks, S.T. Yau
   \paper  The equivariant Dehn's Lemma and Loop Theorem
     \jour  Comment. Math. Helv.
      \vol 56  \yr 1981 \pages 225-239
       \endref 
\ref 
 \key{\MB}
  \by H. Bass, J. W. Morgan
   \book The Smith Conjecture
    \publ Academic Press \yr 1984
     \endref
\ref
 \key{\Mos}
  \by G\. D\. Mostow
   \book Strong rigidity of locally symmetric spaces 
    \bookinfo Ann. of Math. Stud. \vol 78
     \publ  Princeton University Press \publaddr Princeton, NJ \yr 1973
      \endref
\ref 
 \key{\Mum}
  \by D. Mumford
   \book Algebraic Geometry I: Complex projective Varieties
    \publ Springer Verlag \yr 1976 \publaddr Berlin
     \endref
\ref 
 \key{\OtaOne}
  \by J.-P. Otal
   \paper Le th\'eor\`eme d'hyperbolisation pour les
     vari\'et\'es fibr\'ees de dimension 3
     \jour  Ast\'erisque
      \vol 235 \yr 1996
       \endref 
\ref 
 \key{\OtaTwo}
  \by J.-P. Otal
   \paper Thurston's hyperbolization of Haken manifolds
    \jour  Preprint
     \yr 1997
      \endref 
\ref
 \key{\Pe}
  \by S. Peters
   \paper Convergence of Riemannian manifolds
    \jour Compositio. Math.  \vol 62
     \yr 1987
      \pages 3-16 
       \endref
\ref
 \key{\PoOne}
  \by J. Porti
   \paper Torsion de Reidemeister pour les Vari\'et\'es Hyperboliques
    \jour Mem. Amer. Math. Soc. \vol 128 \yr 1997
     \endref
\ref
 \key{\PoTwo}
  \by J. Porti
   \paper Regenerating hyperbolic and spherical
     cone structures from Euclidean ones
     \jour Topology  \vol 37 \pages 365-392 \yr 1998
      \endref
\ref
 \key{\Rag}
  \by M. S. Raghunatan
   \book Discrete Subgroups of Lie Groups
    \publ Springer Verlag
     \publaddr Berlin \yr1972
      \endref
\ref
 \key{\Sak}
  \by T. Sakai
   \book  Riemannian Geometry
    \publ Amer. Math. Soc.
     \publaddr Rodhe Island \yr 1996
      \endref
\ref
 \key{\SOK}
  \by T. Soma, K. Ohshika, S. Kojima
   \paper \, Towards a proof of Thurston's geometrization theorem
     for orbifolds 
    \jour RIMS. K\^oky\^uroku \vol 568 \yr 1985 \pages 1-72
     \endref

\ref
 \key{\Som}
  \by T. Soma
   \paper \, The Gromov invariant of links 
    \jour \, Invent. Math. 
     \yr \, 1981  \vol 64  \pages 445-454
      \endref
\ref
 \key{\Scott}
  \by P. Scott
   \paper The geometries of 3-manifolds
    \jour Bull. London Math. Soc.
     \yr 1983 \pages 401-487 \vol 15
      \endref
\ref
 \key{\Sua}
  \by E. Suarez
   \book Poliedros de Dirichlet de 3-variedades c\'onicas
     y sus deformaciones
     \bookinfo Thesis Univ. Compl. Madrid \yr 1998
      \endref
\ref
 \key{\ToOne}
  \by J. L. Tollefson
   \paper Involutions of sufficiently large 3-manifolds
    \jour Topology
     \yr 1981 \pages 323-352 \vol 20
      \endref
\ref
 \key{\ToTwo}
  \by J. L. Tollefson
   \paper Periodic homeomorphisms of 3-manifolds fibered over $S^1$
    \jour Trans. Amer. Math. Soc.
     \yr 1976 \pages 223-234
\moreref\paper \rm see also erratum \jour Trans. Amer. Math. Soc.
     \yr 1978 \pages 309-310
      \endref
\ref
 \key{\ThuNotes}
  \by W. P. Thurston
   \book The geometry and topology of 3-manifolds
    \publ Princeton Math. Dept. \yr 1979
     \endref
\ref
 \key{\ThuBull}
  \by W. P. Thurston
   \paper Three dimensional manifolds, Kleinian groups and hyperbolic
     geometry
     \vol 6 \jour Bull. Amer. Math. Soc. \yr 1982 \pages 357-381
      \endref
\ref
 \key{\Thu3}
  \by W. P. Thurston
   \paper Hyperbolic Structures on 3-manifolds, I: 
     Deformations of acylindrical manifolds
     \vol 124 \jour Annals of Math. \yr 1986 \pages 203-246
      \endref
\ref
 \key{\Thu4}
  \by W. P. Thurston
   \paper Hyperbolic Structures on 3-manifolds, II: 
     Surface groups and manifolds which fiber over $S^1$ 
    \yr 1986\jour Preprint
     \endref
\ref
 \key{\Thu5}
  \by  W. P. Thurston
   \paper Hyperbolic Structures on 3-manifolds, III: 
     Deformations of 3-man\-ifolds with incompressible boundary
    \yr 1986\jour Preprint
     \endref
\ref
 \key{\Thu6}
  \by  W. P. Thurston
   \paper Three-manifolds with symmetry
    \yr 1982\jour Preprint
     \endref
\ref
 \key{\Yam}
  \by T. Yamaguchi
   \paper Collapsing and pinching under a lower curvature bound
    \yr 1991 \vol 133 \jour Ann. of Math. \pages 317-357
     \endref 
\ref
 \key{\WalOne}
  \by F. Waldhausen
   \paper Eine Klasse von 3-dimensionalen Mannigfaltigkeiten I and II
    \yr 1967 \vol 3 {\rm and} 4 \jour  Invent. Math. \pages 308-333 and 87-117
     \endref 
\ref
 \key{\WalTwo}
  \by F. Waldhausen
   \paper On irreducible 3-manifolds which are suffiently large
    \yr 1968 \vol 87 \jour Ann. of Math. \pages 56-88
     \endref 
\ref
 \key{\Wei}
  \by A. Weil
   \paper Remarks on the cohomology of groups
    \yr 1964 \vol 80 \jour Ann. of Math. \pages 149-157
     \endref 
\ref
 \key{\Zhou}
  \by Q. Zhou
   \book 3-dimensional Geometric Cone Structures 
    \bookinfo Thesis, University of California, L.A. \yr 1989
     \endref
\ref
 \key{\ZhouTwo}
  \by Q. Zhou
   \book the Moduli Space or Hyperbolic Cone Structures 
    \bookinfo Preprint \yr 1997
     \endref

\

\smc\noindent Laboratoire Emile Picard, CNRS UMR 5580,
Universit\'e Paul Sabatier, 118 Route de Narbonne,
  F-31062 TOULOUSE Cedex 4, France

\

 \it e.mail:  \rm
boileau\@picard.ups-tlse.fr, \rm
porti\@picard.ups-tlse.fr
\enddocument